\documentclass[12pt]{amsbook}

\usepackage{pb-diagram}

\usepackage{amsmath}

\usepackage{amsthm}

\usepackage{amsfonts}

\usepackage{amssymb}

\newcommand\R{\mathbb{R}}

\newcommand\Q{\mathbb{Q}}
\newcommand\C{\mathbb{C}}
\newcommand\Z{\mathbb{Z}}
\newcommand\Ha{\mathbb{H}}

\newcommand\kap{\kappa_0 \cup_{\theta} \kappa_1^{-}}

\newcommand\GE{\bar{G}^1_2}
\newcommand\hl{\hat{\lambda}}
\newcommand\bl{\bar{\lambda}}

\newcommand\TDG{TG^{\wedge}}
\newcommand\ph{{\rm ph}}
\newcommand\Sec{{\rm Sec}}

\newcommand\Hom{{\rm Hom}}
\newcommand\Ker{{\rm Ker}}
\newcommand\Cok{{\rm Cok}}

\newcommand\expo{{\rm exp}}
\newcommand\modu{{\rm mod}}
\newcommand\Top{\rm {Top}}
\newcommand\Cat{\rm {Cat}}

\newcommand\mcal[1]{\mathcal{#1}}

\newcommand\LL{\mathcal{L}}

\newcommand\HC{\mathcal{HC}}
\newcommand\HH{\mathcal{H}}
\newcommand\QL{\mathcal{Q}}
\newcommand\QF{\mathcal{F}}
\newcommand\FF{\mathcal{F}}

\newcommand\NN{\mathcal{N}}
\newcommand\NCat{\mathcal{N}^{\Cat}}

\newcommand\SO{\mathcal{S}^O}

\newcommand\SPL{\mathcal{S}^{PL}}
\newcommand\SCat{\mathcal{S}^{\Cat}}

\newcommand\NPL{\mathcal{N}^{PL}}

\newcommand\M{\left( \begin{array}{cc} \bl_G|_{G(\Psi)} & 0 \\ 0 & \hat{b}({\hl}) \end{array} \right) }
\newcommand\tensor{\otimes}

\newcommand\del{\partial}
\newcommand\ahat{\hat{A}}

\newcommand{\lra}{\longrightarrow}

\newcommand\hra{\hookrightarrow}

\newcommand\ra{\rightarrow}

\newcommand\nd{\not  | \;}
\newcommand\ba{\bar{\alpha}}
\newcommand\bg{\bar{\gamma}}

\newcommand\halfp{\frac{p_1}{2}}

\newcommand\bs{\bar{s}_1}

\newtheorem{Theorem}{Theorem}[chapter]

\newtheorem{Definition}[Theorem]{Definition}

\newtheorem{Example}[Theorem]{Example}

\newtheorem{Lemma}[Theorem]{Lemma}

\newtheorem{Proposition}[Theorem]{Proposition}

\newtheorem{Corollary}[Theorem]{Corollary}

\newtheorem{Remark}[Theorem]{Remark}
\newenvironment{remark}{\begin{Remark}\rm}{\end{Remark}}

\newtheorem{Fact}[Theorem]{Fact}

\newtheorem{Diagram}[Theorem]{Diagram}

\newtheorem{theorem}{Theorem}

\newtheorem{thm}{Theorem}

\title{The classification of highly connected manifolds in
dimensions 7 and 15.}
\author{Diarmuid J. Crowley\\
\vskip 2cm
\small{
Submitted to the faculty of the University Graduate School\\
in partial fulfillment of the requirements\\
for the degree\\
Doctor of Philosophy \\
in the Department of Mathematics\\
Indiana University\\
July 23rd 2001}
}
\makeindex

\begin{document}

\maketitle
\tableofcontents
\mainmatter

\chapter{Introduction} \label{chapintro}
An $n$-dimensional manifold $M$ is called highly connected when all of its homotopy
groups vanish below the middle dimension.  That is, $\pi_i(M) = 0$ for all $i < (n-1)/2$.  We begin 
by briefly recalling the history of the classification of highly connected, closed, smooth manifolds in
dimensions greater than $4$.  Throughout all manifolds will be smooth and oriented and all maps
orientation preserving.  In \cite{Wa1} Wall built on the results of Smale, most notably the
$h$-cobordism theorem, to classify highly connected even dimensional manifolds up to the addition of
homotopy spheres.  Later Wall \cite{Wa3} presented a classification of highly connected odd dimensional
manifolds up to the addition of homotopy spheres except in dimensions $5$, $7$ and $15$.  Barden,
a student of Wall, solved the problem of the classification of simply connected $5$-manifolds
\cite{Ba}.  In dimensions $7$ and $15$ Wilkens, another student of Wall, obtained results
but they were not complete as we now explain \cite{Wi1,Wi2}.

Throughout this thesis $P$ shall denote a highly connected, closed manifold of dimension $(8j-1)$ and $L$ shall
denote a handlebody by which we mean a highly connected, compact, even dimensional manifold with a nonempty
and connected boundary.  For every $P$ there is a unique homotopy sphere such $\Sigma_P$ (often the standard
sphere) such that $P \# \Sigma_P$ is the boundary of a handlebody $L$ and such that $\Sigma_P$ bounds a Spin
manifold with vanishing decomposable Pontryagin numbers and signature equal to zero (see Section
\ref{chapprelim}.\ref{sechcb}).  To each handlebody $L$, Wall associated the triple
$(H,\lambda,\alpha)$, where $H=H_{4j}(L;\Z)$ is the middle homology group of $L$, $\lambda: H \times H \ra \Z$
is the intersection form of
$L$ and $\alpha \in H^{4j}(L;\Z)$ is the primary obstruction to the stable triviality of $L$.  We identify
$H^{4j}(L;\Z)$ with
$H^*$, the dual group of $H$, via the Kronecker pairing.  Two triples $(H_0,\lambda_0,\alpha_0)$ and
$(H_1,\lambda_1,\alpha_1)$ are isomorphic if there is an isomorphism $\Theta: H_0 \ra H_1$ such that 
$\Theta^*(\alpha_1) = \alpha_0$ and $\lambda_0 = \lambda_1 \circ (\Theta \times \Theta)$.  If the handlebodies
$L_0$ and $L_1$ have associated triples $(H_0,\lambda_0,\alpha_0)$ and $(H_1,\lambda_1,\alpha_1)$ then Wall
proved in \cite{Wa1} that an isomorphism between these triples is realized by a diffeomorphism.
Wilkens later proved that if $f:\del L_0 \ra \del L_1$ is a diffeomorphism between the boundaries of two
handlebodies, then there are handlebodies $V_0$ and $V_1$ with standard spheres as boundary and a
diffeomorphism $g: L_0 \natural V_0 \ra L_1 \natural V_1$ which extends $f$ (here $\natural$ denotes the
boundary connected sum).   Now, the intersection form of a handlebody is nonsingular if and only if the
handlebody cobounds a homotopy sphere.   Thus the question of classifying the manifolds $P$ up to connected
with sum homotopy spheres becomes the question of the stable classification of triples $(H,\lambda,\alpha)$. 
Here we say that two triples $(H_0,\lambda_0,\alpha_0)$ and $(H_1,\lambda_1,\alpha_1)$ are stably equivalent 
if there are triples $(H_0',\lambda_0',\alpha_0')$ and $(H_1',\lambda_1',\alpha_1')$ such that $\lambda_0'$
and $\lambda_1'$ are nonsingular and $(H_0 \oplus H_0',\lambda_0 \oplus \lambda_0',\alpha_0 \oplus \alpha_0')$
is isomorphic to $(H_1 \oplus H_1',\lambda_1 \oplus \lambda_1',\alpha_1 \oplus \alpha_1')$.

To classify the manifolds $P$ up to almost diffeomorphism Wilkens associated to each $P$ the triple
$(G,b,\beta)$ where $G=H^{4j}(P;\Z)$, $b$ is the linking form of $P$ which is a nonsingular symmetric bilinear
function $b: TG \times TG \ra \Q/\Z$ where $TG$ is the torsion subgroup of $G$ and $\beta \in G$ is the primary
obstruction to the triviality of the stable tangent bundle of $P$.  We shall call this triple the Wilkens
triple or Wilkens invariants of $P$.  Two triples
$(G_0,b_0,\beta_0)$ and $(G_1,b_1,\beta_1)$ are isomorphic if there is an isomorphism $\theta: G_0 \ra
G_1$ such that $\theta(\beta_0) = \beta_1$ and $b_0 = b_1 \circ (\theta \times \theta)$.  The only
constraint on the triples which are realized by highly connected $7$ and $15$ manifolds is that $\beta =
2\gamma$ for some $\gamma \in G$.  In higher dimensions all triples are realized.  If $P = \del L$
and $(H,\lambda,\alpha)$ is the triple associated to $L$ then $(H,\lambda,\alpha)$ induces $(G,b,\beta)$
in a natural way.  Wilkens was able to prove that if the group $G$ contains no $2$-torsion then any two
highly connected $(8j-1)$-manifolds, $j=1,2$, with isomorphic triples are diffeomorphic after the addition of
homotopy spheres.  When $G$ contains $2$ torsion and $b$ is not decomposable as the sum of
smaller nontrivial linking forms, Wilkens showed that up to connected sum with homotopy spheres, there were at
most two highly connected $(8j-1)$-manifolds, $j=1,2$ with isomorphic triples.  For a precise statement of
Wilkens' results we refer the reader to Section \ref{chapprelim}.\ref{secwilkc}.  A $\Z/2$-ambiguity
associated to groups of even order points to the need for a quadratic refinement of the linking form $b$. 
Indeed Wall had left the classification unsolved in dimensions $7$ and $15$ because the existence of bundles
of Hopf-invariant one over the
$4$-sphere and $8$-sphere meant that his definition of the quadratic refinement of the linking form
of a highly connected manifold did not work.  Wilkens, however, did not interpret his results on the
stable algebra of triples $(H,\lambda,\alpha)$ in terms of quadratic refinements of the induced linking
forms.

Now let $G$ be a finite abelian group.  A quadratic linking function on $G$ is a function 
$q:G \ra \Q/\Z$ such that 
$$q(x+y) = q(x) + q(y) + b(x,y)$$
for every $x, y \in G$ and for some linking form, $b$, on $G$.  If $q$ is homogeneous, that is 
$q(-x) = q(x)$ for every $x \in G$, then $q$ is called a quadratic linking form.  We let $\QL(b)$ denote
the set of quadratic linking functions refining a linking form $b$.   There is a free
action of $G$ on $\QL(b)$ defined, for $a \in G$ and $q \in \QL(b)$, by $a\cdot q = q_a$ where for every $x
\in G$,
$$q_a(x) = q(x) + b(x,a).$$  

A nondegenerate, even bilinear form $\lambda: H \times H \ra \Z$ on the free abelian group $H$ induces a
quadratic linking form, $q(\lambda)$, on $G = \Cok(\hl)$ where $\hl: H \ra H^*$ is the adjoint of
$\lambda$ (see Definition \ref{defqlf}).  Both Durfee \cite{Du1,Du2} and
Wall \cite{Wa4} proved that if $\lambda_0$ and $\lambda_1$ are nondegenerate even bilinear forms then
$\lambda_0$ and $\lambda_1$ are stably equivalent if and only if $q(\lambda_0) \cong q(\lambda_1)$, where
stable equivalence means isomorphic after the addition of nonsingular even bilinear forms.  However,
they did not consider the question of the stable equivalence of the triples
$(H,\lambda,\alpha)$.  In particular in dimensions $8$ and $16$ the triples $(H,\lambda,\alpha)$ arising
from handlebodies need not contain even intersection forms.  So these algebraic results were confined to
topological applications in dimensions other than $7$ and $15$.

Thus the question of the classification of $2$-connected $7$-manifolds and $6$-connected $15$-manifolds
stood where Wilkens left it for thirty years.\footnote{This is not quite correct.  Madsen, Taylor and
Williams \cite{MTW} claimed that Wilkens' ambiguous pairs are homeomorphic but this is not so.  However more
on this in Section \ref{chaphom}.\ref{secsphere}.} Then Grove and Ziller showed that the total spaces of
$S^3$-bundles over
$S^4$ admit metrics of non-negative sectional curvature and they called for a classification of these total
spaces
\cite[Problem 5.2]{GZ}.  Of course such total spaces are $2$-connected $7$-manifolds.  This problem was
solved in \cite{CE} where we showed that the invariants $s_1$ and $\bs$ defined in \cite{KS}
\footnote{We explain the invariants $s_1$ and $\bs$ in detail in Section \ref{chapprelim}.\ref{secs_1}.  They
are defined in all dimensions $4k-1$ and in dimensions $7$ and $15$, $s_1$ coincides with the Eels-Kuiper
$\mu$-invariant \cite{EK}.} 
resolved the ambiguity in Wilkens' classification.  Specifically we showed that if $P_0$ and
$P_1$ are the total spaces of $S^3$-bundles over $S^4$ which do not admit a cross section then $P_0$ is
homeomorphic (resp. diffeomorphic) to $P_1$ if and only if their Wilkens invariants are isomorphic and
$\bs(P_0) =
\bs(P_1)$ (resp. $s_1(P_0) =  s_1(P_1)$).  We shall prove that the same statement holds for all
$2$-connected $7$-dimensional rational homology spheres.

Now let $L$ be a handlebody of dimension $8$ or $16$ with associated triple $(H,\lambda,\alpha)$.  From the
point of view of Durfee's and Wall's results above the difficulty in dimensions $8$ and
$16$ is that $\lambda$ need not be an even form which is an immediate consequence of the existence of bundles
of Hopf invariant one.  However, this is compensated by the fact that the existence of bundles of Hopf
invariant one ensures that the tangential invariant, $\alpha$, of $L$ reduces $\modu~2$ to either
$w_4$ or $w_8$, the fourth or eight Stiefel-Whitney class.  Thus $\alpha$ is characteristic for $\lambda$ by
which we mean that  the quadratic function $\kappa(\lambda,\alpha)$ defined by
$$\kappa(\lambda,\alpha)(v) = \lambda(v,v) + \alpha(v) $$
takes on even values for every $v \in H$.  When $\lambda$ is nondegenerate this fact allows the definition
of a  quadratic linking function on the cokernel of $\hl$ as follows.  For every $x \in H^*$ there is a
nonzero integer $r$ and an element $v \in H$ such that $\hl(v) = rx$.  If $[x]$ denotes the image of $x$ in
$\Cok(\hl)$ then 
$$q^c(\kappa(\lambda,\alpha))([x]) := \frac{1}{2r}(x(v) + \alpha(v)) ~\modu~\Z.$$
We note that $q^c(\kappa(\lambda,\alpha))$ is not necessarily homogeneous.   Let $P$ be a highly
connected $(8j-1)$-manifold, $j=1,2$.  Recall that there is a homotopy sphere
$\Sigma_P$ such that $P \# \Sigma_P = \del L$ where $L$ is a handlebody.  Suppose that $L$ defines the
triple $(H,\lambda,\alpha)$.   Then, identifying $H^{4j}(P)$ and  $H^{4j}(P \# \Sigma_P)$ if necessary,
$\Cok(\hl) = H^{4j}(P)$.  If $P$ is a rational homology sphere then $H^{4j}(P)$ is finite and $\lambda$
is nondegenerate so we may define the quadratic linking function of $P$ by
$$q^c(P) := q^c(\kappa(\lambda,\alpha)).$$
We shall show that $q^c(P)$ is independent of the choice of handlebody and refines $b$, the linking form
of $P$.  Also we show that for any quadratic linking form $q$ which refines $b$ there is an element $a \in G$
such that $q^c(P) = q_a$ and that the Wilkens invariant of $P$ is $(G,b,2a)$.

As we first classify the manifolds $P$ up connected sum with homotopy spheres we introduce the notion of
an {\em almost diffeomorphism}.  A homeomorphism $f:M_0 \ra M_1$ between smooth manifolds is an almost
diffeomorphism if it is a diffeomorphism except at a finite number of singular points.  Two manifolds
are called almost diffeomorphic when there is an almost diffeomorphism between them.  Note that $M_0$
and $M_1$ are almost diffeomorphic if and only if there is a single homotopy sphere $\Sigma$ and a
diffeomorphism $f':M_0 \# \Sigma \ra M_1$.  Given any degree one map $f:P_0 \ra P_1$ between two highly
connected $(8j-1)$-manifolds we define the induced Umker map on cohomology, $f_{!}:H^{4j}(P_0) \ra
H^{4j}(P_1)$, where for $x \in H^{4j}(P_0)$, $f_{!}(x) := f_*(x \cap [P_0])/(P_1)$.  Here 
$[P_0] \in H_{8j}(P_0)$ and $(P_1) \in H^{8j}(P_1)$ are the orientation classes of
$P_0$ and $P_1$ and $\cap$ and $/$ denote the cap and slant products respectively.  

\begin{thm} \label{thma}
Let $j=1$ or $2$.  For every finite group $G$, every linking form $b:G \times G \ra \Q/\Z$ and every
quadratic linking function $q \in \QL(b)$, there is a highly connected $(8j-1)$-dimensional rational
homology sphere, $P$, such that $q^c(P) \cong q$.  If $P_0$ and $P_1$ are two highly connected
$(8j-1)$-dimensional rational homology spheres with quadratic linking functions $q^c(P_0)$ and $q^c(P_1)$ 
and Wilkens invariants $(G_0,b_0,\beta_0)$ and
$(G_1,b_1,\beta_1)$, then,
\begin{enumerate}
\item{there is an almost diffeomorphism  $f:P_0 \ra P_1$ such that 
$f_! = \theta: H^{4j}(P_0) \cong H^{4j}(P_1)$ if and only if $q^c(P_0) = q^c(P_1) \circ \theta$,}
\item{there is a diffeomorphism  $f:P_0 \ra P_1$ such that $f_! = \theta: H^{4j}(P_0) \cong H^{4j}(P_1)$
if and only if $q^c(P_0) = q^c(P_1) \circ \theta$, $s_1(P_0) = s_1(P_1)$ and $\Sigma_{P_0} \cong
\Sigma_{P_1}$,}
\item{$P_0$ and $P_1$ are almost diffeomorphic if and only if $(G_0,b_0,\beta_0) \cong
(G_1,b_1,\beta_1)$ and $\bs(P_0) = \bs(P_1)$,}
\item{$P_0$ and $P_1$ are diffeomorphic if and only if $(G_0,b_0,\beta_0) \cong
(G_1,b_1,\beta_1)$, $s_1(P_0) = s_1(P_1)$ and $\Sigma_{P_0} \cong \Sigma_{P_1}$,}
\item{in dimension $7$, 	``almost diffeomorphism'' may be replaced by ``homeomorphism'' in {\em 1},}
\item{in dimension $7$ there is a homotopy equivalence $f:P_0 \ra P_1$ such that 
$f_! = \theta: H^{4j}(P_0) \cong H^{4j}(P_1)$ if and only if $q^c(P_0)_{12a} = q^c(P_1) \circ \theta$
for some $a \in G_0$.}
\end{enumerate}
\end{thm}

\begin{remark}
The proof of Theorem \ref{thma} is distributed throughout the thesis as follows.  Part {\em 1} is a special
case of Theorem \ref{thmb} which is proven in Theorem \ref{thmB}.  Part {\em 2} is a special case of Theorem
\ref{thmsmoothclass}.  Parts {\em 3} and {\em 4} are equivalent to Corollary \ref{corwilkdata}.  Part {\em 5}
is proven in Theorem \ref{thm7homeo} and Part {\em 6} is a special case of Theorem \ref{thmhomotpclass}.
\end{remark}

\begin{remark}
Note that we do not follow the classical surgery path of first classifying the manifolds $P$ up to homotopy and
then up to homeomorphism.  Nor do we use Kreck's \cite{Kr2} modified surgery.  Rather, we use Wall's \cite{Wa1}
classification on the handlebodies $L$ and consider more closely the algebraic structures particular to
dimensions $7/8$ and $15/16$ which are induced on the boundaries $P$.\footnote{This is the route
suggested by Wall just after Theorem 8 in \cite{Wa3} (and followed by Wilkens).}  In Chapter
\ref{chaphom} we use surgery theory and the topological classification of part {\em 5} to deduce the homotopy
classification.
\end{remark}

\begin{remark}
Part {\em 4} of this theorem extends Theorem 6 of \cite{Kr2} to cases where $H^4(P_i)$ is not generated by
$\beta_i$.  After completing this thesis we discovered that F. Bermbach had already undertaken such an
extension in \cite{Be}.  However, Bermbach's general theorem on the classification of simply-connected
$7$-manifolds uses a Wall style homogeneous quadratic refinement which is not easy to calculate.  Moreover, he
does not apply this theorem to the classification of $2$-connected $7$-manifolds but in applications
he is exclusively concerned with certain homogeneous spaces having $\pi_2 \cong \Z^2$.
\end{remark}

We now proceed to the general case where $G=H^{4j}(P)$ might have a nontrivial free subgroup.  If $G$ is
itself free then the classification problem is easily solved by Wall's classification of handlebodies.

\begin{Proposition} \label{propfree}
Let $P_0$ and $P_1$ be highly connected $(8j-1)$-manifolds with Wilkens invariants $(G_0,b_0,\beta_0)$
and $(G_1,b_1,\beta_1)$ and suppose that $G_0$ and $G_1$ are free groups.  Then there is an
almost diffeomorphism $f:P_0 \ra P_1$ such that $f_!=\theta:G_0 \cong G_1$ if and only if
$\theta(\beta_0) = \beta_1$.  Moreover, if $G_0 \cong \Z^l$ then $P_0$ is diffeomorphic to 
$P_0' \# (\#_{i=1}^{l-1} (S^{4j-1} \times S^{4j}))$ where $H^{4j}(P_0') \cong \Z$.
\end{Proposition}
\noindent
If, however, $G = H^{4j}(P)$ has a nontrivial torsion subgroup $TG$ and a nontrival free subgroup then
the classification becomes a more complicated matter because in this case a highly connected
$(8j-1)$-manifold $P$, ($j=1,2$), does not necessarily define a unique quadratic linking function on
$TG$.  We are therefore lead to the notion of a quadratic linking family which we now explain.  Firstly,
let $\tau_G$ denote the canonical homomorphism
$$\begin{array}{cccc}
\tau_G: & G & \ra & G^{**}\\
& g & \mapsto & (x \mapsto x(g))
\end{array}$$
and let $\Sec(\tau_G)$ denote the set of homomorphisms $\Phi: G^{**} \hra G$ such that $\tau_G \circ \Phi =
Id_{G^{**}}$.  A quadratic linking family is a triple $(G,q^*,\beta)$ consisting of the finitely generated
abelian group
$G$, an element 
$\beta \in G$ and a function $q^*: \Sec(\tau_G) \ra \QL(b)$ for some linking form on $TG$.  The
superscript $*$ takes the values $c$ for {\em characteristic} and $ev$ for {\em even}.  When $j$ is greater
than $2$, highly connected $(8j-1)$-manifolds define even quadratic linking families on their $4j$th cohomology
groups.  The function $q^{ev}$ is always a constant function whose value is a quadratic linking form.  Thus an
even quadratic linking family is equivalent to a triple $(G,q,\beta)$ where $q$ is a quadratic linking
form on $TG$. However, in dimensions $7$ and $15$, highly connected manifolds define characteristic
quadratic linking families whose function $q^c$ may be far from constant.  We refer the reader to
Section \ref{chapprelim}.\ref{secqlfam} for the precise definitions of quadratic linking families and isometries
between them.  A triple $(H,\lambda,\alpha)$ (where $\lambda$ may now be a degenerate bilinear form) induces a
quadratic linking family, $Q(\lambda,\alpha)$ on $\Cok(\hl) = G$.  Given $P$, we define $Q(P)$ to be
$Q(\lambda,\alpha)$ where $P \# \Sigma_P = \del L$ for a handlebody with triple $(H,\lambda,\alpha)$
(again see Section \ref{chapprelim}.\ref{secqlfam}).

\begin{thm} \label{thmb}
Let $j=1$ or $2$.  For every finitely generated abelian group $G$ and every characteristic 
quadratic linking family $Q$ defined on $G$, there is a highly connected $(8j-1)$-manifold $P$ such
that $Q(P) \cong Q$.  The quadratic linking family $Q(P)$ of a highly connected $(8j-1)$-manifold $P$ is
a  complete invariant of almost diffeomorphisms.  That is, if $P_0$ and $P_1$ are two highly connected 
$(8j-1)$-manifolds then there is an almost diffeomorphism $f: P_0 \ra P_1$ such that $f_! = \theta:
H^{4j}(P_0) \cong H^{4j}(P_1)$ if and only if $\theta$ defines an isometry of quadratic linking families
from $Q(P_0)$ to $Q(P_1)$.
\end{thm}

\noindent
For comparison we present the analogous theorem in higher dimensions which follows rapidly from a
compilation of Wall's results on highly connected manifolds.

\begin{Theorem} [Wall \cite{Wa1,Wa2,Wa3,Wa4}] \label{thmwallqlfam}
Let $j>2$.  For every finitely generated abelian group $G$ and every even quadratic linking
family $Q$ on $G$ there is a highly connected $(8j-1)$-manifold $P$ such that
$Q(P) \cong Q$.  The quadratic linking family $Q(P)$ of a highly connected $(8j-1)$-manifold $P$ is a 
complete invariant of almost diffeomorphisms. 
\end{Theorem}

The remainder of this thesis is organized as follows.  In Chapter \ref{chapprelim} we develop the
necessary preliminaries.  We first summarize Wall's classification of the handlebodies $L$ and rephrase it in
terms of quadratic functions.  Here a quadratic function on a free abelian group $H$ is a function
$\kappa(\lambda,\alpha):H \ra \Z$ such that $\kappa(v) = \lambda(v,v) + \alpha(v)$ for $\lambda$ a symmetric
bilinear form on $H$ and $\alpha$ a linear function on $H$.  We then show that after the addition of a homotopy
sphere, every manifold $P$ bounds a handlebody $L$ and present Wilkens' Theorem (\ref{thmwilk3.2}) on the
extension of diffeomorphisms of the manifolds $P$ to cobounding handlebodies.  At this stage it is clear that
the almost diffeomorphism classification of the manifolds $P$ is equivalent to the stable classification of
quadratic functions where ``stable" means ``modulo nonsingular quadratic functions".  We therefore devote
the sections which follow to the algebraic invariants that a quadratic function, $\kappa(\lambda,\alpha)$,
induces on its cokernel, $G=\Cok(\hl)$, which are stably invariant.  These invariants are, in increasing
order of sophistication: a linking form on the torsion subgroup of $G$, a quadratic refinement of this linking
form and, in the case where $\kappa$ in degenerate, a quadratic linking family.  In the final two sections
of the chapter we provide categorical statements of our classification and define the invariants $s_1$ and
$\bs$. 

Chapter \ref{chapalg} is an entirely algebraic chapter concerning the stable classification of quadratic
functions.   We prove that quadratic functions are stably classified by the quadratic linking families which
they induce.  We first consider the case of nondegenerate quadratic functions where the key idea is that
quadratic linking functions behave like algebraic boundaries to the quadratic functions which induce them.  In
particular, the Gluing Lemma \ref{chapalg}.\ref{lemmaglue} states that if quadratic functions $\kappa_0$ and
$\kappa_1$ induce quadratic linking functions $q_0$ and $q_1$, then an isometry,
$\theta:q_0 \cong q_1$, may be used to glue $\kappa_0$ to $\kappa_1$ and form a nonsingular quadratic function
which we denote $\kap$.  This procedure is a precise algebraic analogue of gluing two handlebodies together
along a diffeomorphism of their boundaries.  We are then able to mimic the proof of Wilkens' Theorem
\ref{thmwilk3.2} in the algebraic setting and show that the isometry $\theta:q_0 \ra q_1$ above is induced by
an isometry of quadratic functions $\Theta: \kappa_0 \oplus \kap \ra \kappa_1 \oplus \mu$ where $\mu$ is a
nonsingular quadratic function.  Taking a little more care we seen that the somewhat complex notion of a
quadratic linking family is the required algebraic object to generalize this result for degenerate
quadratic functions.  We note that our results in this chapter generalize the results of Wall and
Durfee on the stable classification of nondegenerate even quadratic forms (in a natural way) to the stable
classification of nondegenerate characteristic quadratic functions.

In Chapter \ref{chaptop} we prove our central classification results, Theorem \ref{thmb} and Theorem
\ref{thmwallqlfam}, by applying Wilkens' Theorem \ref{thmwilk3.2} to the algebraic results obtained in Chapter
\ref{chapalg}.  As well as completing the classification left open by Wilkens' this gives a new and simpler
proof of Wall's Theorem \ref{thmwallqlfam}.  In Section \ref{seckap} we use the $h$-cobordism theorem to show
that the algebraic Gluing and Splitting Lemmas (\ref{chapalg}.\ref{lemmaglue} and
\ref{chapalg}.\ref{lemmasplit}) from Chapter \ref{chapalg} are mimicked precisely in the world of
handlebodies.  In the case of rational homology spheres, this gives a second proof of Theorem \ref{thmb}
and Theorem
\ref{thmwallqlfam} which is independent of Theorem
\ref{thmwilk3.2}.  In the last section of the chapter we present the  smooth classification of highly
connected $(8j-1)$-manifolds.

In Chapter \ref{chaplink} we review the classification of linking forms and quadratic linking forms on finite
abelian groups due to Kawauchi and Kojima \cite{KK} and Nikulin \cite{Ni} respectively and we then extend
Nikulin's classification of quadratic linking forms to the classification of quadratic linking functions.  
Specifically, let  $G$ be a finite abelian group and let $q$ be a quadratic linking function on $G$ which
refines the linking form $b(q)$.  There is a quadratic linking form $q^o$ such that $q=q^o_a$ for $a \in
G$.  If we define $\beta(q) := 2a$ then the pair $(b(q),\beta(q))$ is an invariant of $q$ (see
\ref{remwilkinvt}).  In addition to the invariant $(b(q),\beta(q))$ there is a Kervaire-Arf invariant,
$K(q) \in \Q/\Z$, defined for each quadratic linking function $q$ via a Gauss sum (see Section
\ref{chaplink}.\ref{secquadf}).  Our classification result is the following

\renewcommand{\thetheorem}{\ref{thmqf}}

\begin{theorem}
Let $q_0$ and $q_1$ be quadratic linking functions over a finite abelian group $G$.  
Then $q_0$ is isometric to $q_1$ if and only if $(b(q_0),\beta(q_0)) \cong (b(q_1),\beta(q_1))$ and $K(q_0) =
K(q_1)$.
\end{theorem}
\noindent
If $q^c(P)$ is the quadratic linking function of a highly connected $(8j-1)$-manifold ($j=1,2)$, we also show
that $K(q(P)) = -\bs(P)$.  Hence Theorem \ref{thmqf} is the algebraic counter-part to {\em 3} of Theorem
\ref{thma}.

In Chapter \ref{chaphom} we use smoothing theory and surgery theory to deduce respectively the topological and
homotopical classification of highly connected $7$-manifolds from the almost diffeomorphism
classification of these manifolds.  We conclude by applying our results to give the
classification of the total spaces of $3$-sphere bundles over the $4$-sphere which appears in \cite{CE}.

\chapter{Preliminaries} \label{chapprelim}
\section{Manifold essentials} \label{secmaness}
Throughout this thesis all manifolds shall be compact and oriented and all maps between manifolds
shall be assumed orientation-preserving.  If $M$ is a manifold then $-M$ shall denote the manifold with
opposite orientation.  Unless otherwise noted manifolds shall also be smooth.   The letter $L$
shall always denote a $2k$-dimensional handlebody by which we mean a manifold obtained by attaching
some number of $k$-handles to the $2k$-ball.  We review the definition and theory of handlebodies
in Section \ref{secwall}.  For all $k>2$ every highly connected $2k$-manifold with a single boundary
component is diffeomorphic to a handlebody.  The symbol $\del$ shall be used to indicate the
operation of taking the boundary of a manifold.  The boundary of every handlebody is a highly
connected manifold of dimension $(2k-1)$.  The symbol $M_0 \# M_1$ denotes the (interior) connect
sum of the manifolds $M_0$ and $M_1$.  The symbol $L_0 \natural L_1$ denotes the boundary connect
sum of $L_0$ and $L_1$.  We refer the reader to \cite{KM} for definitions.   It is a consequence of
these definitions that
$$\del(L_0 \natural L_1) = \del L_0 \# \del L_1.$$

A simply connected $h$-cobordism is a simply connected compact $(n+1)$-manifold $W$ whose boundary has
two components, $\del W = M_0 \sqcup -M_1$, such that each inclusion $M_i \hra W$ is a homotopy
equivalence for $i=0,1$.  The $h$-cobordism theorem of Smale asserts that when $n>4$ every simply
connected $h$-cobordism is diffeomorphic to the product $M_0 \times [0,1]$ via a diffeomorphism that is
the identity restricted to $M_0$.

A homotopy sphere is a manifold which is homotopy equivalent to the standard sphere.  
We shall always denote homotopy spheres by the letter $\Sigma$.  In dimensions higher than $4$
Kervaire and Milnor \cite{KM} famously classified the groups of diffeomorphism classes of homotopy
spheres under the operation of connected sum.  The group of $n$-dimensional homotopy spheres is denoted
by $\Theta_n$ and the subgroup of those spheres bounding parallelizable manifolds by $bP_{n+1}$.  There
is an exact sequence
$$ 0 \lra bP_{n+1} \stackrel{i}{\lra} \Theta_n \lra \Cok(J_n) $$
where $\Cok(J_n)$ denotes the cokernel of the $J$-homomorphism, $J_n: \pi_n(O) \ra \pi_n^S$ (see
\cite{Wh} pg 502-4).  When $n=4k-1$, $\Theta_{4k-1} \cong bP_{4j} \oplus \Cok(J_n)$ and 
Brumfiel defined a splitting
$s: \Theta_{4k-1} \ra bP_{4k}$ of $i$ as follows.  The group $bP_{4k} \cong \Z/|bP_{4k}|$ is cyclic (of
order determined by the $k$th Bernoulli number).  Brumfiel \cite{Bru} showed that every homotopy sphere
$\Sigma
\in \Theta_{4k-1}$ bounds a spin manifold $W$ with vanishing decomposable Pontryagin numbers.  The
homomorphism
$s$ is defined by $\Sigma \mapsto \sigma(W) (\modu |bP_{4k}|)$ where $\sigma(W)$ denotes the signature of
the $4k$-manifold $W$.  For the dimensions of primary concern in this thesis, namely $7$ and $15$, the
groups of exotic spheres are:
$$  \Theta_7 = bP_8 \cong \Z/28 \hskip 2cm \Theta_{15} \cong bP_{16} \oplus \Z/2 \cong \Z/8128 \oplus
\Z/2.$$

Let $M_0$ and $M_1$ be smooth manifolds.  A homeomorphism $f: M_0 \ra M_1$ which is a
diffeomorphism except at a finite number of singular points shall be called an almost diffeomorphism. 

\begin{Proposition}  \label{propad}
Let $f: M_0 \ra M_1$ be an almost diffeomorphism of $n$-dimensional manifolds.  Then there is an
embedded disc $D^n \hra M_0$ containing all the singular points of $f$, a homotopy sphere $\Sigma$ and
a diffeomorphism $f':M_0 \# \Sigma \ra M_1$ such that $f'|_{M_0-D^n} = f|_{M_0-D^n}$. (Here it is
understood that the connected sum of $M_0$ and $\Sigma$ takes place inside $D^n$.)
\end{Proposition}

\begin{proof}
Let $p_1, \dots p_l$ be the singular points of $f$ and let $I_i:D^n_i\hra M_0$ be $l$ pairwise disjoint
embeddings of discs with $I_i(0) = p_i$.  The image of each $D^n_i$ under $f$ is a subset of
$M_1$ homeomorphic to the $n$-disc and restricting the smooth structure of $M_1$ to $f(D^n_i)$ yields a
smooth manifold which is diffeomorphic to the $n$-disc.  So let $J_i:D^n_{i1}$ be embeddings of discs
in $M_1$ such that $J_i(D^n_{i1}) = f(I_i(D^n_i))$ and $J_i(0) = f(p_i)$.  There are then
diffeomorphisms $f_i:(D^n-0) \ra (D^n_{i1}-0)$ defined by $f_i = J_i^{-1} \circ f \circ I_i|_{D^n_i -
0}$ and we construct the homotopy spheres $\Sigma_i := -D^n_i \cup_{f_i} D^n_{i1}$.  Using the obvious
embedding of $-D^n_i$ into $\Sigma_i$ and $I_i$ we define the connected sum 
$M_0 \# \Sigma_1 \# \dots \# \Sigma_l$.  There is then a diffeomorphism
$$\begin{array}{cccc}
f^{''}: & M_0 \# \Sigma_1 \# \dots \# \Sigma_l & \ra & M_1 \\
& M_0 - \{p_0, \dots, p_l \} \ni x & \mapsto & f(x) \\
& D^n_{i1} \ni x & \mapsto & J_i(x). 
\end{array}$$
Now let $\Sigma = \Sigma_1 \# \dots \# \Sigma_l$ and let $D^n$ be a disc in $M_0$ containing
$\cup_{i=1}^l D^n_i$ in its interior.  There is a diffeomorphism 
$f^{'''}:M_0 \# \Sigma \ra M_0 \# \Sigma_1 \# \dots \# \Sigma_l$ such that $f^{'''}|_{M^0-D^n} =
Id_{M^0-D^n}$.  We may then take $f'$ to be the diffeomorphism $f^{''}\circ f^{'''}$.
\end{proof}

\noindent
Now let $f_0:M_0 \ra M_1$ and $f_1:M_1 \ra M_2$ be a pair of almost diffeomorphisms with singular
points $p_0,\dots, p_{l_0}$ and $q_0,\dots, q_{l_1}$ respectively.  Their composition, $f_1 \circ f_0$
is a homeomorphism with singular points $p_0,\dots,p_{l_0},f_0^{-1}(q_0),\dots,f_0^{-1}(q_{l_1})$ (some
of which may be repeated) and hence an almost diffeomorphism.  Thus, for any class of smooth manifolds
there is a corresponding category whose morphisms are almost diffeomorphisms.  We work primarily in the
almost diffeomorphism category because this is the level at which the algebra of the manifolds $L$ and
$P$ is the simplest.  

By an almost smooth manifold we shall mean a $PL$ $n$-manifold, $M$, with a smooth 
structure except at a finite number of singular points, $p_1, \dots, p_l$.  We may assume that there
are mutually disjoint discs $D^n_i$ with $p_i \in {\rm int}(D^n_i)$.  The manifold $M -
(\sqcup_{i=1}^l D^n_i)$ is a smooth manifold with boundary $l$ disjoint homotopy spheres
$\Sigma_1, \dots, \Sigma_l$ and the $PL$ structure on $D^n_i$ is obtained by coning on $\Sigma_i$. 
Almost smooth manifolds shall arise in this thesis in two ways.  Firstly, we shall sometimes consider a
handlebody $L$ with boundary an exotic sphere $\Sigma$ and form the almost smooth manifold $L \cup
D^{n}$ with singular point at the center of the disc
$D^{n}$.  At other times we shall consider two handlebodies $L_0$ and $L_1$ and an almost
diffeomorphism $f:\del L_0 \ra \del L_1$.  The adjunction space $L_0 \cup_f -L_1$ has the structure of
an almost smooth manifold with singular points corresponding to the singular points of $f$.  An almost
smooth $h$-cobordism between closed simply connected smooth manifolds $M_0$ and $M_1$ is a compact
almost smooth manifold $W$ with boundary $M_0 \sqcup -M_1$, singular points contained in
$W - (M_0 \sqcup M_1)$ and with the inclusions $M_i \hra W$ homotopy equivalences.  We now demonstrate
that the boundaries of such an almost smooth $h$-cobordism are almost diffeomorphic.  Let the singular
points of $W$ be $p_0, \dots, p_n$ and suppose that they are contained in pair-wise disjoint closed
discs $D^n_i$ whose boundaries are diffeomorphic to $\Sigma_i = \del D^n_i$ with the smooth structure
they inherit from $W$.  In addition let $\gamma_i: [0,1] \times {\rm int}(D^{n-1}) \hra W$ be
pair-wise disjoint smooth embeddings such that $\gamma_i^{-1}(D^n_i) = 1 \times {\rm
int}(D^{n-1}_i)$ and $\gamma_i^{-1}(\del W) = 0 \times {\rm int}(D^{n-1}_i) \subset M_0$. 
The manifold $W - (\bigsqcup_i {\rm Im}(\gamma_i))$ is a smooth simply connected $h$-cobordism between
with boundary $M_0 \# (\Sigma_1 \# \dots \# \Sigma_n) \sqcup -M_1$ and hence $M_0$ and
$M_1$ are almost diffeomorphic by the $h$-cobordism theorem.




\section{Wall's classification of the handlebodies $L$} \label{secwall}
We adopt the notation of \cite{Sm} and consider manifolds $L=L(D^{2k},\phi_1, \dots, \phi_l;k)$ where
$\phi_i:(\del D^{k} \times D^k)_i \ra D^{2k}$ are embeddings and $L$ is obtained from $D^{2k} \sqcup
\sqcup_{i=1}^l (D^k \times D^k)_i$ by identifying points using the
$\phi_i$ and then smoothing corners.  The embedding 
$$\phi=\sqcup_{i=1}^l \phi_i : \sqcup_{i=1}^l (S^{k-1} \times D^k)_i \hra \del D^{2k}$$
is called a presentation of $L$.  Smale denoted
the set of manifolds formed in this way by $\HH(2k,l,k)$ and defined $\HH(k) :=
\bigcup_{l=0}^{\infty}\HH(2k,l,k)$.  The isotopy extension theorem entails that the
diffeomorphism type of $L(D^{2k},\phi_1, \dots , \phi_l;k)$ is determined by the isotopy class of the
embedding $\phi$.  Two pieces of data determine the isotopy class of $\phi$.  Firstly if one considers
$\phi$ restricted to the attaching spheres of each handle one obtains a link 
$$\bar{\phi} = \sqcup_{i=1}^l\bar{\phi}_i :\sqcup_{i=1}^l (S^{k-1} \times 0)_i \hra S^{2k-1}.$$ 
When $k>2$, Smale proved in \cite{Sm} that the isotopy class of $\bar{\phi}$ is completely determined by
the linking numbers of the attaching spheres.  Setting $S^{k-1}_i = (S^{k-1} \times 0)_i$, the linking
number of $\bar{\phi}_i(S^{k-1}_i)$ with $\bar{\phi}_j(S^{k-1}_j)$ is $c_{ij}$ where
$\bar{\phi}_i(S^{k-1}_i)$ is homologous to $c_{ij}$ times $\phi(0 \times \del D^k)_j$ in 
$H_{k-1}(\del D^{2k-1} - \bar{\phi}_j(S^{k-1}_j))$.  Secondly, for each $\phi_i$ there are framing
invariants, $\bar{\alpha}_i \in \pi_{k-1}(SO(k))$, which are defined as follows.  The subspace $N :=
\phi_i(S^{k-1} \times D^k)$ is a tubular neighborhood of $\bar{\phi}_i(S^{k-1}_i)$ in
$S^{2k-1}$.  There is a framing of $N$, $\psi_i: N \cong \bar{\phi}_i(S^{k-1}_i) \times D^k$, uniquely
defined up to homotopy, such that $\psi$ extends to a framing of the normal bundle of a $k$-disc which
is embedded in $D^{2k}$ and which bounds
$\bar{\phi}_i(S^{k-1}_i)$.   After isotopy of
$\phi_i$ we may assume that the composition $\psi_i \circ \phi_i$ has the form
$$\begin{array}{cccc}
\psi_i \circ \phi_i:& (S^{k-1} \times D^k)_i & \ra & \bar{\phi}_i(S^{k-1}_i) \times D^k\\
& (x,y) & \mapsto & (\bar{\phi}_i(x),a(\phi_i)(x)(y))
\end{array}$$
for some map $a(\phi_i):S^{k-1} \ra SO(k)$.  While the map $a(\phi_i)$ depends upon the choices made
the homotopy class of $a(\phi_i)$, which we denote by $\bar{\alpha}_i$, does not.

\begin{Lemma}[Wall \cite{Wa1}] \label{lemma1wall1}
The invariants $c_{ij}, (i \neq j)$ and $\bar{\alpha}_i$ form a complete set of isotopy invariants of
the presentation $\phi$.
\end{Lemma}

Wall then proceeds to identify abstract invariants for handlebodies $L$ which are counter parts 
to the invariants of a presentation but which do not depend upon giving a presentation for $L$. 
Viewed abstractly, the manifold $L$ has invariants a free abelian group $H = H_k(L)$, an
intersection pairing $\lambda: H \times H \ra \Z$, and a map $\bar{\alpha}: H
\ra \pi_{k-1}(SO(k))$.  Geometrically $\lambda$ is defined as follows.  By the Hurewicz theorem
$H_k(L) \cong \pi_{k}(L)$ and so every homology class $v \in H$ is represent by a map $\bar{v}:S^k
\ra L$.  By the embedding theorems of Whitney this map may be chosen to be an embedding.  Given
$v,w \in H$ we may assume further that the embeddings $\bar{v}: S^k \ra L$ and $\bar{w}: S^k \ra L$ are
in general position and so intersect at isolated points.  The integer $\lambda(v,w)$ is the signed sum
of these intersection points where the sign comes from comparing the orientation of $L$ with
orientation induced by $\bar{v}$ and $\bar{w}$ at a particular intersection point.  The map 
$\bar{\alpha}$ evaluated on a homology class $v$ is given by mapping the normal bundle of
 $\bar{v}:S^k \hra L$ to the homotopy class of its clutching function which is an element of
$\pi_{k-1}(SO(k))$.  If $\del \iota_k$ denotes the element of $\pi_{k-1}(SO(k))$ corresponding to the
tangent bundle of the $k$-sphere, then Wall proves:

\begin{Lemma}[Wall \cite{Wa1}] \label{lemma2wall1}
$$\bar{\alpha}(v + w) = \bar{\alpha}(v) + \bar{\alpha}(w) + \lambda(v,w)\cdot \del \iota_k.$$
\end{Lemma}
\noindent
This allows us to translate between the abstract invariants of $L$ and the invariants of the
presentation for $L = L(D^{2k},f_1,\dots ,f_n;k)$.  If one attaches a disc $D^k_i$ inside $D^{2k}$
to the core of each handle $(D^k \times 0)_i$ then one obtains a $k$-sphere $\bar{v_i}$ which one may
assume is embedded in $L$.  The homology classes $v_i$ corresponding to these embedded spheres
form a basis for $H = H_k(L)$.  Moreover the spheres $\bar{v_i}$ intersect only inside $D^{2k}$
and their intersection number $\lambda(v_i,v_j)$ is coincides with linking number $c_{ij}$.  The self
linking numbers $\lambda(v_i,v_i)$ are given by the Euler numbers $E(\bar{\alpha}_i)$ of the framings
from the presentation invariants.  Finally the framings 
$\bar{\alpha}_i$ are evidently given by $\bar{\alpha}(v_i)$ and Lemma \ref{lemma2wall1} ensures
that the framings $\bar{\alpha}_i$ determine $\bar{\alpha}$.  Wall is able to complete this translation
between presentations and the abstract point of view by proving

\begin{Theorem}[Theorem 1 \cite{Wa1}] \label{thm1wall1}
Suppose $L \in \HH(k)$, and choose a basis for $H_k(L)$.  Then $L$ has a presentation
corresponding to this basis.
\end{Theorem}

For much of the rest of this thesis $k$ shall be equal to $4j$ and we shall be concerned almost
exclusively  with $8j$-dimensional handlebodies and more particularly with the boundaries of
$8j$-dimensional handlebodies.  It will be algebraically convenient to consider the stabilization of the
tangential invariant $\bar{\alpha}$ which is $\alpha := S \circ \bar{\alpha}$, where $S:
\pi_{4j-1}(SO(4j)) \ra \pi_{4j-1}(SO) \cong \Z$ is the stabilization map.  Since $S(\del \iota_k) = 0$,
Lemma \ref{lemma2wall1} implies that $\alpha: H \ra \Z$ is a homomorphism.  In fact, 
Kervaire computed that $\alpha = a_j(2j-1)!p_j$ where $p_j$ denotes the $j$th Pontryagin class and
$a_j = 1$ if $j$ is even and $2$ if $j$ is odd, \cite{Ke}.  Note that we lose nothing by considering
$\alpha$ over $\bar{\alpha}$ since the unstable part of $\bar{\alpha}(v)$ is given by the Euler number of
the normal bundle of $\bar{v}$ and this is equal to $\lambda(v,v)$.  Hence $\bar{\alpha}$ may be
recovered from
$\alpha$ and $\lambda$.  We conclude this section with a summary in the form of a Corollary.

\begin{Corollary} \label{corsechandle}
Let $k>2$ and let $L$ be a highly connected $2k$-manifold with a non-empty and connected boundary. 
Then $L$ is diffeomorphic to a handlebody.  If $L_0$ and $L_1$ are two handlebodies of dimension $2k$
with invariants $(H_0,\lambda_0,\alpha_0)$ and $(H_1,\lambda_1,\alpha_1)$ then there is diffeomorphism
$f:L_0 \ra L_1$ such that $f_* = \Theta$ if and only if there is an isomorphism $\Theta: H_0 \ra H_1$ 
such that $\lambda_1(\Theta(v),\Theta(w)) = \lambda_0(v,w)$ and  $\alpha_1(\Theta(v)) = \alpha_0(v)$ for
all $v,w \in H_0$.
\end{Corollary}

\begin{remark} \label{rempltophb}
Note that the $PL$-invariance of rational Pontryagin classes (see \cite{MS} \S 20) entails that the
words words ``diffeomorphic'' and ``diffeomorphism'' in Corollary \ref{corsechandle} may be
replaced by ``$PL$-homeomorphic'' and``$PL$-homeomorphism''.  Novikov's result on the topological
invariance of the rational Pontryagin classes entails that they may be replaced by  ``homeomorphic''
and  ``homeomorphism''.
\end{remark}


\section{Quadratic functions} \label{secqf}
In the previous section we saw Wall's classification of the handlebodies $L$ via
their intersection form $\lambda(L):H_{4j}(L) \times H_{4j}(L) \ra \Z$ and their stable
tangential invariant $\alpha(L) \in H^{4j}(L)$.  In this section we give an abstract summary of
the algebraic invariants of $L$.  

To begin let $H$ be a finitely generated free abelian group.   A linear form on $H$ is simply a
group homomorphism $\alpha:H \ra \Z$.  The set of all such forms with
the operation of pointwise addition of homomorphisms forms a free abelian group  called the dual
group of $H$ and denoted by $H^*$.  A symmetric bilinear form on $H$ is a bilinear map $\lambda: H
\times H \ra \Z$ such that $\lambda(v,w) = \lambda(w,v)$ for all $v,w \in H$.  There are two ways one
might combine a bilinear form $\lambda: H \times H \ra \Z$ and a linear form $\alpha: H \ra \Z$.   One
may follow Wall and simply consider the pair $(\lambda,\alpha)$.  A morphism of pairs
$(\lambda_0,\alpha_0)$ $(\lambda_1,\alpha_1)$ is a homomorphism $\Theta: H_0 \ra H_1$ such that
$\alpha_0=\alpha_1 \circ \Theta$ and $\lambda_0 = \lambda_1 \circ (\Theta \times \Theta)$.  Another way
to combine $\lambda$ and $\alpha$ is to take the quadratic form associated to $\lambda$ and add $\alpha$
to it.  

\begin{Definition}
A function $\kappa: H \ra \Z$ is called a {\em quadratic function} if there exists a bilinear form
$\lambda$ on $H$, a linear functional $\alpha \in H^*$ such that
$$\kappa(v) = \lambda(v,v) + \alpha(v)$$ 
for all $v \in H$.  
\end{Definition}

\begin{Proposition}
\label{propqf}
The form $\lambda$ and dual element $\alpha$ are uniquely recovered from $\kappa$ using the
following relations:
$$\begin{array}{c}
\kappa(v + w) = \kappa(v) + \kappa(w) + 2\lambda(v,w)\\
 \alpha(v) = \kappa(v) - \lambda(v,v).
\end{array}$$
\end{Proposition}

\begin{proof}
The first relation shows that $\lambda$ is determined by $\kappa$ and the second shows how to
recover $\alpha$ from $\kappa$ given that $\lambda$ is known.
\end{proof}

\noindent
When we wish to specify the components of a quadratic function we shall write $\kappa =
\kappa(H,\lambda,\alpha)$ for a quadratic function with domain $H$, bilinear form $\lambda$
and linear form $\alpha$.  On occasions we shall drop any of $H,\lambda$ or $\alpha$ from the list
when they are understood or not relevant.  Note that $\kappa$ is linear in the case where $\lambda$
is identically zero.   When $\alpha = 0$, the quadratic function $\kappa$ is homogeneous, in which case
$\kappa$ is a {\em quadratic form}.  

We shall see below that what distinguishes the quadratic functions which arise from handlebodies of
dimension $8$ or $16$ is that they take on only even integer values.  We therefore make the following

\begin{Definition}
An element $\alpha \in H^*$ is called characteristic for $\lambda$ when $\kappa(\lambda,\alpha)$
takes on only even values.  That is, $\lambda(v,v) + \alpha(v)$ is even for all $v \in H$.  In
this case we shall call $\kappa(\lambda,\alpha)$ a characteristic quadratic function.  
\end{Definition}

Now let $\kappa_0 = \kappa(H_0,\lambda_0,\alpha_0)$ and $\kappa_1 =
\kappa(H_1,\lambda_1,\alpha_1)$ be quadratic functions on $H_0$ and $H_1$ respectively and let
$\Theta: H_0 \ra H_1$ be a homomorphism.  We define $\Theta^*\kappa_1$ by
$\Theta^*\kappa_1(v)=\kappa_1(\Theta(v))$.  A homomorphism $\Theta$ is called an {\em isometric
embedding} if it is injective and $\Theta^*\kappa_1(v) = \kappa_0(v)$ for every $v \in H_0$.  A bijective
isometric embedding is called an {\em isometry}.  We observe that an isometry of quadratic functions
induces an isometry of both the homogeneous quadratic and linear parts and record this fact as a
proposition.

\begin{Proposition} \label{propqfrecov}
Let $\lambda_0$ and $\lambda_1$ be symmetric bilinear forms on free abelian groups $H_0$ and $H_1$
respectively and let $\alpha_0$ and $\alpha_1$ be elements of $H_0^*$ and $H_1^{*}$ respectively.  Then
the the quadratic functions $\kappa(H_0,\lambda_0,\alpha_0)$ and $\kappa(H_1,\lambda_1,\alpha_1)$ are
isometric if and only if the pairs $(\lambda_0,\alpha_1)$ and $(\lambda_1,\alpha_1)$ are isomorphic.
\end{Proposition}
\begin{proof}
This follows immediately from Proposition \ref{propqf}.
\end{proof}
\noindent
Thus the two methods for combining $\lambda$ and $\alpha$ which we mentioned earlier are equivalent.  

Now let $\kappa = \kappa(H,\lambda,\alpha)$ be a quadratic function.   The adjoint of the bilinear form
$\lambda$ is the linear transformation induced by $\lambda$ from $H$ to its dual,
$$
\begin{array}{rcl}
\hat{\lambda}: H & \ra & H^*\\
v & \mapsto & [w \mapsto \lambda(v,w)].
\end{array}
$$
Thus for all $v,w \in H$, $\lambda(v,w) = \hat{\lambda}(v)(w)$.  We define the {\em radical} of $\lambda$
to be the group $F:=\Ker(\hl) = \{v \in H | \lambda(v,w) = 0~\forall w \in H \}$ and the quotient of
$\lambda$ to be the group $G:=\Cok(\hl) = H^*/\hl(H)$.  The groups $F$, $H$, $H^*$ and $G$ fit into the
following exact sequence:
$$0 \lra F \lra H \stackrel{\hl}{\lra} H^* \stackrel{\pi}{\lra} G \lra 0$$
where $\pi$ is the canonical projection from $H^*$ to $H^*/\hl(H) = G$.  We call this sequence the {\em
fundamental sequence} of $\lambda$ and $\kappa$ as it will occupy us for most of the remainder of this
thesis.  If $F=0$ then $\lambda$ and $\kappa$ are called {\em nondegenerate}.  If in addition $G=0$ then
$\lambda$ and $\kappa$ are called {\em nonsingular}.  

Given two quadratic functions $\kappa_0~=~\kappa(H_0,\lambda_0,\alpha_0)$ and $\kappa_1 =
\kappa(H_1,\lambda_1,\alpha_1)$, the direct sum $\kappa_0 \oplus \kappa_1$ is defined on $H_0 \oplus
H_1$ via the equation
$$(\kappa_0 \oplus \kappa_1)[(v_0,v_1)] = \kappa_0(v_0) + \kappa_1(v_1).$$

We conclude this section by making the relation between handlebodies and quadratic functions explicit.

\begin{Definition} \label{defqfhb}
Let $L \in \HH(4j)$ be a handlebody with intersection form $\lambda(L)$ and tangential invariant
$\alpha(L)$ defined on $H_{4j}(L)$, then $\kappa(L) := \kappa(H_{4j}(L),\lambda(L),\alpha(L))$
shall be called the {\em quadratic function of $L$}.  Conversely, if $\kappa = \kappa(H,\lambda,\alpha)$
is a quadratic function such that there is an $8j$-dimensional handlebody $L$ with $\kappa(L) \cong
\kappa(H,\lambda,\alpha)$ then $\L^{8j}(\kappa)$ shall denote the $8j$-dimensional handlebody defined up
to diffeomorphism by $\kappa$.
\end{Definition}
\noindent
Now let $P$ denote the boundary of $L$, let $PD:H_{4j}(L,P) \cong H^{4j}(L)$ denote the Poincar\'{e}
duality isomorphism and let  $j:H_{4j}(L) \ra H_{4j}(L,P)$ denote the map induced by the map of pairs
$(L,\phi) \ra (L,P)$.  The following fact is fundamental and well known.  If we use the Kronecker pairing
to identify $H^{4j}(L)$ with $\Hom(H_{4j}(L),\Z)$ then $PD \circ j: H_{4j}(L) \ra H^{4j}(L)$
is the adjoint of the intersection pairing on $L$.  Thus, given Definition \ref{defqfhb} we have the
following pair of isomorphic sequences 
\[
\divide\dgARROWLENGTH by2
   \begin{diagram}
\node{0} \arrow{e}  
      \node{H_{4j}(P)} \arrow{e,t}{i_*} \arrow{s,<>}
       \node{H_{4j}(L)} \arrow{e,t}{PD \circ j} \arrow{s,<>}
        \node{H^{4j}(L)} \arrow{e,t}{i^*} \arrow{s,<>}
									\node{H^{4j}(P)} \arrow{e} \arrow{s,<>}
          \node{0}  \\
     \node{0} \arrow{e}
					 \node{F} \arrow{e}
       \node{H} \arrow{e,t}{\hat{\lambda}}
        \node{H^*} \arrow{e,t}{\pi} 
         \node{G} \arrow{e}
	         \node{0}
   \end{diagram}
\]
where $i_*$ and $i^*$ denote the maps induced by the inclusion $i:P \hra L$.  This gives a topological
interpretation of the fundamental sequence of $\kappa(L)$ which we shall use throughout this thesis.
 
It is apparent that the boundary connected sum of handlebodies gives rise to the direct sum of the
associated quadratic functions so that $\kappa(L_0 \natural L_1) = \kappa(L_0) \oplus
\kappa(L_1)$.  If we reverse the orientation on $L$ then this will change the intersection form
of $L$ by reversing its sign.  However, reversing the orientation on $L$ does not affect the
tangential invariant.  This is because $\alpha(L)$ is the primary obstruction to
trivializing the tangent bundle and we noted above that this class is a multiple of the Pontryagin 
class which is not affected by a change of orientation. Thus if
$\kappa=\kappa(H,\lambda,\alpha)$ we define $\kappa^{-}$ to be the quadratic function
$\kappa(H,-\lambda,\alpha)$ and with this definition $\kappa(-L) = \kappa(L)^{-}$.

\section{Wall's classification restated} \label{secwrs}
We may now restate Wall's theorem \ref{thm1wall1} as follows: for every isometry $\Theta :
\kappa(L_0) \cong \kappa(L_1)$ there corresponds a diffeomorphism
$f(\Theta)$ whose induced map on homology is $\Theta$.  We shall now dress this statement up in
categorical terms by defining the following categories.  Let $\HH^{8j}:=\HH(4j)$ be the category whose
objects are highly connected $8j$-manifolds $L$ and whose morphisms are equivalence classes of
diffeomorphisms of handlebodies.  Two diffeomorphisms $g,h:L_0 \ra L_1$ are equivalent if they induce
the same isomorphism on middle homology, $g_* = h_*:H_{4j}(L_0) \ra H_{4j}(L_1)$.  We define $\QF$
to be the category whose objects are quadratic functions on free abelian groups and whose morphisms are
isometries of quadratic functions.  Finally we define $\QF^{8j}$ to be the full
sub-category of $\QF$ whose objects are quadratic functions which are isomorphic to the quadratic
function of some handlebody $L$.  A symmetric monoidal category is one with an associative and
commutative addition operation on its objects that is functorially well behaved (see \cite{MacL} VII 7). 
We shall not distract the reader with the precise definitions but assure him that $\HH^{8j}$ and $\FF^*$
are symmetric monoidal categories with respect to the operations of boundary connected sum of
handlebodies and of direct sum of quadratic functions.

\begin{Theorem}[Restatement of \cite{Wa1} Theorem 1] \label{thmwallrs}
There is an equivalence of (symmetric monoidal) categories 
$$\begin{array}{cccc}
\kappa^j:& \HH^{8j} &  \cong & \QF^{8j} \\
&  \{L,[g], \natural \} & \ra & \{\kappa(L),g_*,\oplus\}.
\end{array}$$
\end{Theorem}

\begin{proof}
Recall that a skeleton of a category ${\mathcal C}$ is a full subcategory of ${\mathcal
C}$ which contains precisely one object for every isomorphic class of objects in ${\mathcal C}$.  Every
category ${\mathcal C}$ is equivalent to every one of its skeletons (\cite{MacL} IV 4).  So let
$\underline{\HH^{8j}}$ be a skeleton of $\HH^{8j}$.  By Corollary \ref{corsechandle}, 
$\kappa^j(\underline{\HH^{8j}})$ is a skeleton of $\QF^{8j}$ and therefore that $\kappa$ is an
equivalence of categories.  Finally, for any pair of handlebodies $L_0$ and $L_1$ and any pair of
diffeomorphisms
$f_0$ and $f_1$, $\kappa(L_0 \natural L_1) = \kappa(L_0) \oplus \kappa(L_1)$ and $(f_0 \natural f_1)_* =
f_{0*} \oplus f_{1*}$.
\end{proof}


We now come to a fundamental observation which Wall and Wilkens did not fully exploit.

\begin{Proposition} \label{propqfc}
When $j \neq 1,2$ the category $\FF^{8j}$ has as objects all quadratic functions
$\kappa(\lambda,\alpha)$ with even bilinear forms $\lambda$ (and no restriction on
$\alpha$).  When  $j=1,2$ the category $\FF^{8j}$ has objects all characteristic quadratic functions
(with no restriction on the type of the bilinear form).
\end{Proposition}

\begin{proof} 
We have seen that a handlebody $L$ is completely described by a presentation 
$\{ (v_1, \dots ,v_l),\{c_{ij}\},\bar{\alpha}(1), \dots, \bar{\alpha}(l) \}$ consisting of a
basis of handles, linking numbers between the handles and the normal bundle of each handle.  With
respect to this basis the intersection matrix of $L$ has components $\lambda_{ij}$ where
$$\lambda_{ij} = \left\{\begin{array}{lr} 
c_{ij} & i \neq j \\
E(\alpha(i)) & i = j
\end{array} \right\} $$
and where $E(\alpha(i))$ is the Euler number of $\bar{\alpha}(i) \in \pi_{4j-1}(SO((4j))$.  When $j \neq
1,2$ Adams' result on the non-existence of Hopf invariant one bundles ensures that each $E(\ba(i))$
is even and hence that $\lambda$ is even.  However, since every set of linking numbers may be realized by
a link of $(4j-1)$-spheres in $S^{4j-1}$ we see that any even form may be produced.  The
stable tangential invariant is determined by the stabilizations
$S(\bar{\alpha}(i)) \in \pi_{4j-1}(SO)$ which are the values of $\alpha$ on the basis elements $v_i$. 
That these can be chosen arbitrarily follows from the result (see \cite[171]{Wa1}) that the image
of the homomorphism
$$S \oplus E:\pi_{4j-1}(SO(4j)) \ra \pi_{4j-1}(SO) \oplus \Z$$
has index $2$.  Since $E$ takes only even values for $j \neq 1,2$ the map
$S$ must be onto $\pi_{4j-1}(SO)$.  When $j=1,2$ then there are bundles with every possible Euler number
and thus any intersection form one desires may be produced.  Moreover, for every $\ba \in
\pi_{4j-1}(SO(4j))$, we now have that $S(\ba)$ and $E(\ba)$ are of the same parity.
Hence for $i=1, \dots, l$,
$$\lambda(v_i,v_i) = E(\ba(i)) \equiv S(\ba(i)) ~(\modu~2) \equiv \alpha(v_i) ~(\modu~2).$$
As $\{v_1, \dots, v_l\}$ is a basis for $H_{4j}(L)$ it follows that $\lambda(v,v) \equiv \alpha(v)
~\modu~(2)$ for every $v \in H_{4j}(L)$.
\end{proof}

As a consequence of Proposition \ref{propqfc} we define three categories of quadratic
functions as follows.
\begin{itemize}
\item{
$\QF^{c}$ shall denote the
full sub-category of $\FF$ whose objects are characteristic quadratic functions with arbitrary
bilinear forms.}
\item{$\QF^{ev}$ shall denote the full sub-category of $\FF$ whose objects are quadratic functions with
even bilinear forms and arbitrary linear terms.}
\item{$\QF^o$ shall denote the full sub-category of $\FF$ whose objects are even quadratic
forms.}
\end{itemize}
Proposition \ref{propqfc} states that
$\QF^{8j} = \QF^c$ when $j=1,2$ and $\QF^{8j} = \QF^{ev}$ when $j>2$.  For all values of $j$ the
category $\QF^o$ corresponds to the category of stably parallelizable handlebodies.  Within each
category $\FF^*$ there are full subcategories of nondegenerate and nonsingular quadratic functions.  They
correspond respectively to the categories of handlebodies cobounding rational homology spheres and
handlebodies cobounding homotopy spheres.  As we are primarily concerned with the classification of the
boundaries of handlebodies up to connected sum with homotopy spheres we define the notion of stable
$\FF^*$-equivalence of quadratic functions.  

\begin{Definition} \label{defstabfequiv}
We call two quadratic functions $\kappa_0$ and $\kappa_1$ belonging to the category $\FF^*$,
stably $\FF^*$-equivalent if there are nonsingular quadratic functions $\mu_0$ and $\mu_1$ also
belonging to
$\FF^*$ and an isometry
$$\Theta: \kappa_0 \oplus \mu_0 \cong \kappa_1 \oplus \mu_1.$$
\end{Definition}

\noindent
We shall denote the stable $\FF^*$-equivalence class of $\kappa \in \FF^*$ by $[\kappa]_{\FF^*}$.  Note
the direct sum of stable $\FF^*$-equivalence classes of quadratic functions is well defined,
$$[\kappa_0]_{\FF^*} \oplus [\kappa_1]_{\FF^*} := [\kappa_0 \oplus \kappa_1]_{\FF^*}.$$

\section{Highly connected bordism} \label{sechcb}
Our perspective on highly connected $(8j-1)$-manifolds is to view them as the boundaries of
handlebodies.  However, not every highly connected $(8j-1)$-manifold $P$ bounds a handlebody.  In
this section we digress somewhat to introduce an invariant of $P$ which is the obstruction to $P$
bounding a handlebody.  With this invariant we shall show that every $P$ bounds a handlebody after the
addition of the right homotopy sphere.  We remind the reader that an $n$-manifold is highly connected
if all of its homotopy groups vanish below its middle dimension which is $[n/2]$.  Now let
$\Omega_{n}^{HC}$ be the set of highly connected bordism classes of highly connected, closed (smooth,
oriented) manifolds.  A highly connected bordism between highly connected
$n$-manifolds $M_0$ and $M_1$ is a compact highly connected $(n+1)$-manifold
$W$ and a diffeomorphism $\del W \cong M_0 \sqcup -M_1$.  The relation of being bordant via a highly
connected manifold defines an equivalence relation.  We shall write $[M]$ for the highly connected
bordism class of a highly connected manifold $M$ and $0$ for the class of the empty manifold. 
Note that $[P] = 0$ if and only if $P$ is the boundary of a handlebody $L$.

\begin{Lemma} \label{lemmahcb1}
When $n$ is even or when $n = 2k-1$ is odd and $k \equiv 0,4,6,7~ ({\rm mod}~8)$ the operations of
connected sum and reversing orientation give rise to well defined operations on $\Omega_{n}^{HC}$ under
which $\Omega_n^{HC}$ is an abelian group.
\end{Lemma}

\begin{proof}
Given highly connected bordism classes $[M_0]$ and $[M_1]$ we define 
$$[M_0] + [M_1] := [M_0 \# M_1] \hskip 1cm \text{and} \hskip 1cm -[M_0] := [-M_0].$$
We now verify that these operations are well defined and that $[M \# -M] = 0$.  If $W$ is a
highly connected bordism from $M_0$ to $M_0'$ then $-W$ is a bordism from  $-M_0$ to $-M_0'$ and so
reversing orientation is well defined on highly connected bordism.  Now, 
$M \# -M = \del ((M-D^n) \times [0,1])$ and if the dimension of $M$ is even, $(M-D^n) \times [0,1]$ is a
highly connected zero bordism of $M \# -M$.  When $n=2k-1$ is odd the additional hypothesis on $k$
entails that $\pi_{k-2}(SO) = 0$.  This ensures that every stable bundle over $S^{k-1}$ is trivial and
hence that $(M-D^n) \times [0,1]$ is
$(k-1)$-parallelizable in the language of \cite{Mi}.  We may therefore perform surgeries in the interior
of $(M-D^n) \times [0,1]$ to kill $\pi_{k-1}$.  The resulting manifold is a highly connected zero
bordism of $M \# -M$.

Suppose that $W_0$ is a highly connected bordism from $M_0$ to $M_0'$.  To show that connected sum is
well defined over highly connected bordism we must show that there is a highly connected bordism from
$M_0 \# M_1$ to $M_0' \# M_1$.  We first note that there is a highly connected bordism $W_1$ from $M_1$
to $-M_1$.  When $n$ is even we take $M_1 \times [0,1]$ and when $n = 2k-1$ we argue as before that
surgery may be performed in the interior of $M_1 \times [0,1]$ to produce a highly connected bordism. 
The band connected sum of $W_0$ and $W_1$ (see \cite{KM} \S 2) is then a highly connected bordism from
$M_0 \# M_1$ to $M_0' \# M_1$.  It is clear that the connected sum operation on highly connected bordism
is associative and commutative.  It follows that $\Omega_n^{HC}$ is an abelian group in the indicated
dimensions.
\end{proof}

\begin{Lemma} \label{lemmahcb2}
Let $k>1$ and $k \equiv 0,4,6,7~ ({\rm mod}~ 8)$.  Then every member of
$\Omega_{2k-1}^{HC}$ is represented by a homotopy sphere and
$$\Omega_{2k-1}^{HC} \cong \Theta_{2k-1}/bP_{2k}.$$
\end{Lemma}

\begin{proof}
We let $$\Omega :\Theta_{2k-1} \ra \Omega_{2k-1}^{HC}$$ be the homomorphism we sends a homotopy sphere
to the highly connected bordism class which it represents.  It follows from calculations of Wall in
\cite{Wa1} that the kernel of $\Omega$ is precisely $bP_{2k}$.  It remains to show that $\Omega$ is onto
but this is a simple consequence of the success of simply connected surgery in odd dimensions as we now
show. One may always take the $(k-1)$-skeleton of a closed highly connected $(2k-1)$-manifold $M$ to be
a wedge of $(k-1)$-spheres.  From the assumption on $k$ it follows that that $\pi_{k-2}(SO) = 0$ and thus
every stable bundle over a wedge of $(k-1)$-spheres is trivial.  In particular the tangent
bundle of $M$ restricted to the $(k-1)$-skeleton of $M$ is trivial.  Thus $M$ is $(k-1)$-parallelizable
and we may perform a series of surgeries on $(k-1)$ spheres in $M$ to kill $\pi_{k-1}(M)$.  The
resulting $(2k-1)$-manifold is a homotopy sphere $\Sigma$ and the trace of the surgeries performed is a
highly connected bordism from $M$ to $\Sigma$.  
\end{proof}

\begin{Corollary} \label{corhcbhb}
Let $P$ be a highly connected $(8j-1)$-manifold.  There is a unique homotopy sphere $\Sigma_P$ such that
$P \# \Sigma_P$ is the boundary of a handlebody and $s(\Sigma_P) = 0 \in bP_{8j-1}$.  Moreover, if $P_0$
and $P_1$ are two highly connected $(8j-1)$-manifolds then $\Sigma_{P_0 \# P_1} = \Sigma_{P_0} \#
\Sigma_{P_1}$.
\end{Corollary} 

\begin{proof}
We note first that $8j-1$ is a dimension to which both Lemmas \ref{lemmahcb1} and \ref{lemmahcb2}
apply.  By Lemma \ref{lemmahcb2}, $-P$ is highly connected bordant to a homotopy sphere $\Sigma$
which is determined up to connected sum with homotopy spheres in $bP_{8j}$.  There is thus a unique
homotopy sphere $\Sigma_P$ such that $s(\Sigma_P)=0$ and $[\Sigma_P]=-[P] \in \Omega_{8j-1}^{HC}$. 
Now by  Lemma \ref{lemmahcb1}, $[P \# \Sigma_P] =0$ and therefore $P \# \Sigma_P$ bounds a handlebody.  
To prove the second assertion we need only note that if $L_i$ are handlebodies such that $P_i \#
\Sigma_{P_i} = \del L_i$ for $i=0,1$ then 
$$(P_0 \# P_1) \# (\Sigma_{P_0} \# \Sigma_{P_1}) = \del (L_0 \natural L_1) \text{~~~and~~~}
s(\Sigma_{P_0} \# \Sigma_{P_1}) = s(\Sigma_{P_0}) + s(\Sigma_{P_1}) = 0.$$
\end{proof}

\noindent
We remind the reader that in the dimensions of primary concern in this thesis
$\Theta_7/bP_8 = 0$ and $\Theta_{15}/bP_{16} = \Z/2$.  In particular, in dimension $7$, the homotopy
sphere $\Sigma_P$ is always the standard sphere and every manifold $P$ bounds a handlebody.

\section{Decomposition of quadratic functions} \label{secdoqf}
Quadratic functions are of two essentially different types.  There are linear or ``totally degenerate''
quadratic functions at one extreme and at the other there are nondegenerate quadratic functions.  In
this section we consider the decomposition of a general quadratic function into the sum of linear and
nondegenerate quadratic functions.

Let $\kappa = \kappa(H,\lambda,\alpha)$ be a quadratic function and recall the fundamental sequence of
$\kappa$,
$$0 \lra F \lra H \stackrel{\hl}{\lra} H^* \stackrel{\pi}{\lra} G \lra 0.$$
By definition $F$ is the subspace of $H$ on which $\kappa$ is linear.  It is clear that
$F$ is a summand of $H$ for  suppose that $r.v \in F$ for some non-zero integer $r$.  This means that
$\lambda(rv,w) = r\lambda(v,w) = 0$ for every $w \in H$.  Since $r$ is non-zero we conclude that
$\lambda(v,w) = 0$ for every $w \in H$ and that $v$ belongs to $F$.  Thus $F$ is a summand of $H$ and
$H/F$ is a free group.  The bilinear form $\lambda$ descends to define a nondegenerate bilinear form on
$H/F$ in the familiar way:
$$\begin{array}{cccc}
\lambda_{/F}: & H/F \times H/F & \ra & \Z \\
& (v+F,w+F) & \mapsto & \lambda(v,w).
\end{array}$$
The function $\lambda_{/F}$ is well defined because different coset representatives differ by an element
of $F$ but $\lambda(u,v)=0$ for any $u \in F$ and any $v \in H$.  Moreover, if $\lambda{/F}(v+F,w+F) = 0$
for every $w+F \in H/F$ then $v \in F$, $v+F = 0$ and so $\lambda_{/F}$ is nondegenerate.  
Now let $\rho: H \ra H/F$ denote the canonical projection from $H$ onto $H/F$
and let
$$\Sec(\rho) := \{\Psi: H/F \ra H|\rho \circ s = Id_{H/F} \}$$
denote the set of sections of $\rho$.  The set $\Sec(\rho)$ is nonempty and there is a free and
transitive action of $\Hom(H/F,F)$ on $\Sec(\rho)$ defined by
$$(\psi \cdot \Psi)(v+F) = \Psi(v+F) + \psi(v+F)$$
where $\psi \in \Hom(H/F,F)$, $\Psi \in \Sec(\rho)$ and $v+F \in H/F$.  For every $\Psi \in \Sec(\rho)$,
we define the projections
$$
\begin{array}{cccc}
p_F(\Psi): & H & \ra & F \\
& v & \mapsto & v - \Psi(\rho(v))
\end{array} \hskip 1cm \text{and} \hskip 1cm
\begin{array}{cccc}
p(\Psi): & H & \ra & \Psi(H/F) \\
& v & \mapsto & \Psi(\rho(v))
\end{array} $$
and the subgroups 
$$\begin{array}{cc}
H(\Psi) := \Psi(H/F) &  H^*(\Psi) := p(\Psi)^*(H(\Psi)^*)\\
F^*(\Psi) : = p_F(\Psi)^*(F^*) & G(\Psi) = \pi(F^*(\Psi)) \subset G.
\end{array}$$ 
It is straightforward to check that there are direct sum decompositions
$H = F \oplus H(\Psi)$, $H^* = F^*(\Psi) \oplus H^*(\Psi)$ and $G=G(\Psi) \oplus TG$ where $TG$ denotes
the torsion subgroup of $G$.  Moreover, we show below that $H^*(\Psi) \subset H^*$ does not depend upon
$\Psi$ and so we shall henceforth identify $H(\Psi)^*$ and $H^*(\Psi)$ via $p(\Psi)^*$.  The reason
for doing this is that it is technically convenient to think of $H(\Psi)^*$ as a subgroup of $H^*$.  We
denote the restrictions of $\alpha$, $\lambda$ and $\kappa$ as follows 
$$\alpha(\Psi): = \alpha|_{H(\Psi)}, ~~~ \lambda(\Psi):=\lambda|_{H(\Psi) \times H(\Psi)}
\text{~~~and~~~}
\kappa(\Psi):=\kappa|_{H(\Psi)}.$$
With this notation $\kappa(\Psi) = \kappa(H(\Psi),\lambda(\Psi),\alpha(\Psi))$, $\alpha =
\alpha(\Psi) + p_F(\Psi)^*(\alpha|_F)$.  We also obtain the following decomposition of the fundamental
sequence of $\kappa$:
$$\begin{array}{ccccccccccc}
0 & \lra & F & \stackrel{=}{\lra} & F & \stackrel{0}{\lra} & F^*(\Psi) & \lra & G(\Psi) & \lra & 0 \\
& & & & \bigoplus & & \bigoplus & & \bigoplus \\
& & 0 & \lra & H(\Psi) & \stackrel{\hl(\Psi)}{\lra} & H^*(\Psi) & \lra & TG & \lra & 0
\end{array}$$
where $\hl(\Psi)$ denotes the adjoint of $\lambda(\Psi)$.

\begin{Proposition} \label{propdoqf}
Let $\kappa(H,\lambda,\alpha)$ be a quadratic function and $\lambda_{/F}$, $\rho$
and $Sec(\rho)$ be as above.  Then, for every section $\Psi \in \Sec(\rho)$ and for every $\psi \in
\Hom(H/F,F)$,
\begin{enumerate}
\item{$\Psi^*(\lambda) = \lambda_{/F}$,}
\item{$\Psi$ and $\rho|_{H(\Psi)}$ define inverse isometries  
$\Psi:\Psi^*(\kappa) \cong \kappa(\Psi)$ and $\rho|_{H(\Psi)}: \kappa(\Psi) \cong \Psi^*(\kappa)$,}
\item{$\kappa(\Psi)$ is nondegenerate and $\kappa = \kappa(\Psi) \oplus \alpha|_F $,}
\item{$H^*(\Psi) = \rho^*((H/F)^*)$ and $TG = H^*(\Psi)/\hl(\Psi)(H(\Psi))$,}
\item{$\alpha(\psi\cdot\Psi) = \alpha(\Psi) + \rho^*(\psi^*(\alpha|_F))$.}
\end{enumerate}
\end{Proposition}

\begin{proof}
{\em (1.)}  Let $v,w \in H$ be arbitrary and let $v' = \Psi(v+F)$ and $w' = \Psi(w+F)$ so that $v'+F =
v+F $ and $w'+F = w+F$.  We calculate that
$$\begin{array}{cl}
\Psi^*(\lambda)(v+F,w+F) & = \lambda(\Psi(v+F),\Psi(w+F)) \\
& = \lambda(v',w') \\
& = \lambda_{/F}(v+F,w+F).
\end{array}$$
\noindent {\em (2.)}
This is true by definition.  Certainly $\Psi \circ \rho|_{H(\Psi)} = Id_{H(\Psi)}$ and $\rho|_{H(\Psi)}
\circ \Psi = Id_{H/F}$.  Moreover, for $v \in H(\Psi)$ we calculate that
$$\begin{array}{cl}
\kappa(\Psi)(\Psi(v+F)) & = \kappa(\Psi(v+F)) \\
& = \Psi^*(\kappa)(v+F)
\end{array} \text{~~~and~~~}
\begin{array}{cl}
\Psi^*(\kappa)(\rho|_{H(\Psi)}(v)) & = \kappa(\Psi(\rho|_{H(\Psi)}(v)) \\
&= \kappa(v) \\
& = \kappa(\Psi)(v).
\end{array}$$

\noindent{\em (3.)}
We showed above the that $\lambda_{/F}$ is nondegenerate and so by {\em 1} we see that
$\Psi^*(\kappa)$ is nondegenerate and then by {\em 2~} we have that $\kappa(\Psi)$ is nondegenerate.  For
each $v \in H$, we decompose $v$ as $v=v-\Psi(\rho(v)) + \Psi(\rho(v))$ where $v-\Psi(\rho(v)) \in F$
and $\Psi(\rho(v))
\in H(\Psi)$.  Therefore
$$\begin{array}{cl}
\kappa(v) & = \lambda(\Psi(\rho(v)),\Psi(\rho(v))) + \alpha(\Psi(\rho(v))) + \alpha(v-\Psi(\rho(v)))\\
& = \kappa(\Psi)(\Psi(\rho(v)) + \alpha|_F(v-\Psi(\rho(v))) \\
& = (\kappa(\Psi) \oplus \alpha|_F)[(\Psi(\rho(v)),v-\Psi(\rho(v)))].
\end{array}$$
\noindent {\em (4.)}  
Since this result justifies our convention identifying $H^*(\Psi)$ and $H(\Psi)^*$ we briefly drop that
convention.  For every $\Psi$, $\rho^*:(H/F)^* \ra H^*$ factors as 
$$(H/F)^* \stackrel{(\rho|_{H(\Psi)})^*}{\lra} H(\Psi)^* \stackrel{p(\Psi)^*}{\lra} H^*(\Psi) \subset
H^*$$ as the following calculation demonstrates.
For any $z \in (H/F)^*$ and any $v \in H$
$$\begin{array}{cl}
p(\Psi)^*\circ (\rho|_{H(\Psi)})^*(z)(v) & = (\rho|_{H(\Psi)})^*(z)((\Psi \circ \rho)(v)) \\
& = z((\rho|_{H(\Psi)} \circ \Psi \circ \rho)(v)) \\
& = z(\rho(v)) \\
& = \rho^*(z)(v).
\end{array}$$
As $(\rho|_{H(\Psi)})^*$ and $p(\Psi)^*$ are isomorphisms it follows that 
$H^*(\Psi) = p(\Psi)^*(H(\Psi)^*) = \rho^*((H/F)^*)$.  The second assertion follows from the splitting
of the fundamental sequence of $\kappa$ just prior to the statement of this proposition.  Equivalently
it may demonstrated as follows.  For every $v \in H$, 
$\hl(v) = \hl(\Psi(\rho(v)))$ and for every $w \in H$,
$$\begin{array}{cl}
\lambda(v,w) & = \lambda(v - \Psi(\rho(v)) + \Psi(\rho(v)),w) \\
& = \lambda(\Psi(\rho(v)),w).
\end{array}$$
It follows that $\hl(H) = \hl(\Psi)(H(\Psi))$ and that 
$$\Cok(\hl) = F^*(\Psi) \oplus H^*(\Psi)/\hl(\Psi)(H(\Psi)) = F^*(\Psi) \oplus \Cok(\hl(\Psi)).$$  
Moreover, as $\lambda(\Psi)$ is nondegenerate, $\Cok(\hl(\Psi))$ is a torsion group and thus it is the
torsion subgroup of $\Cok(\hl)$.  

\noindent{\em (5.)}
Given $v \in H$, $\alpha(\Psi)(v) = \alpha(\Psi(\rho(v)))$ and thus
$$\begin{array}{cl}
\alpha(\psi \cdot \Psi)(v) & = \alpha(\Psi(\rho(v))+\psi(\rho(v))) \\
& = \alpha(\Psi(\rho(v))) + \alpha|_F(\psi(\rho(v))) \\
& = \alpha(\Psi)(v) + \rho^*(\psi^*(\alpha|_F))(v).
\end{array}$$

\end{proof}

\section{Linking forms and Kronecker pairings} \label{secblf}
In this section we review the algebraic structures that a quadratic function $\kappa =
\kappa(H,\lambda,\alpha)$ induces on its quotient $G$ and radical $F$.  The most important of these
is the linking form induced on the torsion subgroup of $G$.  By a linking form on a finite
abelian group $TG$ we shall mean a function
$$b: TG \times TG \ra \Q/\Z$$
with the following properties.
\begin{enumerate}
\item{Symmetry: $b(x,y) = b(y,x)$.}
\item{Bilinearity: $b(x,y+z) = b(x,y) + b(x,z)$.}
\item{Nonsingularity: $b(x,y) = 0 ~~~\forall x \in TG \Rightarrow y = 0$.}
\end{enumerate}
In Chapter \ref{chaplink} we shall describe in detail the classification of linking forms on
finite abelian groups but for now we report some fundamental facts concerning linking forms.  If
$(TG_0,b_0)$ and $(TG_1,b_1)$ are linking forms then an isomorphism $\theta:TG_0 \ra TG_1$ is an
isometry from $b_0$ to $b_1$ if $b_0 = b_1 \circ (\theta \times \theta)$.  There is a linking form $(TG_0
\oplus TG_1, b_0 \oplus b_1)$ defined by
$$(b_0 \oplus b_1)[(x_0,x_1),(y_0,y_1)] = b_0(x_0,y_0) + b_1(x_1,y_1).$$
If $(TG,b)$ is a linking form and $TG_0$ is a subgroup of $TG$ then the annihilator of $TG_0$ is the
set $$TG_0^{\perp} := \{x \in TG: b(x,y) = 0 ~\forall y \in TG_0 \}.$$  

\begin{Lemma} \label{lemmablf}
Let $(TG,b)$ be a linking form.  If $TG_0$ is a finite subgroup of $TG$ such that \\
$b|_{TG_0 \times TG_0}:TG_0 \times TG_0 \ra \Q/\Z$ is nonsingular then $b|_{TG_0^{\perp} \times
TG_0^{\perp}}$ is  nonsingular and 
$$b \cong b|_{TG_0 \oplus TG_0} \oplus b|_{TG_0^{\perp} \times TG_0^{\perp}}.$$
\end{Lemma}

\begin{proof}
We postpone the proof which is straight forward until chapter \ref{chaplink}.
\end{proof}

\noindent
A linking form $(TG,b)$ is called {\em indecomposable} if it cannot be written as a non-trivial
direct sum of sublinking forms.  In \cite{Wa2} Wall proved that the groups which arise as the domains of
indecomposable linking forms are $\Z/p^k$ and $\Z/2^k \oplus \Z/2^k$ for $p$ any prime and $k$ any
positive integer.  

Every quadratic function $\kappa(H,\lambda,\alpha)$ defines a linking form on the torsion subgroup 
of its quotient group $TG \subset G=\Cok(\hl)$ via the bilinear form $\lambda$.  We begin by
assuming that $\lambda$ is nondegenerate in which case the fundamental sequence of $\lambda$ is
$$0 \lra H \stackrel{\hat{\lambda}}{\lra} H^* \stackrel{\pi}{\lra} TG \lra 0.$$
\noindent
If we tensor this sequence with $\Q$ and $\hat{\lambda}$ with $Id_{\Q}$ then $\hat{\lambda} \tensor
Id_{\Q}$ is an isomorphism.  We let $\hat{\lambda}^{-1}$ denote the restriction of the inverse
of $\hat{\lambda} \tensor Id_{\Q}$ to $H^*$, 
$$\hat{\lambda}^{-1}:=(\hat{\lambda} \tensor
Id_{\Q})^{-1}|_{H^*}: H^* \ra H \tensor \Q.$$
The linear map $\hat{\lambda}^{-1}:H^* \ra H \tensor \Q$ defines a bilinear
pairing
$$\begin{array}{cccc}
\lambda^{-1}: & H^* \times H^* & \ra & \Q \\
& (x,y) & \mapsto & y(\hat{\lambda}^{-1}(x))
\end{array}$$
where $y(v \tensor \frac{p}{q}) = \frac{p}{q}y(v)$.  The symmetry of $\hl^{-1}$ follows immediately
from the symmetry of $\lambda$.  The pairing $\hl^{-1}$ may be defined alternatively as
follows.  For any $x,y \in H^*$ there is a non-zero $r \in \Z$ and some $v \in H$ such that $rx =
\hl(v)$.  We have the following identity since $\hat{\lambda}^{-1}(x) = \frac{1}{r}v$,
$$\lambda^{-1}(x,y) = \frac{1}{r}y(v).$$
The bilinear pairing $\lambda^{-1}$ induces a nonsingular linking form on $TG$ as follows.  Let
$[x] = \pi(x)$ and $[y]=\pi(y) \in TG.$ 
$$\begin{array}{ccccc}
b = b(\lambda):& TG \times TG & \ra & \Q/\Z \\
& ([x],[y]) & \mapsto & \lambda^{-1}(x,y) & \modu~ \Z
\end{array}$$
\noindent
The linking form $b$ is well defined in $\Q/\Z$, for suppose that $x + \hat{\lambda}(v)$ and $y +
\hat{\lambda}(w)$ are different lifts of $[x]$ and $[y]$.  Then
$$\begin{array}{rl}
\lambda^{-1}(x+\hat{\lambda}(v),y+\hat{\lambda}(w)) & = \lambda^{-1}(x,y) +
\lambda^{-1}(x,\hat{\lambda}(w)) +
\lambda^{-1}(\hat{\lambda}(v),y) \\
& ~~ + \lambda^{-1}(\hat{\lambda}(v),\hat{\lambda}(w))~ \modu ~ \Z\\  
& = \lambda^{-1}(x,y) + x(w) + y(v) + \lambda(v,w) ~ \modu ~ \Z\\ 
& = \lambda^{-1}(x,y) ~ \modu ~ \Z.
\end{array}$$
The function $b(\lambda)$ is evidently bilinear and symmetric.  It is also nonsingular as we now
show.  We may choose an indivisible $v \in H$ and an integer $r$ such that $r x = \hl(v)$ for if
$v' \in H$ and $r' \in \Z$ satisfy $x = r'\hl(v')$ then if $v' = s.v$ for $s\in \Z$
and $v$ indivisible then $x = r's\hl(v)$.  Now, if  If $b([x],[y]) = 0$ for every $[y] \in TG$ then
$\frac{1}{r}y(v) \in \Z$ for every $y \in H^*$.  Since $v$ is indivisible there is a $y \in
H^*$ such that $y(v) = 1$.  Thus $r=\pm1$ and hence $[x] = 0 \in TG$.  We may also define
$b$ by its adjoint map to the torsion dual of $TG$, $\TDG = \Hom(TG,\Q/\Z)$,
$$\begin{array}{cccc}
\hat{b} : & TG & \ra & \TDG\\
& [x] & \mapsto & ([y] \mapsto b([x],[y])).
\end{array}$$
As $TG$ and $\TDG$ are finite groups of the same order $\hat{b}$ is an isomorphism if and only if $b$ is
nonsingular.

If the quadratic function $\kappa$ is degenerate then by Proposition \ref{propdoqf}, each
section $\Psi$ defines a nondegenerate quadratic function $\kappa(\Psi)$ where $TG =
H^*(\Psi)/\hl(\Psi)(H(\Psi))$ and $H^*(\Psi) = \rho^*((H/F)^*)$ are is independent of $\Psi$.  In fact
the following is true.

\begin{Lemma} \label{lemmaindeppsi}
The bilinear pairing $\lambda(\Psi)^{-1}: H^*(\Psi) \times H^*(\Psi) \ra \Q$ is independent of $\Psi$.
\end{Lemma}

\begin{proof}
Suppose that $\Psi'$ is another section of $\rho$, that $x,y \in H^*(\Psi') = H^*(\Psi)$ and
that $rx = \hl(\Psi)(v)$ for $v \in H(\Psi)$ and a non-zero integer $r$.  If we set $v':=
\Psi'(\rho(v))$, then $\hl(\Psi')(v) = rx$, $v' - v \in F$ and hence $y(v'-v) = 0$.  We
calculate that
$$\begin{array}{cll}
\lambda^{-1}(\Psi)(x,y) & = \frac{1}{r}y(v) \\
& = \frac{1}{r}y(v + v' - v) \\
& = \frac{1}{r}y(v') \\
& = \lambda^{-1}(\Psi')(x,y). 
\end{array}$$
\end{proof}
It follows from  Lemma \ref{lemmaindeppsi} that $b(\lambda(\Psi)) = b(\lambda(\Psi'))$ and so given a
degenerate bilinear form $\lambda$, we define 
$$b(\lambda): = b(\lambda(\Psi))$$ for any section $\Psi$.

In addition to the linking form $b(\lambda)$, a bilinear form $\lambda$ also defines a Kronecker 
pairing between $F = \Ker(\hl)$ and $G=\Cok(\hl)$ by
$$\begin{array}{cccc}
\bar{\lambda}:& F \times G & \ra & \Z \\
& (v,[x]) & \mapsto & x(v).
\end{array}$$
The pairing $\bl$ is well defined because any two choices of lifts of $[x]$ differ by $\hl(w)$ for
some $w \in H$ and $\hl(w)(v) = 0$.  We label the adjoint homomorphisms associated to $\bl$ by $\bl_F: F
\ra G^*$ and $\bl_G:G \ra F^*$ where $\bl_F(v)([x]) = \bl(v,[x])$ and $\bl_G([x])(v) = \bl(v,[x])$. 
Here we have extended the $^*$ notation to any abelian group $G$, so that $G^* := \Hom(G,\Z)$ denotes 
the group of homomorphism from $G$ to $\Z$.  Recall the canonical homomorphism 
$$\begin{array}{cccc}
\tau_G: & G & \ra & G^{**}\\
& g & \mapsto & (x \mapsto x(g)).
\end{array}$$

\begin{Lemma} \label{lemmaadjs}
The composition $\bl_F^* \circ \tau_G: G \ra F^*$ is equal to $\bl_G$. The homomorphisms and 
$\bl_G: G \ra F^*$ and $\bl_F:F \ra G^*$ are respectively surjective and bijective.  
\end{Lemma}

\begin{proof}
The first assertion is a simple matter to check.  Let $[x] \in G$ and $v \in F$,
$$\begin{array}{cl}
\bl_F^*(\tau_G([x]))(v) & = \tau_G([x])(\bl_F(v)) \\
& = \bl_F(v)([x]) \\
& = x(v) \\
& = \bl_G([x])(v).
\end{array}$$
In Section \ref{secdoqf} we noted that choosing a section $\Psi \in \Sec(\rho)$ defines a splitting 
$H^* = F^*(\Psi) \oplus H^*(\Psi)$ and the inclusion $p_F(\Psi)^*:F^* \cong F^*(\Psi) \hra
H^*$.  Given $x \in F^*$ it is a matter of definitions that $\bl_G(\pi(p_F(\Psi)^*(x))) = x$
and therefore $\bl_G$ is onto.  Since $\bl_G = \bl_F^* \circ \tau_G$, it follows that $\bl_F^*$ in onto. 
Now note that $F$ and $G$ have the same rank by the Rank-Nullity Theorem applied to $\hl$.  Thus
$G^{**}$ and $F^*$ have the same rank and so $\bl_F^*$ is a surjective homomorphism between free
abelian groups of the same rank and thus an isomorphism.  This entails that $\bl_F$ is itself an
isomorphism.
\end{proof}

We now review and relate the linking form and Kronecker pairing of a bilinear form from a slightly more
abstract perspective.  Let $A$ and $B$ be free abelian groups and $f:A \ra B$ a homomorphism with kernel
$F$ and cokernel $G$.  Again we let $\rho: A \ra A/F$ denote the canonical projection.  Then there is a
map  $\hat{b}(f):\rho^*((A/F)^*) \ra \TDG$ constructed in a similar manner as before.  If
$[b] \in TG$ then $[b] =\pi(b)$ for some $b \in B$ and $rb = f(a)$ for some
non-zero integer $r$ and some $a \in A$.  For $x \in \rho^*((A/F)^*)$, $x(a) = 0$ for every $a \in F$,
and so we define may 
$$\begin{array}{cccc}
\hat{b}(f): & \rho^*((A/F)^*) & \ra & \TDG\\
& x & \mapsto & ([b] \mapsto \frac{1}{r}x(a) \in \Q/\Z).
\end{array}$$
Now choose any section $\Psi: A/F \ra A$ so that we have the decomposition 
$A^* = \rho^*((A/F)^*) \oplus F^*(\Psi)$.  With respect to this splitting define $\hat{\pi}: A^* \ra F^*
\oplus \TDG$ by $\hat{\pi}(x,y) = (y|_F,\hat{b}(x))$.  It is easy to check that we  now have the
following pair of exact sequences:
$$
\begin{diagram}
\divide\dgARROWLENGTH by2
\node{0} \arrow{e}
\node{F} \arrow{e}
 \node{A} \arrow{e,t}{f} 
 \node{B} \arrow{e,t}{\pi}
  \node{G} \arrow{e} \node{0} \\
\node{0} \arrow{e}
 \node{G^*} \arrow{e,t}{\pi^*}
  \node{B^*} \arrow{e,t}{f^*}
   \node{A^*} \arrow{e,t}{\hat{\pi}}
    \node{F^* \oplus \TDG} \arrow{e} \node{0}
\end{diagram}
$$
Observe also that $\hat{b}(f)$ is identically zero on the image of $f$ and
thus defines a function which we also denote $\hat{b}(f):(A/F)^*/f(A) \ra \TDG$.  If we
return to the case where $B=A^*$ and $f=\hl$ is the adjoint of a bilinear form on $A$ then
$\hat{b}(\hl)$ is the adjoint map of the linking form $b(\lambda)$ which we defined above.    
Setting $G(\Psi) := \pi(F^*(\Psi)) \subset G$ the above pair of sequences becomes
linked into the following commutative diagram where each vertical arrow denotes and isomorphism.
$$
\begin{diagram}
\divide\dgARROWLENGTH by2
\node{0} \arrow{e}
\node{F} \arrow{e} \arrow[2]{s,l}{\bl_F}
 \node{A} \arrow{e,t}{\hat{\lambda}} \arrow[2]{s,l}{\tau_A}
 \node{A^*} \arrow{e,t}{\pi} \arrow[2]{s,l}{Id}
  \node{G(\Psi) \oplus TG} \arrow{e} \arrow[2]{s,r}{\M}
   \node{0} \\ \\
\node{0} \arrow{e}
 \node{G^*} \arrow{e,t}{\pi^*}
  \node{A^{**}} \arrow{e,t}{\hat{\lambda}^*}
   \node{A^*} \arrow{e,t}{\hat{\pi}}
    \node{F^* \oplus \TDG} \arrow{e} \node{0}
\end{diagram}
$$

Recall that if a handlebody $L$ has quadratic function $\kappa = \kappa(H,\lambda)$ and
boundary $P$ then the topological correlates of $F=\Ker(\hl)$ and $G=\Cok(\hl)$ are $H_{4j}(P)$ and
$H^{4j}(P)$ respectively.  The pairing $\bl$ corresponds to the Kronecker pairing $H_{4j}(P) \times
H^{4j}(P) \ra
\Z$ and this gives an alternate explanation of Lemma \ref{lemmaadjs}.  The linking form $b(\lambda)$ is
the linking form of
$P$ which is a defined on the torsion subgroup of $H^{4j}(P)$.  Recall also that for every manifold $P$
there is a specific homotopy $\Sigma_P$ such that $P \#
\Sigma_P$ bounds a handlebody, $L$ say.  Now there is a canonical homeomorphism $P \ra P \# \Sigma_P$
and therefore a canonical identification of $H^{4j}(P)$ with $H^{4j}(P \# \Sigma_P)$. 

\begin{Definition} \label{defb}
Let $P$ be a highly connected $(8j-1)$-manifold and $L$ a handlebody cobounding $P \# \Sigma_P$
with quadratic function $\kappa(L) = \kappa(\lambda)$.  The linking form of $P$  is defined to be the
linking form induced by $\lambda$;
$$b(P) := b(\lambda).$$
\end{Definition}

\begin{remark}
The linking form of $P$ has at least two other definitions: by way of Poincar\'{e} duality with $\Q/\Z$
coefficients (\cite{KM}) and by letting $r$ times a cycle representing a homology class bound and
computing intersection $(\modu~r)$ (\cite{Wa3}).  We leave it as an exercise to verify the well known
fact that all three definitions give the same linking form.  That our definition of $b(P)$ is
independent of the choice of handlebody will follow from results of Wilkens' thesis which we review in
Section \ref{secwilkc}.  
\end{remark}

Now let $\lambda_0$ and $\lambda_1$ be a pair of bilinear forms.  It is clear from the definitions
that $b(\lambda_0 \oplus \lambda_1) = b(\lambda_0) \oplus b(\lambda_1)$.  We claim for two
highly connected manifolds $P_0$ and $P_1$ that 
$$b(P_0 \# P_1) = b(P_0) \oplus b(P_1).$$  
For firstly by Corollary \ref{corhcbhb}, $\Sigma_{P_0 \# P_1} = \Sigma_{P_0} \# \Sigma_{P_1}$ 
and thus if $\del L_i = P_i \# \Sigma_{P_i}$ for $i=0,1$, then $(P_0 \# P_1) \# \Sigma_{P_0 \# P_1} =
\del (L_0 \natural L_1)$.  But now $\lambda(L_0 \natural L_1) = \lambda(L_0) \oplus \lambda(L_1)$ and
our claim follows.

\section{Wilkens' classification of the manifolds $P$} \label{secwilkc}
The fundamental topological result in Wilkens' thesis is the following theorem (whose proof we shall
present in Chapter \ref{chaptop}).

\begin{Theorem}[Wilkens' Thesis 3.2]
\label{thmwilk3.2}
Let $L_0$ and $L_1$ be $2k$-dimensional handlebodies $(k>2)$ with boundaries $\del L_0 = P_0$ and
$\del L_1 = P_1$.  Suppose that $f:P_0 \ra P_1$ is a diffeomorphism.  Then there
are handlebodies $V_0, V_1$ with boundaries the standard sphere and a diffeomorphism $g$
which makes the following diagram commute.
\[
\divide\dgARROWLENGTH by 2
\begin{diagram}
\node{P_0} \arrow{e,t}{f} \arrow{s,j} 
 \node{P_1} \arrow{s,J}  \\
\node{L_0 \natural V_0} \arrow{e,t}{g}
 \node{L_1 \natural V_1}
\end{diagram}
\]
\end{Theorem}

\noindent
In the light of Theorem \ref{thmwilk3.2} and recalling Definition \ref{defstabfequiv} which defined the
stable $\FF^*$-equivalence class of a quadratic function, we make the following

\begin{Definition} \label{defsqf}
Let $P$ be a highly connected $(8j-1)$-manifold and $L$ a handlebody such that $\del L =  P \# \Sigma_P$.
We define the stable quadratic function of $P$, $[\kappa(P)]$, to be the stable $\FF^*$-equivalence 
class of $\kappa(L)$.
\end{Definition}

\begin{Corollary}[To Wilkens' \ref{thmwilk3.2} and Wall's \ref{thmwallrs}] \label{corwilkwall}
The stable quadratic function of a highly connected manifold $P$, $[\kappa(P)]$, is a well defined and
complete almost diffeomorphism invariant of $P$.
\end{Corollary}

\begin{proof}
We show that $[\kappa(P)]$ is both well defined and invariant under almost diffeomorphisms at the same
time.  Let $P_0$ and $P_1$ be highly connected manifolds, $\Sigma$ a homotopy sphere and $f:P_0 \ra P_1
\# \Sigma$ a diffeomorphism.  Then $\Sigma_{P_0} = \Sigma_{P_1 \#
\Sigma} = \Sigma_{P_1} \# \Sigma_{\Sigma}$ and there is a diffeomorphism
$$f':P_0 \# \Sigma_{P_0} \ra P_1 \# \Sigma \# \Sigma_{P_0} = P_1 \# \Sigma_{P_1} \# \Sigma_{\Sigma} \#
\Sigma.$$
Let $L_0$, $L_1$ and $L_{\Sigma}$ be handlebodies bounding respectively $P_0 \# \Sigma_{P_0}$, $P_1 \#
\Sigma_{P_1}$ and $\Sigma \# \Sigma_{\Sigma}$.  Applying Theorem \ref{thmwilk3.2} to $f'$, $L_0$ and $L_1
\natural L_{\Sigma}$ we deduce that there are handle bodies $V_0$ and $V_1$ bounding standard spheres and a
diffeomorphism
$$g: L_0 \natural V_0 \ra L_1 \natural L_{\Sigma} \natural V_1$$
which extends $f'$.  So by Wall's classification of handlebodies there is an isometry
$$\Theta(g): \kappa(L_0) \oplus \kappa(V_0) \ra \kappa(L_1) \oplus \kappa(L_{\Sigma}) \oplus 
\kappa(V_1).$$
But the boundaries of $V_0$, $V_1$ and $L_{\Sigma}$ are all homotopy spheres, so
$\kappa(V_0)$, $\kappa(V_1)$ and $\kappa(L_{\Sigma})$ are all nonsingular quadratic functions 
in $\FF^*$ and thus $\Theta(g)$ defines a stable $\FF^*$-equivalence from $\kappa(L_0)$ to
$\kappa(L_1)$.  That is, $[\kappa(P_0)] = [\kappa(L_0)]_{\FF^*} = [\kappa(L_1)]_{\FF^*} =
[\kappa(P_1)]$.  If we take
$P_0 = P_1 = P$, $\Sigma = S^{8j-1}$ and $f=Id_P$ then this argument shows that $[\kappa(P)]$ is
independent of the handlebody chosen and hence well defined.

Finally we show that $[\kappa(P)]$ is a complete invariant.  Suppose that $[\kappa(P_0)] =
[\kappa(P_1)]$.  Then there are nonsingular quadratic functions $\mu_0$ and $\mu_1$ belonging to $\FF^*$
and an isometry
$$\Theta: \kappa(L_0) \oplus \mu_0 \ra \kappa(L_1) \oplus \mu_1.$$
Again by Wall's classification of handlebodies there are handlebodies $L(\mu_0)$ and $L(\mu_1)$ which
cobound homotopy spheres $\Sigma_0$ and $\Sigma_1$ and there is a diffeomorphism implementing $\Theta$
$$h(\Theta): L_0 \natural L(\mu_0) \ra L_1 \natural L(\mu_1).$$
When restricted to the boundary of $L_0 \natural L(\mu_0)$, $h(\Theta)$ yields a diffeomorphism
$$h': P_0 \# \Sigma_{P_0} \# \Sigma_0 \ra P_1 \# \Sigma_{P_1} \# \Sigma_1$$
and therefore $P_0$ and $P_1$ are almost diffeomorphic.
\end{proof}

\noindent
As Wilkens worked with triples $(H,\lambda,\alpha)$ as opposed to quadratic functions
$\kappa(H,\lambda,\alpha)$ he would not have stated Corollary \ref{corwilkwall} as it is stated here. 
However, the equivalent statement for triples $(H,\lambda,\alpha)$ was known to Wilkens and was indeed
central to his thesis.  In the language of this thesis, the algebraic task which Wilkens undertook was to
find further algebraic invariants of stable quadratic functions in $\FF^c$ which would allow him to
determine when two quadratic functions in $\FF^c$ were stably equivalent.  He was partially
successful.  Given a triple $(H,\lambda,\alpha)$ Wilkens extracted the triple $(G,b,\beta)$.  Where $G =
\Cok(\hl)$, $b$ is the linking form induced on $G$ and $\beta$ is the image of $\alpha \in H^*$ under the
canonical map $H^* \ra H^*/\hl(H) = \Cok(\hl)$.  The topological interpretation of this triple is by now
largely familiar.  We suppose that $L$ is a $8$ or $16$ dimensional handlebody with quadratic function
$\kappa(H,\lambda,\alpha)$ and with boundary $P$.  Then $G = H^{4j}(P)$, $b=b(P)$ and if $i: P \hra L$
denotes the inclusion then $\beta = i^*(\alpha)$ is the tangential invariant of $P$ which is the primary
obstruction to the stable triviality of the tangent bundle of $P$.  

Wilkens' first achievement was to show that the manifolds $P$ split via connected sum according to 
the splitting of their linking forms.  To be precise, Wilkens defined
$P$ to be {\em decomposable} if $P$ can be expressed as a connected sum of manifolds neither of which
is a homotopy sphere and he defined $P$ to be {\em indecomposable} if is it not decomposable.

\begin{Proposition} [Wilkens Thesis] \label{propwilkc}
A manifold $P$ is indecomposable if and only if either $H^{4j}(P) \cong \Z$ or $H^{4j}(P)$ is finite and
the linking form of $P$ is indecomposable.  To every decomposition of $H^{4j}(P)$ as a sum of infinite
cyclic groups and finite groups over which the linking form of $P$ is indecomposable there corresponds a
decomposition of $P$ by indecomposable manifolds.
\end{Proposition}

\noindent
Wilkens' second achievement was to largely classify the indecomposable manifolds.

\begin{Theorem} [Wilkens' Thesis] \label{thmwilkc}
Let $P$ be indecomposable and let $P$ have invariants $(G,b,\beta)$.  If $G$ contains no $2$-torsion
then the triple
$(G,b,\beta)$ is a complete almost diffeomorphism invariant of $P$.  If $G$ contains $2$-torsion then
there are at most two almost diffeomorphism classes of manifolds with the same invariants.  A triple
$(G,b,\beta)$ with $b$ indecomposable is ambiguous (i.e. is realized by a pair of not almost
diffeomorphic manifolds) if and only if $|G| \leq 4$ and $\beta = 0$ or if
$|G|>4$ and $\beta$ is not divisible by $4$.
\end{Theorem}

\section{Quadratic linking functions}  \label{secqlf}
In this section we define the quadratic linking functions which are the invariants
required to settle the ambiguity in Wilkens classification.  We also prove elementary lemmas concerning
$\QL(b)$, the set of quadratic linking functions which refine a given linking form $b$, and then describe
how a nondegenerate quadratic function induces a quadratic linking function on its quotient group.

\begin{Definition} \label{defqlf}
A quadratic linking function $q:TG \ra \Q/\Z$ on a finite abelian group $TG$ is
a function such that the function $b(q): TG \times TG \ra \Q/\Z$ defined by
$$b(q)(x,y) := q(x+y) - q(x) - q(y) \in \Q/\Z $$ 
is a linking form on $TG$.  We say that $q$ is a quadratic refinement
of $b(q)$.  If $q$ satisfies $q(x) = q(-x)$ for all $x \in TG$ then we say that $q$
is homogeneous and we call homogeneous quadratic linking functions quadratic
linking forms.
\end{Definition}

\begin{Lemma} \label{lemmaqlfu1}
Let $(TG,b)$ be a bilinear form.  Then 
\begin{enumerate}
\item{there is a quadratic linking form $q:TG \ra \Q/\Z$ which refines $b$,}  
\item{if the order of $G$ is odd then $q$ is unique.}  
\end{enumerate}
\end{Lemma}

\begin{proof}
Choose a basis (minimal set of generators each of order $p^k$, $p$ a prime and $k$ a
positive integer) 
$\{ [x_1] \dots, [x_n] \}$ for $TG$.  If the order of $[x_i]$ is a power of $2$ then let
$q([x_i])$ be either of the members of $\Q/\Z$ such that $2q([x_i]) = b([x_i],[x_i])$.  If the
order of $[x_i]$ is odd choose the unique fraction with an odd denominator from amongst the two
choices.  The values of $q([x_i])$ and $b([x_i],[x_j])$  determine $q$ by Definition \ref{defqlf}.  The
constraint placed on the choice of $q$ for odd order elements arises because $q(nx) = n^2x$ (as
we show below).  If $[x_i]$ is of odd order then $p^k[x_i] = 0$ for $p$ an odd prime and
some positive integer $k$.  Thus $0=p^{2k}q([x_i])$ and hence the denominator of $q([x_i])$ must be
odd.  Now we show $q(nx) = n^2q(x)$ since firstly
$$\begin{array}{cl}
b([x],[x]) &  = -b([x],-[x]) \\
& = -(q(0) -q([x]) - q(-[x]))\\
& = 2q([x]).
\end{array}$$
From this it follows that $q((n+1)[x]) -q(n[x]) - q([x]) = b(n[x],[x]) = 2nq([x])$ and so
by induction $q(n[x]) = n^2q([x])$.  
\end{proof}

\begin{Lemma} \label{lemmaqlfu2}
Let $b:TG \times TG \ra \Q/\Z$ be a linking form, let $\QL(b)$ be the set of all quadratic linking
functions refining $b$ and let $q \in \QL(b)$.  For any $a \in TG$  define
$$q_{a}(x) := q(x) + b(x,a) = q(x+a) - q(a).$$  Then,
\begin{enumerate}
\item{$q_{a}$ is a quadratic linking function refining $b$,}
\item{$q_{a}(x) - q_{b}(x) = b(x,a-b)$,}
\item{$TG$ acts freely and transitively on $\QL(b)$ by the action $a \cdot q = q_a$,}
\item{if $q$ is homogeneous then $q_a$ is also homogeneous if and only if $2a=0$,}
\item{if $q=q^o_a$ where $q^o$ is a quadratic linking form, then $\beta(q):=2a$ is independent of the
choice of
$q_o$ and $a$ and thus an invariant of $q$.}
\end{enumerate}
\end{Lemma}

\begin{proof}
{\em (1.)}  Let $x,y$ be elements of $TG$. 
$$\begin{array}{cl}
q_{a}(x+y) - q_{a}(x) - q_{a}(y) & = q(x+y) + b(x+y,a)\\
&  ~~- q(x) -b(x,a) - q(y) - b(y,a)\\ 
& = q(x) + q(y) + b(x,y) -q(x) -q(y)\\
& = b(x,y)\\
\end{array}$$
Thus $q_{a}$ is a refinement of $b$. \\ 
{\em (2.)}  The second statement is clear.  \\
{\em (3.)}  It is evident for $a,b \in G$ that $q_{a+b} = (q_a)_b$ and hence we have an action.
Suppose now that $q'$ is another quadratic refinement of $b$.  We must show that $q' = q_{a}$ for
some $a \in TG$ hence that the action is transitive.  Consider then the function on 
$TG$ defined by $r(x) = q'(x) - q(x)$.  Clearly
$r(0) = 0$ and 
$$\begin{array}{cl}
r(x+y) & = q'(x+y) - q(x+y)\\
& = q'(x) + q'(y) + b(x,y) - q(x) - q(y) - b(x,y)\\
& = (q'(x) - q(x)) + (q'(y) - q(y))\\
& = r(x) + r(y).
\end{array}$$ 
Thus $r$ is linear and since $b$ is nonsingular there exists an element $a \in TG$ such that
$r(x) = b(x,a)$.  So we see that
$$\begin{array}{cl}
q_a(x) & = q(x) + b(x,a)\\
& = q(x) + r(x)\\
& = q'(x).
\end{array}$$
Finally note that {\em 2~} and the nonsingularity of $b$ entail that $q_{a} = q_{b}$ if and
only if $a = b$ and hence the action is free.\\
{\em (4.)}  Assume that $q$ is homogeneous.  Then
$$\begin{array}{cl}
q_a(x) - q_a(-x) & = q(x) + b(x,a) - (q(-x) + b(-x,a)) \\
& = b(2x,a).
\end{array}$$
Now $b(2a,x) = 0$ for all $x \in G$ if and only if $2a = 0$ and hence $q_a$ is homogeneous if and
only if $2a=0$.

\noindent
{\em (5.)} By Lemma \ref{lemmaqlfu1}, there is a quadratic form $q^o$ which refines $b$ and so by {\em
3}, there is an $a$ such that $q=q_a^o$.  Now let $q^{o'}$ be another quadratic form refining $b$ and
let $a'$ be such that $q=q^{o'}_{a'}$.  By {\em 3} and {\em 4}, $q^{o'} = q^o_{\epsilon}$ where $2
\epsilon = 0$ and by {\em 3}, $a=a'+\epsilon$.  Thus $2a'=2a$ as required.

\end{proof}

We now describe how a nondegenerate quadratic function $\kappa(H,\lambda,\alpha)$ induces a
quadratic linking function on the torsion subgroup of $G=\Cok(\hat{\lambda})$.  Recall the fundamental 
sequence of $\kappa$
$$0 \lra H \stackrel{\hat{\lambda}}{\lra} H^* \stackrel{\pi}{\lra} TG \lra 0$$
and that for $x \in H^*$ we denote $\pi(x)$ by $[x]$.
  
\begin{Definition} \label{defqlf}
Let $\kappa=\kappa(H,\lambda,\alpha)$ be a nondegenerate quadratic function.  If $\kappa \in
\QF^{ev}$ then $\lambda$ is even and we define\\
$$
\begin{array}{cccl}
q^{ev}(\kappa): & TG & \ra & \Q/\Z \\
& [x] & \mapsto & \lambda^{-1}(x,x)/{2} ~ \modu ~ \Z.
\end{array}
$$
\noindent
If $\kappa \in \QF^c$ then we define 
$$
\begin{array}{cccl}
q^c({\kappa}): & TG & \ra & \Q/\Z \\
& [x] & \mapsto & (\lambda^{-1}(x,x) + \lambda^{-1} (x,\alpha))/{2} ~ \modu ~ \Z.
\end{array}
$$
\end{Definition}

\begin{remark}
I am grateful to Stephan Stolz who first suggested the definition of $q^c$ to me
and also use of Gauss sums (see Chapter \ref{chaplink}) in the classification of
quadratic functions.  
The definition of $q^{ev}$ has been know at least since
\cite{Wa3} and its properties have been thoroughly investigated by Wall
\cite{Wa3}, Durfee \cite{Du1} and \cite{Du2} and Nikulin \cite{Ni}.  Much of this thesis is devoted
to proving that the results of these three authors concerning even quadratic
forms and quadratic linking forms generalize to characteristic quadratic
functions and quadratic linking functions.
\end{remark}

\begin{Lemma} \label{lemmaqlfrb}
Let $\kappa = \kappa(H,\lambda,\alpha)$ be a nondegenerate quadratic function in $\FF^*$.
\begin{enumerate}
\item{If $\kappa \in \FF^{ev}$ then $q^{ev}(\kappa)$ is a homogeneous refinement of $b(\lambda)$.}
\item{If $\kappa \in \FF^c$ then $q^c(\kappa)$ is a quadratic refinement of $b(\lambda)$.}
\end{enumerate}
\end{Lemma}

\begin{proof}
The calculations are entirely routine.  Since definition the of $q^{ev}(\kappa)$ is well known, we verify
only that $q^c(\kappa)$ is well defined.  Thus let $\kappa \in \FF^*$, let $x \in H^*$ and let $x +
\hl(v)$ be another representative of $[x]$.  We calculate that
$$\begin{array}{l}
\lambda^{-1}(x+\hl(v), x+\hl(v) + \alpha)/2 \text~\modu ~\Z \\
\hskip 1cm = \lambda^{-1}(x,x+\alpha)/2 + [2\lambda^{-1}(x,\hl(v)) + \lambda^{-1}(\hl(v),\hl(v)) +
\lambda^{-1}(\hl(v),\alpha)]/2 ~\modu~\Z \\
\hskip 1cm  = \lambda^{-1}(x,x+\alpha)/2 + [2x(v) + \lambda(v,v) + \alpha(v)]/2 ~\modu~\Z \\
\hskip 1cm = \lambda^{-1}(x,x+\alpha)/2 ~\modu~\Z, 
\end{array}$$
where the last equality holds because $\alpha$ is characteristic for $\lambda$.

Next, given $x,y \in H^*$ and $[x], [y] \in TG$ we calculate that
$$\begin{array}{clc}
q^{ev}(\kappa)([x]+[y]) & = \lambda^{-1}(x+y,x+y)/2 & \modu~ \Z \\
& = \lambda^{-1}(x,x)/2 + \lambda^{-1}(y,y)/2 + \lambda ^{-1}(x,y) & \modu~\Z \\
& = q^{ev}(\kappa)([x]) + q^{ev}([y]) + b(\lambda)([x],[y]).\\
\end{array}$$
Hence $q^{ev}$ refines $b(\lambda)$.  Moreover it is clear that $q^{ev}(-[x]) = q^{ev}([x])$ and hence
we have proven {\em 1}.  In the characteristic case we calculate that
$$\begin{array}{clc}
q^{c}(\kappa)([x]+[y]) & = \lambda^{-1}(x+y,x+y)/2 + \lambda^{-1}(x+y,\alpha)/2 & \modu~ \Z \\
& = (\lambda^{-1}(x,x)+\lambda^{-1}(x,\alpha))/2 \\
& ~~+ (\lambda^{-1}(y,y)+\lambda^{-1}(y,\alpha))/2 + \lambda ^{-1}(x.y) & \modu~\Z \\ 
& = q^{c}(\kappa)([x]) + q^{c}([y]) + b(\lambda)([x],[y])\\
\end{array}$$
and this proves {\em 2}.
\end{proof}

Given two quadratic linking functions $q_0: TG_0 \ra \Q/\Z$ and $q_1:TG_1 \ra \Q/Z$ an isomorphism
$\theta:TG_0 \ra TG_1$ is called an isometry if $q_0 = q_1 \circ \theta$.  In this situation we shall
write that $q_0 \cong q_1$.  We consider next the question
of which quadratic linking functions are isometric to  $q^*(\kappa)$ for some $\kappa$.

\begin{Theorem}[Wall \cite{Wa3}] \label{thmwallrealise}
Let $q:TG \ra \Q/\Z$ be a quadratic linking form on a finite abelian group.  Then there is a
nondegenerate even quadratic form $\kappa(\lambda,0)$ such that $q \cong \delta^{ev}(\kappa)$.
\end{Theorem} 

\begin{Corollary} \label{corcharrealise}
\begin{enumerate}
\item{ Let $b:TG \times TG \ra \Q/\Z$ be a linking form on a finite abelian group, then there is a
bilinear form $\lambda$ such that $b \cong b(\lambda)$.}
\item{Let $q: TG \ra \Q/\Z$ be a quadratic linking function of a finite abelian group, then
there is a characteristic quadratic function $\kappa$ such that $q \cong \delta^c(\kappa)$.}
\end{enumerate}
\end{Corollary}

\begin{proof}
{\em (1.)}  By Lemma \ref{lemmaqlfu1} there is a quadratic linking form $q$ which refines $b$. 
Now apply Theorem \ref{thmwallrealise} to deduce that there is a nondegenerate quadratic form
$\kappa(\lambda,0)$ such that $q \cong \delta^{ev}(\kappa)$.  From Lemma \ref{lemmaqlfrb} we see that $b
\cong b(\lambda)$.
\\ {\em (2.)}  We begin by stating  and proving a simple lemma which will be useful elsewhere.

\begin{Lemma} \label{lemmaqlfpert}
Let $\kappa = \kappa(H,\lambda,\alpha)$ be a nondegenerate characteristic quadratic function and let
$\epsilon \in H^*$ be divisible by $2$.  Then
$$q^c(\kappa(H,\lambda,\alpha+\epsilon)) = q^c(\kappa)_{[\epsilon/2]}.$$
\end{Lemma}

\begin{proof}
Let $[x]$ be any element in $TG=\Cok(\hl)$.
$$\begin{array}{cll}
q^c(\kappa(\lambda,\alpha+\epsilon))([x]) & = \hl^{-1}(x,x+\alpha + \epsilon)/2 & \modu~\Z\\
& = \hl^{-1}(x,x+\alpha)/2 + \hl^{-1}(x,\epsilon/2) & \modu~\Z\\
& = q^c(\kappa)([x]) + b(\lambda)([x],[\epsilon/2])\\
& = q^c(\kappa)_{[\epsilon/2]}([x]).
\end{array}$$
\end{proof}

\noindent
Now let $b$ be the linking form induced by $q$ and let $q'$ be a homogeneous refinement of $b$.  
By Theorem \ref{thmwallrealise} there is an even quadratic form $\kappa(H,\lambda,0)$ and an
isomorphism $\theta: \Cok(\hl) \ra TG$ such that
$q' \circ \theta = \delta^{ev}(\kappa)$.  By Lemma \ref{lemmaqlfu2} there is an element 
$a \in TG$ such that $q = q'_{a}$.  So choose $\epsilon \in H^*$ such that $\theta([\epsilon/2]) = a$. 
As $\lambda$ is even, $\kappa(\lambda,0)$ is characteristic and we may apply Lemma \ref{lemmaqlfpert} to
$\kappa(\lambda,0)$ and $\epsilon$ to conclude that 
$\delta^c(\kappa(\lambda,\epsilon)) = q'_a  \circ \theta = q \circ \theta$.
\end{proof}

Observe that any quadratic linking function $q$ over a finite abelian group $TG$, has for invariants the
bilinear for it refines, $b(q)$, and the ``linear element'', $\beta(q) \in TG$.  We shall call the pair 
$(b(q),\beta(q))$ the ``Wilkens invariants" of $q$.  Chapter \ref{chaplink} deals with question of
when quadratic functions $q_0$ and $q_1$ with isomorphic Wilkens invariants are isometric.  (The answer
is that there is one extra invariant, the Kervaire-Arf Gauss sum invariant $K$, such that $q_0
\cong q_1$ if and only if $(b(q_0),\beta(q_0) \cong (b(q_1),\beta(q_1))$ and $K(q_0) = K(q_1)$).  We
conclude this section with the following obvious but important remark concerning $\beta(q)$.

\begin{remark} \label{remwilkinvt} 
Let $q$ be a quadratic function on a finite abelian group $TG$.  The element $\beta(q) \in TG$ is even. 
If $q=\delta^c(\kappa)$, where $\kappa$ is the nondegenerate quadratic function
$\kappa=\kappa(\lambda,\alpha)$, then $\beta(q)=[\alpha] \in \Cok(\hl) = TG$.
\end{remark}

\section{Quadratic linking families} \label{secqlfam}
In this section we define the notion of a quadratic linking family on a finitely generated abelian group
$G$ and describe how a quadratic function $\kappa(\lambda,\alpha)$ induces a quadratic linking family on
its quotient group $G=\Cok(\hl)$.

We begin with some terminology and preliminary definitions before our main definition.  An element
$\beta$ in the abelian group $G$ shall be called {\em $n$-divisible} for a natural number $n$ if $\beta
= n \gamma$ for some $\gamma \in G$.  If in addition $\gamma$ is indivisible $\beta$ shall be called {\em
precisely $n$-divisible}.  We shall call $2$-divisible elements {\em even}.   Recall the homomorphism
$\tau_G: G \ra G^{**}$ from $G$ to its double dual.  Of course the torsion subgroup of $G$, $TG$, is the
kernel of $\tau_G$.  The set of sections of $\tau_G$ is the set of homomorphisms
$$\Sec(\tau_G) := \{ \Phi: G^{**} \ra G | \tau_G \circ \Phi = Id_{G^{**}} \}.$$
There is a free and transitive action of $\Hom(G^{**},TG)$ on $\Sec(\tau_G)$ defined
for $\phi \in \Hom(G^{**},TG)$, $\Phi \in Sec(\tau_G)$ and $z \in G^{**}$ as follows:
$$(\phi \cdot \Phi) (z) = \Phi(z) - \phi(z).$$

\noindent
Now for any linking form $b$ defined on $TG$ recall from Lemma \ref{lemmaqlfu2} that $TG$ acts freely
and transitively on $\QL(b)$ which is the set of all quadratic linking functions refining $b$.  Given an
even element $\beta \in G$ we let $\Hom(G^{**},TG)$ act on $\QL(b)$ by 
$\phi \cdot q = q_{\phi(\tau_G(\beta)/2)}$.    

\begin{Definition} \label{defqlfam}
Let $G$ be a finitely generated abelian group with torsion subgroup $TG$.  A quadratic linking family on
$G$ is a triple $(G,q^*,\beta)$ where $q^*:\Sec(\tau_G) \ra \QL(b)$ is a function for some
linking form $b$ defined on $TG$ and $\beta$ is an element of $G$.
\begin{enumerate}
\item{
A characteristic quadratic linking family on $G$ is a quadratic linking family $Q^c=(G,q^c,\beta)$
with $\beta$ even and satisfying the following properties:
\begin{enumerate}
\item[(A)]{$q^c$ is equivariant with respect to the actions of $\Hom(G^{**},TG)$, that is,\\
$q^c(\phi \cdot \Phi) = q^c(\Phi)_{\phi(\tau_G(\beta)/2)}$ for every $\phi \in \Hom(G^{**},TG)$ 
and every $\Phi \in \Sec(\tau_G)$,}
\item[(B)]{$\beta = \Phi(\tau_G(\beta)) + \beta (q^c(\Phi))$ for every $\Phi \in \Sec(\tau_G)$, (recall
that $\beta(q^c(\Phi))$ is the even element in $TG$ defined in Lemma \ref{lemmaqlfu2}, {\em 5}).}
\end{enumerate}}
\item{An even quadratic linking family on $G$ is a quadratic linking family $Q^{ev}=(G,q^{ev},\beta)$
where $q^{ev}$ is a {\em constant} function whose value is a quadratic linking form refining $b$ and
$\beta$ may be any element of $G$.}
\end{enumerate}
\end{Definition}

In order to motivate the somewhat complicated definition of a characteristic quadratic linking family we
proceed now to show how a quadratic function induces a quadratic linking family on its quotient
group.  So suppose that $\kappa(H,\lambda,\alpha)$ is a possibly degenerate quadratic function with
$\pi: H^* \ra G = \Cok(\hl)$ the projection onto the quotient of $\lambda$ and with $F =
\Ker(\hl)$ the radical of $\lambda$.  In Section \ref{secdoqf} we
defined the canonical projection $\rho: H \ra H/F$, the sections of $\rho$, 
$\Sec(\rho) = \{\Psi: H/F \ra H | \rho \circ \Psi = Id_{H/F} \}$ and the projections 
$p_F(\Psi): H \ra F$.  Every section $\Psi \in \Sec(\rho)$ defines a section 
$\Phi(\Psi) \in \Sec(\tau_G)$
which is the composition
$$\Phi(\Psi):=\pi\circ p_F(\Psi)^* \circ \bl_F^*: G^{**} \stackrel{\bl_{F}^*}{\lra} F^*
\stackrel{p_F(\Psi)^*}{\lra} F^*(\Psi) \stackrel{\pi|_{F^*(\Psi)}}{\lra} G$$
where $\bl_F$ was defined in Section \ref{secblf}. 
Now by Lemma \ref{lemmaadjs} every $x \in G^*$ is of the form of $\bl_F(v)$ for some $v \in F$ and for 
$z \in G^{**}$ we calculate that
$$\begin{array}{cl}
\tau_G(\Phi(\Psi)(z))[x] & = x[\Phi(\Psi)(z)]\\
& = \bl_F(v)[(\pi|_{F^*(\Psi)} \circ p_F(\Psi)^* \circ \bl_F^*)(z)]\\
& = (p_F(\Psi)^* \circ \bl_F^*(z))[v]\\
& = \bl_F^*(z)[p_F(\Psi)(v)]\\
& = \bl_F^*(z)[v]\\
& = z[\bl_F(v)] \\
& = z(x).
\end{array}$$
Hence $\tau_G \circ \Phi(\Psi) = Id_{G^{**}}$ and $\Phi(\Psi)$ is indeed a section of $\tau_G$.  There is a
free and transitive action of $\Hom(H/F,F)$ on $\Sec(\rho)$ and for $\psi \in \Sec(\rho)$ we define
$\phi(\psi) \in \Sec(\tau_G)$ to be the composition
$$\phi(\psi): G^{**} \stackrel{\bl_F^*}{\lra} F^* \stackrel{\psi^*}{\lra} (H/F)^*
\stackrel{\pi\circ \rho^*}{\lra} TG.$$

\begin{Lemma} \label{lemmaqlfamwd1}
Consider the mappings
$$\begin{array}{cccc}
\delta_H:& \Hom(H/F,F) & \ra & \Hom(G^{**},TG) \\
& \psi & \mapsto & \phi(\psi) 
\end{array} \hskip 1cm
\begin{array}{cccc}
\delta_S:& \Sec(\rho) & \ra & \Sec(\tau_G) \\
& \Psi & \mapsto & \Phi(\Psi).
\end{array}$$
\begin{enumerate}
\item{$\delta_H$ is a surjective homomorphism,}
\item{for every $\Psi \in \Sec(\rho)$ and every $\psi \in \Hom(H/F,F)$,
$\Phi(\psi\cdot\Psi) = \phi(\psi)\cdot(\Phi(\Psi))$, }
\item{$\delta_S$ is onto and 
$\delta_S^{-1}(\Phi(\Psi)) = \{ \psi \cdot \Psi | \psi \in \Ker(\delta_H) \}.$}
\end{enumerate}
\end{Lemma}

\begin{proof}
{\em (1.)} We show first that $\delta_H$ is onto.  Let $\{ z_1, \dots, z_l \}$ be a basis for $G^{**}$. 
Every section $\phi \in \Sec(\tau_G)$ is determined by $\phi(z_1), \dots, \phi(z_l)$.  We choose 
$y_1, \dots, y_l \in (H/F)^*$ so that $\pi \circ \rho^*(y_i) = \phi(z_i)$ for each $i$ which we may do
since $\pi \circ \rho^*$ is onto $TG$.  Now as $\bl_F^*:G^{**} \ra F^*$ is an isomorphism 
$\{ \bl_F^*(z_1), \dots, \bl_F^*(z_l) \}$ is a basis for $F^*$ and we choose $\psi:H/F \ra F$ so that
$\psi^*(\bl_F^*(z_i)) = y_i$ for each $i$.  It follows that 
$$\phi(\psi)(z_i) = (\pi \circ \rho^*)(\psi^*(\bl_F^*(z_i)) = \phi(z_i)$$ for each $i$.  Thus 
$\phi(\psi)=\phi$ and so $\delta$ is onto.  

The group operation in both $\Hom(H/F,F)$ and $\Hom(G^{**},TG)$ is pointwise addition.  For
$\psi_0, \psi_1 \in \Hom(H/F,F)$,
$$\begin{array}{cl}
\phi(\psi_0+\psi_1) & = \pi \circ \rho^* \circ (\psi_0^* + \psi_1^*) \circ \bl_F^* \\
& = \pi \circ \rho^* \circ \psi_0^* \circ \bl_F^* + \pi \circ \rho^* \circ \psi_1^* \circ \bl_F^*\\
& = \phi(\psi_0) + \phi(\psi_1).
\end{array}$$

\noindent {\em (2.)}
Observe first that for any $y \in F^*$ and any $v \in H$
$$\begin{array}{cl}
p_F(\psi\cdot\Psi)^*(y)(v) & = y(v - (\Psi + \psi)(\rho(v))) \\
& = y(v-\Psi(\rho(v))) - y(\psi(\rho(v))) \\
& = p_F(\Psi)^*(y)(v) - \rho^*(\psi^*(y))(v).
\end{array}$$
We therefore calculate that for $z \in G^{**}$
$$\begin{array}{cl}
\Phi(\psi.\Psi)(z) &= (\pi|_{F(\psi.\Psi)})[p_F(\psi.\Psi)^*(\bl_F^*(z))] \\
& = \pi[p_F(\Psi)^*(\bl_F^*(z)) - \rho^*(\psi^*(\bl_F^*(z)))]\\
& = \Phi(\Psi)(z) - \phi(\psi)(z) \\
& = (\phi(\psi) \cdot \Phi(\Psi))(z).
\end{array}$$

\noindent{\em(3.)}  This follows immediately from the {\em 1}, {\em 2} and the fact that $\Hom(H/F,F)$
and $\Hom(G^{**},TG)$ act freely and transitively on $\Sec(\rho)$ and $\Sec(\tau_G)$ respectively.
\end{proof}

Now recall from proposition Proposition \ref{propdoqf} that every section $\Psi \in \Sec(\rho)$
defines  a splitting $\kappa = \alpha|_F \oplus \kappa(\Psi)$ where $\kappa(\Psi) = \kappa|_{\Psi(H)}$
is nondegenerate.  Recall also Definition \ref{defqlf} which defined, for a nondegenerate quadratic
function $\kappa'(\lambda',\alpha') \in \FF^*$, a quadratic linking function $q^*(\kappa')$ on
$\Cok(\hl')$.  

\begin{Lemma} \label{lemmaqlfamwd2}
Let $\kappa = \kappa(H,\lambda,\alpha)$ be a quadratic function and let $\beta = [\alpha] \in
\Cok(\hl)$.  Then, for every $\psi \in \Hom(H/F,F)$ and every $\Psi \in \Sec(\rho)$,
\begin{enumerate}
\item{if $\kappa \in \FF^{ev}$, $q^{ev}(\kappa(\psi\cdot\Psi)) = q^{ev}(\kappa(\Psi)),$}
\item{if $\kappa \in \FF^c$, $q^c(\kappa(\psi \cdot \Psi)) = q^c(\kappa(\Psi))_a$,
where $a = \phi(\psi)(\tau_G(\beta))/2 \in TG$.}
\end{enumerate}
\end{Lemma}

\begin{proof}
Each $q^*(\kappa(\Psi))$ is defined on $TG = \Cok(\hl(\Psi)) = H^*(\Psi)/\hl(\Psi)(H(\Psi))$ 
and by Proposition \ref{propdoqf} both $H^*(\Psi)$ and $\Cok(\hl(\Psi))=TG$ are independent of
$\Psi$.  By Lemma \ref{lemmaindeppsi} the bilinear pairing 
$\lambda(\Psi)^{-1}:H^*(\Psi) \times H^*(\Psi) \ra \Q$ is also independent of $\Psi$.
Now for $x \in H^*(\Psi)$ and $[x] =\pi(x) \in TG$,
$$q^{ev}(\kappa(\Psi))([x]) = \lambda(\Psi)^{-1}(x,x)/2 ~ \modu ~ \Z$$
and so $q^{ev}(\kappa(\Psi))$ is independent of $\Psi$ and {\em 1} is proven.\\
{\em (2.)}  By the preceding remarks the only effect that altering $\Psi$ may have on
$q^c(\kappa(\Psi))$ comes via altering the linear   part of $\kappa(\Psi)$ which is $\alpha(\Psi)$. 
Now by {\em 5} of Proposition \ref{propdoqf} $\alpha(\psi\cdot\Psi) = \alpha(\Psi) +
\rho^*(\psi^*(\alpha|_F))$.  Since $\kappa$ is characteristic $\alpha|_F$ must be even for if $v \in F$,
$0 = \lambda(v,v) = \alpha(v) ~ (\modu~2)$.  We may therefore apply Lemma \ref{lemmaqlfpert} to conclude
that
$$q^c(\kappa(\psi \cdot \Psi)) = q^c(\kappa(\lambda(\Psi),\alpha(\Psi) + \rho^*(\psi^*(\alpha|_F))) =
q^c(\kappa(\Psi))_a$$
where $a = \pi(\rho^*(\psi^*(\alpha|_F)))/2$ and $\pi$, as usual, is the projection $H^* \ra \Cok(\hl)$. 
Now since $\pi(\alpha) = \beta$, $\alpha|_F = \bl_G(\beta)$ and by Lemma \ref{lemmaadjs}, $\bl_G =
\bl_F^* \circ \tau_G$.  Therefore
$$\begin{array}{cl}
a & = \pi(\rho^*(\psi^*(\alpha|_F)))/2 \\
& = (\pi \circ \rho^* \circ \psi^* \circ \bl_F^*) (\tau_G(\beta))/2 \\
& = \phi(\psi)(\tau_G(\beta))/2.
\end{array}$$
\end{proof}

\begin{Proposition} \label{propdefqlfam}
Let $\kappa=\kappa(H,\lambda,\alpha)$ be a quadratic function in $\FF^*$ with $G=\Cok(\hl)$ and $F =
\Ker(\hl)$.  Let $\tau_G: G \ra G^{**}$, $\Sec(\tau_G)$, $\rho: H \ra H/F$ and $\Sec(\rho)$ be as above
and let $\beta := [\alpha] \in G$. If $\kappa \in \QF^*$ and $q^*(\kappa)$ is defined by
$$\begin{array}{cccc}
q^*(\kappa): & \Sec(\tau_G) & \ra & \QL(b(\lambda)) \\
& \Phi & \mapsto & q^c(\kappa(\Psi))
\end{array}$$
for any $\Psi \in \Sec(\rho)$ such that $\Phi(\Psi) = \Phi$.   Then the triple
$$\delta^*(\kappa) := (G,q^*(\kappa),\beta)$$ 
is a characteristic (resp. even) quadratic linking family on $G$ for $* = c$ (resp. $* = ev$).
\end{Proposition}

\begin{proof}
Suppose that $\kappa \in \FF^c$.  We show first the $q^c(\kappa)$ is well defined.  By Lemma
\ref{lemmaqlfrb} each $q^c(\kappa(\Psi))$ refines $b(\lambda(\Psi))$ which is independent of
$\Psi$ and equal, by definition, to $b(\lambda)$.  Secondly, by Lemma \ref{lemmaqlfamwd1}, if $\Psi$ and
$\Psi'$ are such that $\Phi(\Psi) = \Phi(\Psi')$ then $\Psi' = \psi \cdot \Psi'$ for some $\psi$ such
that $\phi(\psi) = 0$.  If follows now from {\em 2} of Lemma \ref{lemmaqlfamwd2}  that
$q^c(\kappa(\Psi')) = q^c(\kappa(\Psi))$ and so $q^c(\kappa)$ is well defined.  Now let 
$\Phi \in \Sec(\tau_G)$ and $\phi \in \Hom(G^{**},TG)$ and let
$\Psi \in \Sec(\rho)$ and $\psi \in \Hom(H/F,F)$ be such that $\Phi(\Psi) = \Phi$ and $\phi(\psi) =
\phi$. 
$$\begin{array}{cll}
q^c(\phi \cdot \Phi) & = q^c(\kappa(\psi.\Psi)) & \text{by Lemma \ref{lemmaqlfamwd1} {\em 2}} \\
& = q^c(\Psi)_{\phi(\tau_G(\beta))/2} & \text{by Lemma \ref{lemmaqlfamwd2} {\em 1}}.
\end{array}$$ 
Hence $q^c(\kappa)$ is equivariant with respect to the action of $\Hom(G^{**},TG)$.  With
$\Phi$ and  $\Psi$ as before we must show that 
$\beta = \Phi(\tau_G(\beta)) + \beta(q^c(\Phi)) $.  Recall from the proof of Lemma \ref{lemmaqlfamwd2} that
$\alpha|_F = \bl_F^*(\tau_G(\beta))$ and also from Section \ref{secdoqf} that $\alpha =
p_F(\Psi)^*(\alpha|_F) + \alpha(\Psi)$.  We calculate that
$$\begin{array}{cl}
\beta & = \pi(\alpha) \\
& = \pi(p_F(\Psi)^*(\alpha|_F))) + \pi(\alpha(\Psi)) \\
& = (\pi \circ p_F(\Psi)^* \circ \hl_F^*)(\tau_G(\beta)) + \beta(q^c(\kappa(\Psi)))\\
& = \Phi(\tau_G(\beta)) + \beta(q^c(\Phi)).
\end{array}$$
Finally, this decomposition shows that $\beta$ is even for by Remark
\ref{remwilkinvt}, $\beta(q^c(\Psi))$ is even and, as we showed in the proof of Lemma
\ref{lemmaqlfamwd2}, $\alpha|_F$ (and hence $\tau_G(\beta)$) is even.  

If $\kappa \in \FF^{ev}$ then by Lemma \ref{lemmaqlfamwd2}, $q^{ev}(\kappa(\Psi))$ is independent of
$\Psi$ and hence  $q^{ev}(\kappa)$ is a constant function.  By Lemma \ref{lemmaqlfrb} $q^{ev}(\kappa)$
is a homogeneous refinement of $b(\lambda)$ and so the triple $(G,q^{ev}(\kappa),\beta)$ defines an even
quadratic linking family on $G$.
\end{proof}

\begin{remark}
We note that we have been using the symbol $q^*$ in two different but closely related ways.  When
$\kappa= \kappa(\lambda)$ is a degenerate quadratic function the symbol $q^*(\kappa)$
denotes a function which is defined over $\Sec(\tau_G)$ (or equivalently the splittings of 
$G = \Cok(\hl)$) which values in the quadratic linking functions defined on $TG$.  When 
$\kappa = \kappa(\lambda)$ is a nondegenerate characteristic quadratic function the symbol $q^*(\kappa)$
denotes a single quadratic linking function define on $TG = \Cok(\hl)$.  Now $\kappa$ is nondegenerate
if and only if $G=TG$ is torsion.  In this case $G^{**} = \{e\}$ and $\Sec(\tau_G) = \{ 0 \}$ is a
singleton set containing the zero homomorphism from the trivial group to $TG$.  Hence our two uses of
$q^*$ coincide provided we are prepared to identify a function on a singleton set with its range. 
In the characteristic case observe also that when $G=TG$ is a torsion group $\beta = \beta(q^c(0))$ is
determined by $q^c(0)$ and so in the torsion case a characteristic quadratic linking family is
completely determined by $q^c(0)$.
\end{remark}

\begin{remark}
Note that the equivariance property (B) of characteristic quadratic linking families ensures that the
function $q^c$ is determined by its value on any single section $\Phi \in \Sec(\tau_G)$ and the
divisibility of $\tau_G(\beta)$.  In particular if $\tau_G(\beta)$ is precisely $2r$-divisible then 
a simple argument concerning the action of $\Hom(G^{**},TG)$ on $\QL(b)$ shows for any $\Phi \in
\Sec(\tau_G)$ that ${\rm Im}(q^c) = \{q^c(\Phi)_{ra}|a \in TG \}$.   
\end{remark}

\begin{remark}
When $G$ is a free group then $TG= \{e\}$, $\Sec(\tau_G) = \{ Id_G \}$ and $q^c(Id)$ is the
trivial quadratic linking function. 
\end{remark}

We now define isometries between and sums of quadratic linking functions.  Let
$Q_0 = (G_0,q^*_0,\beta_0)$ and $Q_1 = (G_1,q^*_1,\beta_1)$ be a pair of quadratic linking families.  An
isomorphism $\theta: G_0 \ra G_1$ defines a bijection 
$$\begin{array}{cccc}
\tau_{\theta}: & \Sec(\tau_{G_0}) & \ra & \Sec(\tau_{G_1})\\
& \Phi & \mapsto & \theta \circ \Phi \circ (\theta^{**})^{-1}
\end{array}$$
which is illustrated in the following commutative diagram.
\[
\begin{diagram}
\node{G_0} \arrow{e,t}{\theta} \node{G_1} \\
\node{G^{**}_0} \arrow{n,l}{\Phi} \node{G^{**}_1} \arrow{n,r}{\tau_{\theta}(\Phi)}
\arrow{w,t}{(\theta^{**})^{-1}}
\end{diagram}
\]
We define $\theta$ to be an isometry
from $Q_0$ to $Q_1$ if $\theta(\beta_0) = \beta_1$ and if for all $\Phi \in \Sec(\tau_{G_0})$,
$$q^*_0(\Phi) = q^*_1(\tau_{\theta}(\Phi)) \circ \theta_{TG_0}.$$
  
\begin{remark} \label{remqlfiso}
If If $Q_0$ and $Q_1$ are characteristic quadratic linking families then $Q_0$ and $Q_1$ are isometric
if and only if there is an isomorphism $\theta:G_0 \ra G_1$ such that $\theta(\beta_0) = \beta_1$ and
$q^c_1(\Phi_1) \circ \theta_{TG_0} = q^c_0(\Phi_0)$ for some sections $\Phi_0 \in \Sec(\tau_{G_0})$ and 
 $\Phi_1 \in \Sec(\tau_{G_1})$.  However, there is no guarantee that $\theta$ itself is an isometry
from $Q_0$ to $Q_1$.  But there is an automorphism $\theta': G_1 \cong G_1$ such that
$\theta'|_{TG_1} = Id_{TG_1}$ and such that $\theta'\circ \theta$ is an isometry from $Q_0$ to $Q_1$.
\end{remark}

Given homomorphisms $f_0:G_0 \ra F_0$ and $f_1:G_1 \ra G_0$ we let $f_0 \oplus f_1$ denote the obvious
homomorphism from $G_0 \oplus G_1$ to $F_0 \oplus F_1$.  It is clear that $\tau_{G_0 \oplus G_1} =
\tau_{G_0} \oplus \tau_{G_1}$ and that there is an inclusion
$$\begin{array}{cccc}
\Sec(\tau_{G_0}) \oplus \Sec(\tau_{G_1}): & \ra & \Sec(\tau_{G_0 \oplus G_1}) \\
(\Phi_0,\Phi_1) & \mapsto & \Phi_0 \oplus \Phi_1.
\end{array}$$
With the quadratic linking families $Q_0$ and $Q_1$ as above we define $Q_0 \oplus Q_1$ to be the
quadratic linking family $(G_0 \oplus G_1,q^*_0 \oplus q^*_1,(\beta_0,\beta_1))$.  Here 
$$q^*_0 \oplus q^*_1: \Sec(\tau_{G_0 \oplus G_1}) \ra \QL(b_0 \oplus b_1)$$
is the unique function such that $(q^*_0 \oplus q^*_1)(\Phi_0 \oplus \Phi_1) = q^*_0(\Phi_0) \oplus
q^*_1(\Phi_1)$ for every $\Phi_i \in \Sec(\tau_{G_i})$, $i=0,1$ and which is equivariant with respect to
the action of $\Hom(G_0^{**} \oplus G_1^{**},TG_0 \oplus TG_1)$.  Note that in the even case we take
this action to be trivial so that $q^{ev}_0 \oplus q^{ev}_1$ is again a constant function. 

We now show that the isometries and directs sums of quadratic functions induce isometries and direct
sums of the associated quadratic linking families.   Let $\kappa_0 = (H_0,\lambda_0,\alpha_0)$ and
$\kappa_1 = (H_1,\lambda_1,\alpha_1)$ be a pair of quadratic functions with canonical projections
$\rho_i: H_0 \ra H_0/F_0$, $i=0,1$.  For every isometry $\Theta: H_0 \ra H_1$, $\Theta(F_0) = F_1$ and
hence $\Theta$ induces an isomorphism $\Theta_{/F}:H_0/F_0 \ra H_1/F_1$.  The isometry $\Theta$ defines a
bijection 
$$\begin{array}{cccc}
\rho_{\Theta}: & \Sec(\rho_0) & \ra & \Sec(\rho_1)\\
& \Psi & \mapsto & \Theta \circ \Psi \circ (\Theta_{/F})^{-1}
\end{array}$$
which is illustrated in the following commutative diagram.
\[
\begin{diagram}
\node{H_0} \arrow{e,t}{\Theta} \node{T_1} \\
\node{H_0/F_0} \arrow{n,l}{\Psi} \node{H_1/F_1} \arrow{n,r}{\rho_{\Theta}(\Psi)}
\arrow{w,t}{(\Theta_{/F})^{-1}}
\end{diagram}
\]  

\begin{Lemma} \label{lemmaqlfam1}
Let $\Theta: H_0 \ra H_1$ be an isometry between quadratic functions 
$\kappa_0 = \kappa(H_0,\lambda_0,\alpha_0)$ and $\kappa_1 = \kappa(H_1,\lambda_1,\alpha_1)$.  Then,
\begin{enumerate}
\item{ the following diagram commutes,
\[
\begin{diagram}
\divide\dgARROWLENGTH	by2
\node{0} \arrow{e}
\node{F_0} \arrow{e} \arrow{s,l}{\Theta|_{F_0}}
 \node{H_0} \arrow{e,t}{\hat{\lambda}_0} \arrow{s,l}{\Theta}
 \node{H_0^*} \arrow{e,t}{\pi_0} \arrow{s,r}{(\Theta^{*})^{-1}}
  \node{G_0} \arrow{e} 
   \node{0} \\
\node{0} \arrow{e}
 \node{F_1} \arrow{e}
  \node{H_1} \arrow{e,t}{\hl_1}
   \node{H_1^*} \arrow{e,t}{\pi_1}
    \node{G_1} \arrow{e} \node{0}
\end{diagram}
\]}
\item{for $x \in H_0^*$ and $[x] = \pi_0(x) \in G_0$, the map
$$\begin{array}{cccc}
\delta^*(\Theta): & G_0 & \ra & G_1 \\
& [x] & \mapsto & [(\Theta^{*})^{-1}x]
\end{array}$$
is a well-defined isomorphism from $G_0$ to $G_1$,}
\item{for every $\Psi \in \Sec(\rho_0)$, $\Phi(\rho_{\Theta}(\Psi)) = \tau_{\delta^*(\Theta)}(\Phi(\Psi))
\in \Sec(\tau_{G_1})$,}
\item{for every $\Psi \in \Sec(\rho_0)$, $\Theta|_{H(\Psi)}: \kappa_0(\Psi) \ra
\kappa_1(\rho_{\Theta}(\Psi))$ is an isometry,}
\item{if $\kappa_0$ and $\kappa_1$ are nondegenerate then 
$\hl_0^{-1} = (\Theta^{-1} \tensor Id_{\Q}) \circ \hl_1^{-1} \circ (\Theta^*)^{-1}.$}
\end{enumerate}
\end{Lemma}

\begin{proof}
{\em (1.)} Let $v_0,w_0$ be elements of $H_0$.
$$\begin{array}{cl}
\Theta^*(\hat{\lambda}_1(\Theta(v_0)))(w_0) & = \hat{\lambda}_1(\Theta(v_0))(\Theta(v_1))\\
& = \lambda_1(\Theta(v_0),\Theta(v_1))\\
& = \lambda_0(v_0,w_0)\\
& = \hat{\lambda}_0(v_0)(w_0)
\end{array}$$
Thus $\hat{\lambda}_0=\Theta^* \circ \hl_1 \circ \Theta$ and so $(\Theta^{*})^{-1} \circ \hl_0 = \hl_1
\circ \Theta$ and the diagram commutes.  

\noindent
{\em (2.)} This follows from part {\em 1} and an elementary diagram chase.

\noindent
{\em (3.)}  This is a particularly complicated but routine diagram chase.  We include the commutative
diagram and leave the reader to check that it commutes.  
\[
\begin{diagram}
\node{G^{**}_0} \arrow[2]{e,t}{\bl_{F_0}^*} 
\node[2]{F_0^*} \arrow[2]{e,t}{p_{F_0}(\Psi)^*} 
\node[2]{F_0^*(\Psi)} \arrow[2]{e,t}{\pi_0|_{F_0^*(\Psi)}} \arrow{s,r}{(\Theta^*)^{-1}|_{F_0^*(\Psi)}}
\node[2]{G_0} \arrow{s,r}{\delta^*(\Theta)}\\
\node{G^{**}_1} \arrow[2]{e,t}{\bl_{F_1}^*} \arrow{n,l}{\delta^*(\Theta)^{**}} 
\node[2]{F_1^*} \arrow[2]{e,t}{p_{F_1}(\rho_{\Theta}(\Psi))^*} \arrow{n,l}{(\Theta|_{F_0})^*}
\node[2]{F_1^*(\rho_{\Theta}(\Psi))} \arrow[2]{e,t}{\pi_1|_{F_1^*(\rho_{\Theta}(\Psi))}}
\node[2]{G_1}
\end{diagram}
\]
The homomorphism $\Phi(\rho_{\Theta}(\Psi))$ is
the composition of the bottom three horizontal arrows whereas the homomorphism
$\tau_{\delta^*(\Theta)}(\Phi(\Psi))$ moves up $\delta^*(\Theta)^{**}$ along the top of the diagram and
then down $\delta^*(\Theta)$.

\noindent
{\em (4.)}  It is only necessary to observe that 
$$\Theta(H_0(\Psi)) = \Theta(\Psi(H_0/F_0))  = \rho_\Theta(\Psi)(\Theta_{/F}(H_0/F_0)) =
\rho_{\Theta}(\Psi)(H_1/F_1) = H_1(\rho_{\Theta}(\Psi)),$$
because the restriction of an isometry remains an isometry.

\noindent
{\em(5.)} If $\kappa_0$ and $\kappa_1$ are nondegenerate then $\hl_i
\tensor Id_{\Q}$ is invertible and the homomorphism $\hl_i^{-1} = (\hl_i \tensor Id_{\Q})^{-1}|_{H_i^*}$,
$i=0,1$.   We calculate that
$$\begin{array}{cl}
\hl_0^{-1} & = ((\Theta^* \circ \hl_1 \circ \Theta) \tensor Id_{\Q})^{-1}|_{H^*_0} \\
& = (\Theta^{-1} \tensor Id_{\Q}) \circ (\hl_1^{-1} \tensor Id_{\Q})|_{H_1^*} \circ 
((\Theta^*)^{-1} \tensor Id_{\Q})|_{H_0^*} \\
& = (\Theta^{-1} \tensor Id_{\Q}) \circ \hl_1^{-1} \circ (\Theta^{*})^{-1}.
\end{array}$$
\end{proof}

\begin{Lemma} \label{lemmaqlfam2}
Let $\kappa_0 = \kappa(H_0,\lambda_0,\alpha_0)$ and $\kappa_1 = \kappa(H_1,\lambda_1,\alpha_1)$ be a pair
of  quadratic functions in $\QF^*$.
\begin{enumerate}
\item{If $\Theta:\kappa_0 \ra \kappa_1$ is an isometry then
$\delta^*(\Theta):\delta^*(\kappa_0) \ra \delta^*(\kappa_1)$ is an isometry of quadratic linking
families. }
\item{$\delta^*(\kappa_0 \oplus \kappa_1) = \delta^*(\kappa_0) \oplus \delta^*(\kappa_1). $}
\end{enumerate}
\end{Lemma}

\begin{proof}  {\em (1.)}  Let $\delta^*(\kappa_0) = (G_0,q^*_0,\beta_0)$ to
$\delta^*(\kappa_1) = (G_1,q^*_1,\beta_1)$.  As $(\Theta^*)^{-1}(\alpha_0) = \alpha_1$ it follows that 
$$\delta^*(\Theta)(\beta_0) = \delta^*(\Theta)([\alpha_0]) = [\alpha_1] = \beta_1.$$  
Suppose now that $\kappa_0$ and $\kappa_1$ are nondegenerate.  Another diagram chase 
shows that $\delta^*(\Theta)$ is an isometry from $\delta^*(\kappa_0)$ to $\delta^*(\kappa_1)$.  We
include the computation for the characteristic case with $x \in H^*$ a lift of $[x] \in TG_0$.
$$\begin{array}{clcl}
q_1(\delta^c(\Theta)([x])) & = \lambda_1^{-1}((\Theta^{*})^{-1}(x),(\Theta^{*})^{-1}(x) +
\alpha_1)/2 & \modu ~ \Z\\ 
& = (\Theta^{*})^{-1}(x + \alpha_0)[\hl_1^{-1}((\Theta^{*})^{-1}(x))]/2 &  \modu ~ \Z \\
& = (\Theta^{*})^{-1}(x + \alpha_0)[(\Theta \tensor Id_{\Q})(\hl_0^{-1}(x))]/2 &  \modu ~ \Z & \text{by
Lemma \ref{lemmaqlfam2}, {\em 5}} \\ 
& = (x+\alpha_0)[\hl_0^{-1}(x)]/2 &  \modu ~ \Z \\
& = \lambda_0^{-1}(x,x+\alpha_0)/2 & \modu ~ \Z \\
& = q_0([x])
\end{array}$$
Now consider the general case.  Let $\Psi \in \Sec(\rho_0)$ and $\Phi \in
\Sec(\tau_{G_0})$ be such that $\Phi = \Phi(\Psi)$.  By part {\em 4} of Lemma \ref{lemmaqlfam1},
$\Theta|_{H(\Psi)}$ defines an isometry form $\kappa_0(\Psi)$ to $\kappa_1(\rho_{\Theta}(\Psi))$. 
Moreover, $\delta^*(\Theta)|_{TG_0} = \delta^*(\Theta|_{H(\Psi)})$ and so by the calculation for the
nondegenerate case we conclude that 
$q^*(\kappa_1(\rho_{\Theta}(\Psi))) \circ \delta^*(\Theta)|_{TG_0} = q^*(\kappa_0(\Psi))$.  We calculate
that
$$\begin{array}{cll}
q^*_1(\tau_{\delta^*(\Theta)}(\Phi)) \circ \delta^*(\Theta)|_{TG_0} 
& = q^*_1(\Phi(\rho_{\Theta}(\Psi)))
\circ \delta^*(\Theta)|_{TG_0} & \text{by Lemma \ref{lemmaqlfam1} {\em 3}}\\ 
& = q^*(\kappa_1(\rho_{\Theta}(\Psi))) \circ \delta^*(\Theta)|_{TG_0} \\
& = q^*(\kappa_0(\Psi)) \\
& = q^*_0(\Phi)
\end{array}$$
and so $\delta^*(\Theta)$ defines an isometry from $\delta^*(\kappa_0)$ to $\delta^*(\kappa_1)$.  

\noindent
{\em (2.)}  It is clear that $\Cok(\hl_0 \oplus \hl_1) = \Cok(\hl_0) \oplus \Cok(\hl_1)$ and that
$[(\alpha_0,\alpha_1)] = ([\alpha_0], [\alpha_1])$.  It is also clear that the quadratic linking
function induced by the direct sum of nondegenerate quadratic functions is the direct sum of the induced
quadratic linking functions. Now we note that for $\Psi_0 \in \Sec(\rho_0)$ and
$\Psi_1 \in \Sec(\rho_1)$ that $\Phi(\Psi_0 \oplus \Psi_1) = \Phi(\Psi_0) \oplus \Phi(\Psi_1)$.  Setting
$\Phi_0 = \Phi(\Psi_0)$ and $\Phi_1 = \Phi(\Psi_1)$ we calculate that
$$\begin{array}{cl}
q^*(\kappa_0 \oplus \kappa_1)(\Phi_0 \oplus \Phi_0) & = 
q^*(\kappa_0 \oplus \kappa_1)(\Psi_0 \oplus \Psi_1) \\
& = q^*((\kappa_0 \oplus \kappa_1)(\Psi_0 \oplus \Psi_1)) \\
& = q^*(\kappa_0(\Psi_0) \oplus \kappa_1(\Psi_1)) \\
& = q^*(\kappa_0(\Psi_0)) \oplus q^*(\kappa_1(\Psi_1)) \\
& = q^*(\kappa_0)(\Phi_0) \oplus q^*(\kappa_0)(\Phi_1)
\end{array}$$
which shows that $\delta^*(\kappa_0 \oplus \kappa_1) = \delta^*(\kappa_0) \oplus
\delta^*(\kappa_1).$
\end{proof}

\begin{remark} \label{remwrongway}
For $i=0,1$, let $L_i = L(\kappa_i)$ be handlebodies corresponding to quadratic functions
$\kappa_i$ and let $P_i = \del L_i$.  The topological correlates of $H_i$, $H_i^*$ and 
$G_i = \Cok(\hl_i)$ are respectively $H_{4j}(L_i)$, $H^{4j}(L_i)$ and $H^{4j}(P_i)$.  The isometry
$\Theta: H_0 \ra H_1$ therefore defines an isomorphism $H_{4j}(L_0) \ra H_{4j}(L_1)$ which, by Wall's
classification of handlebodies, is realized as the map induced on homology by a diffeomorphism 
$g:L_0 \ra L_1$.  If $f:=g|_{P_0}$ is restriction of $g$ to the boundary of $L_0$, then
$\delta^*(\Theta)$ corresponds to the the ``wrong way'' map $f_{!}:H^{4j}(P_0) \ra H^{4j}(P_1)$.  Recall
that $f_!$ was defined in the introduction as follows.  For $x \in H^{4j}(P_0)$, $f_{!}(x) :=
f_*(x \cap [P_0])/(P_1)$, where $[P_0] \in H_{8j}(P_0)$ and $(P_1) \in H^{8j}(P_1)$ are the orientation
classes of $P_0$ and $P_1$ and $\cap$ and $/$ denote the cap and slant products respectively.  
\end{remark}

\begin{Definition} \label{defqlfamP}
Let $P$ be a highly connected $(8j-1)$-manifold and let $L$ be a handlebody with boundary 
$P \# \Sigma_P$.  We define the quadratic linking family of $P$, $Q(P)$,
to be the quadratic linking family induced by $\kappa(L)$,
$$Q(P) := \delta^*(\kappa(L)).$$
\end{Definition}
\noindent

\begin{Lemma} \label{lemmaQPwd}
The quadratic linking family of a highly connected manifold $P$ is well defined and if $P_0$ and $P_1$
are two highly connected $(8j-1)$-manifolds then
$$Q(P_0 \# P_1) = Q(P_0) \oplus Q(P_1).$$
\end{Lemma}

\begin{proof}
Note that the proof is completely analogous to the proof of Corollary \ref{corwilkwall}.
Suppose that $L$ and $L'$ are two handlebodies with boundary $P \# \Sigma_P$.  Then by Wilkens' Theorem
\ref{thmwilk3.2} applied to $Id_{P \# \Sigma_P}$ there are handlebodies $V$ and $V'$ with boundaries
the standard sphere and a diffeomorphism
$$g: L \natural V \ra L' \natural V'$$
which extends $Id_{P \# \Sigma_P}$.  The induced map $g_*:\kappa(L) \oplus \kappa(V) \ra \kappa(L')
\oplus \kappa(V')$ is an isometry.  So by Lemma \ref{lemmaqlfam2} 
$$\delta^*(g_*): \delta^*(\kappa(L)) \oplus \delta^*(\kappa(V)) = \delta^*(\kappa(L')) \oplus
\delta^*(\kappa(V'))$$ is an isometry.  Now by the preceding remark, 
$\delta^*(g_*) = (g|_{\del( L \natural V)})_{!} = (Id_{P \# \Sigma_P})_{!} =Id$.  Moreover, both
$\delta^*(\kappa(V))$ and $\delta^*(\kappa(V'))$ are the trivial quadratic linking family and hence
$\delta^*(\kappa(L)) = \delta^*(\kappa(L'))$ and so
$Q(P)$ is well defined.  

Now suppose that $P_i \# \Sigma_i = \del L_i$ for $i=0,1$ then 
$(P_0 \# P_1) \# \Sigma_{P_0 \# P_1} = \del( L_0 \natural L_1)$ and we see that
$$\begin{array}{cl}
Q(P_0 \# P_1) & = \delta^*(\kappa(L_0 \natural L_1)) \\
& = \delta^*(\kappa(L_0) \oplus \kappa(L_1)) \\
& = \delta^*(\kappa(L_0)) \oplus \delta^*(\kappa(L_1))\\
& = Q(P_0) \oplus Q(P_1).$$
\end{array}$$
\end{proof}

\begin{remark}
In dimensions $(8j-1)$ with $j>2$, Wall gave an intrinsic definition of the quadratic form
associated to a highly connected manifold $P$ \cite{Wa3}.  At this stage I am still unable to give an
intrinsic definition of $Q(P)$ in dimensions $7$ and $15$ although this would be most desirable.
\end{remark}

\noindent
To conclude this section we summarize the algebraic invariants arising from a pair $(L,P)$ which we
have described in previous sections.

\begin{itemize}
\item{A free abelian group $H = H_{4j}(L) \cong H^{4j}(L,P)$ and its dual $H^* =
H^{4j}(L) \cong H_{4j}(L,P)$.
}
\item{A symmetric intersection pairing $\lambda: H \times H \ra \Z$ with adjoint map
$\hat{\lambda}: H \ra H^*$, $\lambda(v,w) = \lambda(w,v) = \hat{\lambda}(w)(v) =
\hat{\lambda}(v)(w)$. }
 \item{An element $\alpha \in H^*$  which is the stable tangential invariant of $L$.  When
$j=1,2$ $\alpha$ is characteristic for $\lambda$, that is $\lambda(v,v) + \alpha(v)$ is even
for all $v \in H$.  When $j>2$ there is no constraint upon $\alpha$.  }
\item{
A quadratic function $\kappa = \kappa(H,\lambda,\alpha) = \kappa(L)$ which is a synthesis of the
above invariants.}
\item{
A finitely generated abelian group $G = H^{4j}(P) \cong H_{4j-1}(P) $.    }
\item{
A nonsingular, symmetric, bilinear (linking) form $b$ on $TG$,
$ b: TG \times TG \ra \Q/\Z$.    }
\item
{An element $\beta \in G$ which is the stable tangential invariant of $P$.  When $j=2,4$ $\beta$ is
even.}
\item{A function $q^*:\Sec(\tau_G) \ra \QL(b)$ which associates a quadratic linking function to every
section of the canonical map $\tau_G: G \ra G^{**}$ (in the even case $q^{ev}$ is a constant function
with value an quadratic linking form).}
\item{A quadratic linking family, $Q(P) = (G,q^*,\beta) = \delta^*(\kappa(L))$.  }
\end{itemize}

\section{Categorical statements} \label{seccat}
Let ${\mathcal{HC}}^{8j-1}$ be the category whose objects
are closed smooth highly connected $(8j-1)$-manifolds.  An almost diffeomorphism is
a homeomorphism which fails to be a diffeomorphism at at most finitely many points (see Section
\ref{secmaness}).  The morphisms of $\HC^{8j-1}$ shall be equivalence classes of almost
diffeomorphisms $f:P_0 \ra P_1$.  Two almost diffeomorphisms $f$ and $f'$ are equivalent if
they induce the same isomorphism on cohomology, $f_{!}=f'_!:H^{4j}(P_0) \rightarrow
H^{4j}(P_1)$.  We let $\QL^c$ (resp. $\QL^{ev}$) denote the categories whose objects are 
characteristic (resp. even) quadratic linking families.  The morphisms of $\QL^*$ are
isometries of quadratic linking families.  The categories $\mathcal{HC}^{8j-1}$ and $\QL^*$ are symmetric
monoidal categories with respect to the operations of connected sum of manifolds and direct sum of
quadratic linking families respectively.

\begin{Theorem} \label{thmcat}
The functors $Q^j$ which associate to every highly connected $(8j-1)$-manifold its
quadratic linking family and to every almost diffeomorphism the induced isomorphism on
cohomology are equivalences of (symmetric monoidal) categories.
$$\begin{array}{cccc}
{Q^j} & :\mathcal{HC}^{8j-1} & \rightarrow & \QL^c \\
& \{ P,[f] \} & \mapsto & \{Q(P),f_{!} \}\\
\end{array} ~(j=1,2) \hskip .75cm
\begin{array}{cccc}
{Q^{j}} & :\mathcal{HC}^{8j-1} & \rightarrow & \QL^{ev} \\
& \{ P,[f] \} & \mapsto & \{Q(P),f_{!} \}\\
\end{array} ~~(j>2)
$$
\end{Theorem}

\begin{proof}
Let $\underline{\mathcal{HC}}^{8j-1}$ be a skeleton of $\HC^{8j-1}$.  By Theorem \ref{thmb},
$Q^j(\underline{\mathcal{HC}}^{8j-1})$ is a skeleton of $\QL^*$ and hence $Q^j$ defines an equivalence of
categories.  Moreover, for any pair of highly connected $(8j-1)$-manifolds, $P_0$ and $P_1$ and any pair
of almost diffeomorphisms $f_0$ and $f_1$, $Q(P_0 \# P_1)  = Q(P_0) \oplus Q(P_1)$ and $(f_0 \# f_1)_{!}
= (f_0)_{!} \oplus (f_1)_{!}$.
\end{proof}

Of course the equivalences of categories in the restatement of Wall's classification of
handlebodies (\ref{thmwallrs}) are related to the equivalences of categories in Theorem
\ref{thmcat}.  Recall that the functor $\kappa^j$ from Theorem \ref{thmwallrs} associates to
every $8j$-dimensional handlebody its quadratic function and to every equivalence class of
diffeomorphisms the induced map on homology.  We let $\del^j: \HH^{8j} \ra \HC^{8j-1}$ be the functor 
which associates to every handlebody its boundary and to the every equivalence class of diffeomorphisms
the equivalence class of diffeomorphisms restricted to the boundary.

\begin{Proposition} \label{propcommsq}
The functors $\del^j, \delta^*, \kappa^j$ and $Q^j$ fit into the following commuting squares. 
$$
\begin{array}{ccc}
\HH^{8j} & \stackrel{\kappa^j}{\lra} & \QF^c \\
\downarrow \del^j & & \downarrow \delta^c \\
\HC^{8j-1} & \stackrel{Q^{j}}{\lra} & \QL^{c}
\end{array} ~~(j = 1,2) \hskip 2cm
\begin{array}{ccc}
\HH^{8j} & \stackrel{\kappa^j}{\lra} & \QF^{ev} \\
\downarrow \del^j & & \downarrow \delta^{ev} \\
\HC^{8j-1} & \stackrel{Q^{j}}{\lra} & \QL^{ev}
\end{array}~~(j>2).$$
\end{Proposition}

\begin{proof}
This follows from the definition of $Q(\del L)$ since $Q(\del L) = \delta^*(\kappa(L))$.  
\end{proof}

\section{The invariants $s_1$, $\bs$ and $\mu$} \label{secs_1}
Up until this point our preliminary remarks have been preparing for the almost diffeomorphism
classification of highly connected $(8j-1)$-manifolds by their quadratic linking families.  In this
section we define the smooth invariant $s_1$ of highly connected rational homology spheres in preparation
for their smooth classification.  We also define the almost diffeomorphism invariant $\bs$ whose
algebraic analogue is crucial for the classification of quadratic linking functions.  

If a handlebody $L$ has a nondegenerate quadratic function then the boundary of $L$, $P$, is a rational
homology sphere bounding the a spin manifold in which case the Eels-Kuiper
$\mu$-invariant \cite{EK} is well defined.  In fact, Stolz \cite{St} has defined a stronger invariant
but only for homotopy spheres.  However, there is nothing which prevents the
definition of this invariant, which I shall call $s_1$, for simply connected rational homology spheres and
so we now present a generalization of Stolz's $s_1$ invariant.  Suppose then that $M$ is a
simply connected spin rational homology sphere of dimension $4k-1$ bounding a spin manifold $W$.  We 
let $\LL (W)$, $\ahat(W)$ and ${\rm ph}(W)$ denote respectively the $\LL$-genus, the $\ahat$-genus and
the Pontryagin character of $W$ which are rational power series in the Pontryagin classes $p_i$ of $W$,
\cite{Hi} and \cite{St}.  For any power series $T$ in the Pontryagin classes we let $T_k$ denote its
$4k$-dimensional component.  So for example, the first four terms of the power series of $\LL(W)$, 
$\ahat(W)$ and $\ph(W)$ are
$$\begin{array}{cccc}
\LL_0 = 1, & \LL_1 = \frac{1}{3}p_1, & \LL_2 = \frac{1}{45}(7p_2-p_1^2), & 
\LL_3 = \frac{1}{3^2 \cdot 5 \cdot 7}(62p_3 - 13p_2p_1 + 2p_1^3), \\
\ahat_0 = 1, & \ahat_1 = -\frac{2}{3}p_1, & \ahat_2 = \frac{2}{45}(-4p_2 + 7p_1^2), &
\ahat_3 = \frac{4}{3^2 \cdot 5 \cdot 7}(-16p_3 + 44p_2p_1 -31p_1^3), \\
\ph_0 = 0, & \ph_1 = p_1, & \ph_2 = \frac{1}{24}(-4p_2 + p_162), &
\ph_3 = \frac{1}{2^7 \cdot 3^ 3 \cdot 5}(144p_3-28p_2p_1 + p_1^3).
\end{array}
$$

Now following Stolz we define
$$S_k(W) = \frac{1}{8}\LL_k(W) + \frac{|bP_{4k}|}{a_k} \left(c_k\ahat_k(W) + 
(-1)^kd_k(\ahat(W){\rm ph}(W))_k\right)$$
where $a_k = 1$ for $k$ even and $2$ for $k$ odd, and $c_k$ and $d_k$ are integers such that
$$c_k{\rm num}(B_k/4k) + d_k{\rm denom}(B_k/4k) = 1.$$
The symbol $B_k$ denotes the $k$th Bernoulli number and ${\rm num}(B_k/4k)$ and ${\rm denom}(B_k/4k)$
denote respectively the numerator and denominator of the fraction $B_k/4k$ when expressed in lowest
terms.  As $M$ is a rational homology sphere $H^*(W;\Q) \cong H^*(W,M;\Q)$ and we use this isomorphism
to view each Pontryagin class $p_i(W)$ as an element of $H^*(W,M;\Q)$.  Hence we may evaluate $S_k(W)$
against on the fundamental class of $W$, $[W,\del W] \in H_{4k}(W)$, to obtain the
rational number $<S_k(W),[W,\del W]>$.

\begin{Theorem} [c.f. \cite{St} Theorem] \label{thmstolz}
Let $M$ be a $(4k-1)$-dimensional spin rational homotopy sphere bounding a spin
manifold $W$ and let $\sigma(W)$ denote the signature of $W$.   Then the coefficient of $p_k(W)$ in 
$S_k(W)$ is zero and
$$s_1(M) := \frac{1}{|bP_{4k}|} \left(<S_k(W),[W,\del W]> - \frac{1}{8}\sigma(W) \right)~~~\modu~\Z$$
is a well defined smooth invariant of $M$.  If $M = \Sigma$ is a homotopy sphere then
$$-|bP_{4k}|\cdot s_1(\Sigma) = s(\Sigma) \in \Z/|bP_{4k}|$$ 
computes Brumfiel's splitting invariant of $\Sigma$.
\end{Theorem}

\begin{proof}
The first and third assertions of the Theorem are proven by Stolz.  Moreover, Stolz's proof that the
quantity $\frac{1}{8}\sigma(W) - <S_k(W),[W,\del W> (\modu |bP_{4k}|)$ is independent of the choice of
spin coboundary $W$ makes no essential use of the fact that $M$ is a homotopy sphere as opposed to a
rational homology sphere.  We present this proof again to reassure the reader of this claim.  Suppose
that $W$ and $W'$ are two spin manifolds bounding $M$.  Let $X$ be the closed spin manifold obtained by
gluing $W$ to $-W'$ along $M$ and let $i:W \hra X$ and $i':-W' \hra X$ denote the inclusions.  As $M$ is
a rational homology sphere, $i^* \oplus i^{'*}: H^{2k}(X;\Q) \ra H^{2k}(W;\Q) \oplus H^{2k}(W'\Q)$ is an
isomorphism and $(i^* \oplus i^{'*})p_i(X) = (p_i(W),p_i(-W'))$ for $i < k$.  It follows that 
$$\sigma(X) = \sigma(W) - \sigma(W')$$
and 
$$<S_k(X),[X]> = <S_k(W),[W,\del W]> - <S_k(W'),[W',\del W']>.$$
We now calculate that
$$\begin{array}{l}
(1/{|bP_{4k}|})\left(<S_k(W),[W,\del W]> - \frac{1}{8}\sigma(W)\right) \\
\hskip 5cm - ({1}/{|bP_{4k}|}) \left(<S_k(W'),[W',\del W']> - \frac{1}{8}\sigma(W') \right) \\ 
~~~=({1}/{|bP_{4k}|}) \left(<S_k(X),[W,\del X]> \frac{1}{8}\sigma(X) \right) \\ 
~~~= ({1}/{|bP_{4k}|}) (L_k(X) - \frac{1}{8}\sigma(X)) \\ 
~~~~~+ \left( c_k(1/a_k)<\ahat(X),[X]> + (-1)^kd_k(1/a_k)
<(\ahat(X){\rm ph}(X))_k,[X]> \right) \\ 
~~~= 0 ~~~\modu ~\Z
\end{array}$$
because $\sigma(X) = <\LL_k(X),[X]>$ by Hirzebruch's signature theorem the quantities 
$(1/a_k)<\ahat(X),[X]>$ and $(1/a_k)<(\ahat(X){\rm ph}(X))_k,[X]>$ are integers by the
Hirzebruch-Riemann-Roch Theorem. 
\end{proof}

\begin{remark}
We mentioned above that Stolz's $s_1$-invariant generalizes the Eels-Kuiper $\mu$-invariant and we
now explain this remark.  For $W$ and $M$ as before Eels and Kuiper made the following definitions:
$$N_k(W) := \frac{1}{8\cdot|bP_{4k}|}\LL_k(W) + \frac{1}{a_k}\ahat_k(W)$$
and 
$$\mu(M) := <N_k(W),[W,\del W]> - \frac{1}{8\cdot|bP_{4k}|}\sigma(W).$$
Now when ${\rm num}(B_{4k}/4k) = 1$ (which is the case for $k=2$ and when $k=4$) 
then we may take $c_k = 1$ and
$d_k=0$ in the definition of
$S_k(W)$ and $s_1(M)$ and then $s_1(M) = \mu(M)$.  However, when ${\rm num}(B_{4k}/4k) \neq
1$ the $s_1$-invariant is strictly shaper than the $\mu$-invariant.  
\end{remark}

Observe that for two closed spin $(4k-1)$-manifolds $M_0$ and $M_1$ bounding spin manifolds $W_0$ and
$W_1$ that $M_0 \# M_1 = \del (W_0 \natural W_1)$ and hence $s_1(M_0 \# M_1) = s_1(M_0) + s_1(M_1)$. 
Observe also that
$$|bP_{4k}| = {\rm max} \{ {\rm denom}(s_1(\Sigma)) | \Sigma \in \Theta_{4k-1} \}$$ and then define
$$\bs(M) := |bP_{4j}| s_1(M) \in \Q/\Z.$$
It follows immediately that $\bs(\Sigma) = 0$ for any homotopy sphere and therefore that $\bs$ is an
almost diffeomorphism invariant.  In fact, when $k=2$ Kreck and Stolz show that $\bs$
is a topological invariant \cite{KS}.  We shall later identify the $\bs$-invariant of a highly connected
rational homology sphere in dimensions $7$  and $15$ with the Kervaire-Arf invariant of its quadratic
linking function.  To conclude this section we present the calculations of $s_1$ for $2$-connected
$7$-manifolds and $6$-connected $15$-manifolds.  
If $L$ is an $8j$-dimensional handlebody with tangential invariant $\alpha$ then the Pontryagin classes
of $L$ are as follows:
$$p_i(L) = \left\{ \begin{array}{cl} (2j-1))!a_k\alpha & i=j \\ 0 & i \neq j \end{array} \right\}.$$
It follows that the only possibly non-zero term in $S_{2j}(L)$ is that containing $p_j(L)^2$.  Now let
$\LL_{(j,2j)}$ and $\ahat_{(j,2j)}$ denote the coefficient of $p_j^2$ in respectively
$\LL_{2j}$ and $\ahat_{2j}$.  We also need to know the coefficients of $p_j$ in ${\rm ph}(L)_j$ and
$p_j^2$ in ${\rm ph}(L)_{2j}$ but these may be computed using Newton's formula (\cite{MS} Ex 16.A) and
the fact that all but one Pontryagin class vanish.  We calculate that
$$S_{2j}(L) = S_{(j,2j)}(L)p_j^2 = \tilde{S}_{(j,2j)}(L)\alpha^2$$
where $\tilde{S}_{(j,2j)}(L) = ((2j-1)!)^2S_{(j,2j)}(L)$ and
$$S_{(j,2j)}(L) := \frac{1}{8}L_{(j,2j)} + |bP_{8j}|\left( c_{2j}\ahat_{(j,2j)} + d_{2j}
\left( \ahat_{(j,2j)}\frac{(-1)^j}{2j!(2j-1)!} + \frac{1}{(2j-1)!} \right) \right).$$

\begin{Proposition} \label{props_1P}
Let $P$ be a highly connected manifold of dimension $7$ or $15$ and let $P \# \Sigma_P = \del L$ where
the handlebody $L$ has quadratic function $\kappa(L) = \kappa(H,\lambda,\alpha)$.  Then
$${\rm dim}(P) = 8j-1: s_1(P) = (1/|bP_{8j}|)\left( \tilde{S}_{(j,2j)} \lambda^{-1}(\alpha,\alpha) -
\frac{1}{8}\sigma(\lambda)) \right),$$
$$\begin{array}{ccc}
{\rm dim}(P) = 7: & s_1(P) = \frac{1}{28}\cdot\frac{1}{8}(\lambda^{-1}(\alpha,\alpha) -
\sigma(\lambda)),
&  \bs (P) = \frac{1}{8}(\lambda^{-1}(\alpha,\alpha) - \sigma(\lambda)), \\ 
{\rm dim}(P) = 15: & s_1(P) = \frac{1}{8,128}\cdot\frac{1}{8}(\lambda^{-1}(\alpha,\alpha) -
\sigma(\lambda)),
&  \bs (P) = \frac{1}{8}(\lambda^{-1}(\alpha,\alpha) - \sigma(\lambda)). \\
\end{array}$$
\end{Proposition}

\begin{proof}
By definition $s(\Sigma_P) = 0$.  Thus $s_1(\Sigma_P)=0$ and $s_1(P) = s_1(\del L)$.  The term
$\lambda^{-1}(\alpha,\alpha)$ appears since to compute $<\alpha^2,[L,\del L]>$ we must transport
$\alpha \in H^{4j}(L;\Q)$ to $\hl^{-1}(\alpha) \in H^{4j}(L,\del L;\Q)$ and $\hl^{-1}(\alpha)^2 =
\lambda^{-1}(\alpha,\alpha).$  The specific computations for dimensions $7$ and $15$ come from from
knowledge of $\LL_{(1,2)}$, $\LL_{(2,4)}$, $\ahat_{(1,2)}$, $\ahat_{(2,4)}$, $B_2$ and $B_4$ which may
all be found in \cite{Hi}.  Alternatively, since $B_2 = B_4 = 1/30$, $s_1 = \mu$ in dimensions $7$ and
$15$ and these calculations are in \cite{EK}. 
\end{proof}

\chapter{Stable algebra of quadratic functions} \label{chapalg} 
In this chapter we establish the algebraic results we shall need to classify the manifolds $P$ up to
almost diffeomorphism.  In Chapter \ref{chaptop} we shall apply these results.  Recall that two quadratic
functions $\kappa_0$ and
$\kappa_1$ belonging to a category $\FF^*$ are stably
$\FF^*$-equivalent if there are nonsingular quadratic functions $\mu_0$ and $\mu_1$ also belonging to
$\FF^*$ such that
$\kappa_0 \oplus \mu_0 \cong \kappa_1
\oplus \mu_1$.   In this chapter we establish that stable $\FF^*$-equivalence is closely related to
two further relations between quadratic functions.  The first of these is that of isometry between
their induced quadratic linking families and the second is orthogonality which we shall define presently.

\begin{Definition}
An isometric embedding $i_0 : \kappa_0(H_0) \hra \kappa(H)$ between nondegenerate quadratic functions is
said to be {\em primitive} if $i_0(H_0)$ is a summand of $H$.  The {\em orthogonal complement} of a
primitive embedding
$i_0$, denoted by $i_0(\kappa_0)^{\perp}$, is defined to be
$\kappa|_{i_0(H_0)^{\perp}}$ where
$$i_0(H_0)^{\perp} := \{ v \in H | \lambda(v,i_0(v_0)) = 0 ~ \forall v_0 \in H_0 \}.$$
\end{Definition}

\begin{Example}  Let $\kappa=\kappa(H,\lambda,0)$ and $\kappa_0=\kappa(H_0,\lambda,0)$ be
quadratic forms defined on $H=\Z^2[v,w]$ and $H_0=\Z[w_0]$, by the following equations:
$\lambda(v,v)=0,\lambda(v,w)=1, \lambda(w,w)=n$ and $\lambda_0(w_0,w_0) = n$,
where $n \in \Z-\{-1,0,1\}$.  Then $i_0:H_0 \hra H, w_0 \mapsto w$, is a primitive embedding of a
nondegenerate but singular form into a nonsingular form.  In this case 
$i_0(H_0)^{\perp} = \Z[nv-w]$ and $i_0(\kappa_0)^{\perp}$ is defined by 
$$i_0(\kappa_0)^{\perp}(nv-w) = n.$$
\end{Example}

\begin{Definition} \label{deforgequiv}
Let $\kappa_0=\kappa(H_0,\lambda_0,\alpha_0)$ and $\kappa_1=\kappa(H_1,\lambda_1,\alpha_1)$ be quadratic
functions belonging to the category $\FF^*$.  
\begin{enumerate}
\item{If $\kappa_0$ and $\kappa_1$ are nondegenerate then they are {\em $\FF^*$-orthogonal} if there
exists a nonsingular quadratic function $\kappa$ belonging to $\FF^*$ and primitive
embeddings $i_0 : \kappa_0 \hra \kappa$ and $i_1: \kappa_1 \hra \kappa$ such that 
$i_0(\kappa_0)^{\perp}=i_1(\kappa_1)$.}
\item{In general, $\kappa_0$ and $\kappa_1$ are {\em $\FF^*$-orthogonal} if there exits splittings 
(see Section \ref{secqlf})
$$\kappa_0 = \alpha_0|_{F_{0}} \oplus \kappa_0(\Psi_0) \hskip 1cm \text{and} \hskip 1cm
\kappa_1 = \alpha_1|_{F_{1}} \oplus \kappa_1(\Psi_1)$$
such that $\alpha_0|_{F_{0}} \cong \alpha_1|_{F_{1}}$ and such that
$\kappa_0(\Psi_0)$ is $\FF^*$-orthogonal to $\kappa_1(\Psi_1)$.}
\end{enumerate}
In both cases we write $\kappa_0 \perp_{\FF^*} \kappa_1$.
\end{Definition}

We now state the theorem whose proof is the subject of the chapter.

\begin{Theorem} \label{thmstabequivs}
Let $\kappa_0$ and $\kappa_1$ be two quadratic functions in the same family $\FF^*$.  Then the
following are equivalent.
\begin{enumerate}
\item{$\kappa_0$ is stably $\FF^*$-equivalent to $\kappa_1$, $\kappa_0 \sim_{\FF^*} \kappa_1$.}
\item{The quadratic linking families induced by $\kappa_0$ and $\kappa_1$ are
isometric, $\delta^*(\kappa_0) \cong \delta^*(\kappa_1)$.}
\item{$\kappa_0$ and $\kappa_1^-$ are $\FF^*$-orthogonal, $\kappa_0 \perp_{\FF^*} \kappa_1^{-}$.}
\end{enumerate}
\end{Theorem}

\begin{remark}
In the nondegenerate case when $\FF^* = \FF^o$, the
equivalence of {\em 1} and {\em 2} was proven by Durfee \cite{Du1}, and Wall \cite{Wa4}.  Though
he did not state it this way, Wilkens \cite{Wi1} proved the equivalence of {\em 1} and {\em 2}
for $\QF^c$ in the case where the induced quadratic linking family is indecomposable.   Nikulin proved 
the equivalence of {\em 2} and {\em 3} (\cite{Ni} Corollary 1.6.2) in the nondegenerate case where
$\FF^* = \FF^o$.  Of course, the family $\QF^c$ is of primary concern for the classification of
manifolds in dimensions $7$ and $15$, however the proofs given below may shed some light on the results
when $\FF^* = \FF^o$ or $\FF^{ev}$.  In particular, Durfee's methods use quite delicate $p$-adic
techniques, Wall's proof is somewhat opaque (to this reader) though similar to our own and Nikulin's
proofs are extremely brief.  Our perspective on these results is that they are algebraic counterparts to
the topology of handlebodies and their boundaries.  
\end{remark}

\section{Gluing and splitting nondegenerate quadratic functions} \label{secgluesplit}
The gluing lemma which follows is inspired by the following topological idea: suppose that $f:P_0
\ra P_1$ is an almost diffeomorphism of highly connected rational homology spheres with
coboundaries $L_0$ and $L_1$.  We could then form the almost smooth manifold $M=L_0 \cup_f -L_1$ which
would be a highly connected $8j$-manifold and there would be a Mayer-Vietoris decomposition for
$H_{4j}(W)$ relating the domains of the quadratic functions of $M$, $L_0$ and $L_1$.  
$$
\begin{diagram}
\divide\dgARROWLENGTH by 4
\node{0} \arrow{e} \node{H_{4j}(L_0) \oplus H_{4j}(-L_1)}
\arrow{e} \arrow{s} \node{H_{4j}(M)}
\arrow {e} \arrow{s} \node{H_{4j-1}(P)} \arrow{s} \arrow{e} \node{0} \\
\node{0} \arrow{e} \node{H_0
\oplus H_1} \arrow{e} \node{H} \arrow{e} 
\node{TG} \arrow{e} \node{0}
\end{diagram}
$$
\noindent
The idea then is to manufacture the algebraic data for a manifold $M$ from the algebraic data
corresponding to the almost diffeomorphism $f$.  If $\kappa_0=\kappa(L_0)$ and $\kappa_1=\kappa(L_1)$ 
are the quadratic functions of the handlebodies $L_0$ and $L_1$ then this data is an isomorphism
$\theta: \delta^*(L_0) \cong \delta^*(L_1)$ of the induced quadratic linking families.  In the
following lemma we show how to define a quadratic function in the class $\QF^*$ from $\theta, \kappa_0$
and $\kappa_1$ which we shall denote by $\kap$.  All of the notation follows the conventions
established in Chapter \ref{chapprelim}.  In particular the fundamental sequence of a non-degenerate
quadratic function $\kappa_i = \kappa(H_i,\lambda_i,\alpha_i)$ is
$$0 \lra H_i \stackrel{\hl_i}{\lra} H_i^* \stackrel{\pi_i}{\lra} TG_i \lra 0$$
and for $x_i \in H^*_i$ we shall write $[x_i]$ for $\pi_i(x_i)$.

\begin{Lemma}[Gluing Lemma]  \label{lemmaglue}
Let $\QF^*$ be one of the families $\QF^c$, $\QF^{ev}$ or $\QF^o$ and let $\kappa_0 =
\kappa(H_0,\lambda_0,\alpha_0)$ and $\kappa_1 = \kappa(H_1,\lambda_1,\alpha_1)$ be nondegenerate
quadratic functions in the same family $\QF^*$.  Suppose that $\theta :\Cok(\hl_0) \ra
\Cok(\hl_1)$ is an isometry of the induced quadratic linking families;
$$\delta^*(\kappa_0)  \stackrel{\theta}{\lra} \delta^*(\kappa_1).$$ 
Let $H \subset H_0^* \oplus H_1^*$ be the set $\{(x_0,x_1):\theta([x_0]) = [x_1]) \}$ and
$\lambda$ the function
$$\begin{array}{cccc}
\lambda: & H \times H & \ra & \Q \\
& [(x_0,x_1),(y_0,y_1)] & \mapsto & \lambda_0^{-1}(x_0,y_0) - \lambda_1^{-1}(x_1,y_1).
\end{array}$$
Then
\begin{enumerate}
\item[1.]{$\lambda$ is a nonsingular integral bilinear form on the free abelian group H,}
\item[2.]{there is a canonical isometric embedding $(\lambda_0 \oplus -\lambda_1) \hra \lambda$,}
\item[3.]{the element $(\alpha_0, \alpha_1) \in H$,}
\item[4.]{if $\alpha:= \hl((\alpha_0,\alpha_1))$, then the quadratic function $\kap :=
\kappa(H,\lambda,\alpha)$ is of type $\QF^*$,}
\item[5.]{there is a canonical isometric embedding $\kappa_0 \oplus \kappa_1^{-} \hra \kap$.}
\end{enumerate}
\end{Lemma}

\begin{proof}
{\em (1. \& 2.)}  Observe that $H$ is the kernel of the homomorphism $ H_0^* \oplus H_1^* \ra TG_1$ which
takes $(x_0,x_1)$ to $\theta([x_0]) - [x_1]$.  Thus $H$ is a subgroup of the free abelian group $H_0^*
\oplus H_1^*$  so is itself a free abelian group.  Now let
$(x_0,x_1)$ and $(y_0,y_1) \in H$.  We show first that
$\lambda[(x_0,x_1),(y_0,y_1)]$ is an integer.
$$\begin{array}{cll}
\lambda[(x_0,x_1),(y_0,y_1)] ~ \modu ~\Z &= 
\lambda_0^{-1}(x_0,y_0)-\lambda_1^{-1}(x_1,y_1) & ~ \modu ~ \Z \\
&= b_0([x_0],[y_0]) - b_1([x_1],[y_1]) & ~ \modu ~\Z \\
&= b_0([x_0],[y_0]) - b_1(\theta([x_0]),\theta([y_0])) & ~ \modu ~\Z \\
&= 0 & ~ \modu ~ \Z
\end{array}$$
\noindent
The last equality holds as $q^*(\kappa_0) = \theta^*q^*(\kappa_1)$ entails that $b_0  =
b_1 \circ (\theta \times \theta)$.  We next show that $\lambda$ is nonsingular.   Consider
first the inclusions 
$$\begin{array}{cccl}
i_0: & H_0 & \hra & H\\
& v_0 & \mapsto & (\hat{\lambda}_0(v_0),0)
\end{array} \hskip 2cm
\begin{array}{cccl}
i_1: & H_1 & \hra & H\\
& v_1 & \mapsto & (0,-\hat{\lambda}_1(v_1)).
\end{array}
$$
The sum of $i_0$ and $i_1$, $i_0 \oplus i_1:H_0 \oplus H_1 \hra H$ defines an embedding $(\lambda_0
\oplus -\lambda_1) \hra \lambda$ for we calculate that
$$\begin{array}{cl}
\lambda[(\hat{\lambda}_0(v_0),-\hat{\lambda}_1(v_1)),(\hat{\lambda}_0(w_0),-\hat{\lambda}_1(w_1))]
&= \lambda_0^{-1}(\hat{\lambda}_0(v_0),\hat{\lambda}_0(w_0)) \\
& \hskip 3cm - \lambda_1^{-1}(-\hat{\lambda}_1(v_1),-\hat{\lambda}_1(w_1))\\  
&= \lambda_0(v_0,w_0) - \lambda_1(v_1,w_1)\\ 
&= (\lambda_0 \oplus -\lambda_1)[(v_0,w_0),(v_1,w_1)].
\end{array}$$
 
\noindent
We now define $TG$ to be the subgroup of $TG_0 \oplus TG_1$ consisting of pairs $([x_0],\theta[x_0])$
and we let $\pi: H \ra TG$ be the map $(x_0,x_1) \ra ([x_0],[x_1])$.  The linking form $b := b(\lambda_0)
\oplus b(-\lambda_1)$ is defined on $TG_0 \oplus TG_1$ with adjoint map $\hat{b}: TG_0 \oplus TG_1 \ra
(TG_0 \oplus TG_1)^{\wedge}$.  Setting $R_{TG}: (TG_0 \oplus TG_1)^{\wedge} \ra TG^{\wedge}$ to be the
restriction map we define $\Pi: H_0^* \oplus H_1^* \ra TG^{\wedge}$ to be the composition
$R_{TG} \circ \hat{b} \circ \pi_0 \oplus \pi_1$.  It is clear that $i(H) \subset\Ker(\Pi)$ where
$i:H \ra H_0^* \oplus H_1^*$ is the inclusion $i(x_0,x_1) = (x_0,x_1)$.  Moreover, $\Pi|_{TG_0 \oplus 0}$
is onto and since the groups involved are finite we conclude that $\Ker(\Pi) = i(H)$.  We shall show that
the following diagram is commutative with exact rows and columns (where we have identified 
$(H_0 \oplus H_1) = (H_0 \oplus H_1)^{**}$.)

\begin{Diagram} \label{diag1}
\[
   \begin{diagram}
     \node[2]{0} \arrow{s} \node{0} \arrow{s} \\
     \node{0} \arrow{e} \node{H_0 \oplus H_1} \arrow{e,t}{i_0 \oplus i_1} \arrow{s,l}{i^*} 
       \node{H} \arrow{e,t}{\pi} \arrow{s,r}{i} \arrow{sw,t}{\hl}
         \node{TG} \arrow{e} \node{0} \\
    \node{0} \arrow{e}
     \node{H^*} \arrow{e,b}{i_0^* \oplus i_1^*} \arrow{s,l}{\hat{\Pi}}
						\node{H_0^* \oplus H_1^*} \arrow{e,t}{\hat{\pi}}
                         \arrow{s,r}{\Pi}
							\node{TG^{\wedge}} \arrow{e} \node{0}\\
				 \node[2]{TG^{\wedge \wedge}} \arrow{s} \node{TG^{\wedge}} \arrow{s} \\
     \node[2]{0} \node{0}
\end{diagram}
\]
\end{Diagram}

\noindent
It is evident that $i_0\oplus i_1$ is injective, that ${\rm Im}(i_0 \oplus i_1) = \Ker(\pi)$ and that
$\pi$ is onto.  It is similarly clear that $i$ is injective and we have already observed that
$\Pi$ is onto and that $\Ker(\Pi) = {\rm Im}(i)$.  Hence the first row and second column are exact.
The exactness of first column and second row follow from the duality we reviewed in section
\ref{secblf}.  We must show that $i^*:(H_0 \oplus H_1)^{**} \ra H^*$ is the same as
$\hl \circ (i_0 \oplus i_1)$.  For $v_i \in H_i$ let $v_i^{**} \in
H_i^{**}$ denote the element $\phi \mapsto \phi(v_i)$ for $\phi \in H_i^{*}$.  We compute that for
$(x_0,y_0) \in H$
$$\begin{array}{cl}
i^*(v_0^{**},v_1^{**})[(x_0,x_1)] & = (v_0^{**},v_1^{**})[(x_0,x_1)]\\
& = x_0(v_0) + x_1(v_1) \\
& = \lambda_0^{-1}(\hl(v_0),x_0) - \lambda_1^{-1}(-\hl_1(v_1),x_1)\\
& = \lambda[(\hl(v_0),-\hl(v_1)),(x_0,x_1)]\\
& = \hl[(i_0 \oplus i_1)(v_0,v_1)][(x_0,x_1)]
\end{array}$$
\noindent
and thus $i^* = \hl \circ (i_0 \oplus i_1)$ after we have identified $H_i =
H_i^{**}$ by $v_i \mapsto v_i^{**}$.  A similar computation shows that $i = (i_0^* \oplus i_1^*) \circ
\hl$.  We now require the following lemma whose proof is left to
the reader.

\begin{Lemma} \label{lemmacok}
Let $A \stackrel{f}{\ra} B \stackrel{g}{\ra} C$ be the composition of a pair of
injective homomorphisms of abelian groups.  Then there is a short exact sequence
$$0 \ra {\rm Cok}(f) \ra {\rm Cok}(g \circ f) \ra {\rm Cok} (g) \ra 0.$$
\end{Lemma}

\noindent
If we apply the lemma to $\hl\circ (i_0 \oplus i_1)$ we see that 
$$0 \ra TG^{\wedge \wedge}  \ra TG \ra \Cok(\hl) \ra 0$$
is exact, that $\Cok(\hl) = 0$ and therefore that $\lambda$ is nonsingular.\\
{\em $(3.)$}  By the definition of an isometry of a quadratic linking family 
$\theta([\alpha_0])=[\alpha_1]$ and hence $(\alpha_0,\alpha_1) \in H$.\\
{\em $(4.)$}  Suppose that $\QF^* = \QF^c$.  We must show that $\alpha = \hl((\alpha_0,\alpha_1))$ is
characteristic for $\lambda$ given that each $\alpha_i$ is characteristic for $\lambda_i$.  Let $x =
(x_0,x_1) \in H$ and observe that
$$\begin{array}{clc}
[\lambda(x,x) + \alpha(x)]/2 ~ \modu ~ \Z & = [\lambda(x,x) + \hl(\alpha_0,\alpha_1)(x)]/2 & \modu~\Z \\
& = [\lambda_0^{-1}(x_0,x_0+\alpha_0) - \lambda_1^{-1}(x_1,x_1+\alpha_1)]/2 & ~ \modu ~ \Z \\
& = q^c(\kappa_0)([x_0]) - q^c(\kappa_1)(\theta([x_0])) & ~ \modu ~ \Z \\
& = 0 & ~ \modu ~ \Z.
\end{array}$$
Thus $\hl(\alpha)$ is characteristic for $\lambda$.  If $\QF^* = \QF^{ev}$ then we must show that
$\lambda$ is even.  We again let $x = (x_0,x_1) \in H$ and compute
$$\begin{array}{clc}
\lambda(x,x)/2 ~ \modu ~ \Z & = [\lambda_0^{-1}(x_0,x_0) -
\lambda_1^{-1}(x_1,x_1)]/2 & ~ \modu ~ \Z \\ 
& = q^{ev}(\kappa_0)([x_0]) - q^{ev}(\kappa_1)(\theta([x_0])) & ~ \modu ~ \Z 
\\ & = 0 & ~ \modu ~ \Z.
\end{array}$$
Thus $\lambda$ is even.  If $\QF^* = \QF^o$ then we must show that $\lambda$ is even and that
$\alpha = 0$.  However, $\alpha = \hl((\alpha_0,\alpha_1)) = \hl((0,0)) = 0$ and 
$\lambda$ is even by the computation for the case $\QF^* = \QF^{ev}$.\\
{\em $(5.)$}  We have already proven this result for the quadratic forms $(\lambda_0 \oplus
-\lambda_1) \hra \lambda$ so we need only show that the inclusions $i_0$ and $i_1$ also commute 
with respect to $\alpha_0, \alpha_1$ and $\alpha:= \hl((\alpha_0,-\alpha_1))$.  We do this now for
$(v_0,v_1) \in H_0 \oplus H_1$.
$$\begin{array}{cl}
\alpha[(i_0 \oplus i_1)(v_0,v_1)] &= \lambda[(\alpha_0,\alpha_1),(\hl_0(v_0),-\hl_1(v_1))] \\
&= \hl^{-1}_0(\alpha_0,\hl_0(v_0)) - \hl^{-1}_1(\alpha_1,-\hl_1(v_1)) \\
&= \alpha_0(v_0) + \alpha_1(v_1) \\
\end{array}$$

\end{proof}

\begin{remark} When $\kappa_0 = \kappa_1$ and $\theta = {\rm Id}$ then
we shall simply write $\kappa_0 \cup \kappa_0^{-}$ for the ``double" of $\kappa$.  
\end{remark}

The Splitting Lemma which follows shows that every primitive embedding $i_0: \kappa_0 \hra
\kappa$ of a nondegenerate quadratic function into a nonsingular quadratic function gives rise
to a decomposition of $\kappa$ as in the Gluing Lemma: $\kappa \cong i_0(\kappa_0) \cup_{\theta}
i_0(\kappa_0)^{\perp-}$.  It is topologically motivated by considering the embedding of a handlebody
$L_0$ with boundary a rational homology sphere into a closed highly connected $8j$-manifold $L$.  The
complement of this embedding is itself a handlebody with boundary the same rational homology sphere
(though the orientation is reversed) and $L$ is obtained by gluing $L_0$ to $(L - L_0)$ along their common
boundary.

\begin{Lemma}[Splitting Lemma] \label{lemmasplit}
Let $\kappa = \kappa(H,\lambda,\alpha)$ be a nonsingular quadratic function in the family $\FF^*$
and let $H_0$ be a summand of $H$ such that $\kappa_0 := \kappa|_{H_0}$ is nondegenerate.  Then
\begin{enumerate}
\item[1.]{the group $H_1 := H_0^{\perp} = \{ v \in H: \lambda(v,v_0) = 0~\forall v_0 \in H_0 \}$ is a
summand of $H$,}
\item[2.]{the quadratic functions $\kappa_1 := (\kappa|_{H_1})^{-}$ and
$\kappa_0$ are  nondegenerate members of $\FF^*$,}
\item[3.]{there is an isomorphism $\theta_{\kappa}: \delta^*(\kappa_0)
\cong \delta^*(\kappa_1),$} 
\item[4.]{$\kappa$ and $\kappa_0 \cup_{\theta_{\kappa}} \kappa_1^{-}$ are isometric.}
\end{enumerate}
\end{Lemma}

\begin{proof}
{\em (1.)} If $rv$ belongs to $H_1$ for a nonzero integer $r$ then $rv$ pairs trivially with all of
$H_0$ but this entails that $v$ does also.  Thus $v \in H_1$ and so $H_1$ is a summand.  \\
\noindent
{\em (2.)} By definition a quadratic function is nondegenerate if and only if the associated bilinear
form is nondegenerate so we consider the bilinear forms $\lambda$, $\lambda_0 := \lambda|_{H_0 \times
H_0}$ and $\lambda_1 := -\lambda|_{H_1 \times H_1}$.  We let a subscript $\Q$ denote the forms obtained
after tensoring with the rationals so that for example $\lambda_{\Q}: (H \tensor \Q) \times (H \tensor
\Q) \ra \Q$.  As $\lambda_{0\Q}$ is a nonsingular subform of the nonsingular form $\lambda_{\Q}$ it
follows that $\lambda_{\Q} = \lambda_{0\Q} \oplus -\lambda_{1\Q}$ and therefore that $\lambda_{1\Q}$ is
nonsingular. Thus $\lambda_1$ and $\kappa_1$ are nondegenerate.  Now let $\alpha_i = \alpha|_{H_i}$,
$i=0,1$.  That both $\kappa_0 =\kappa(H_0,\lambda_0,\alpha_0)$ and
$\kappa_1=\kappa(H_1,\lambda_1,\alpha_1)$  belong to
$\FF^*$ follows from the fact that the properties of being homogeneous ($\FF^o$), having an even
quadratic form ($\FF^{ev}$) and being even ($\FF^c$) are all preserved upon restricting to a summand of
the domain of a quadratic function.  \\
\noindent
{\em (3.)}  We shall present the proof in the case where $\FF^* = \FF^c$ as the other cases are
similar and easier.  From the fact that $\lambda_{\Q} = \lambda_{0\Q} \oplus -\lambda_{1\Q}$, we deduce
that the rational vector space $H \tensor \Q = (H_0 \oplus H_1) \tensor \Q$ and thus $H_0
\oplus H_1$ is a subgroup of $H$ with finite index.  We therefore define $\pi$ to be the projection
from $H$ onto this finite group
$TG:=H/(H_0 \oplus H_1)$.  We now define a series of groups and homomorphisms which shall fit into an
extensive commutative diagram.  Let $\pi'$ be the projection from $H^*$ to the finite group
$\hl(TG):=H^*/\hl(H_0 \oplus H_1))$ which is canonically isomorphic to $TG$ via the map induced by
$\hl$.  Let
$i_0: H_0 \hra H$ and $i_1: H_1 \hra H$ denote the obvious inclusion maps and as usual let $\pi_i$
denote the projection from $H^*_i$ to $H^*_i/\hl(H_i) = TG_i$ for $i = 0,1$.  In addition define 
$J: H \ra H_0^* \oplus H_1^*$ to be the map which sends $v \in H$ to $(\hl(v)|_{H_0},\hl(v)|_{H_1})$,
$I$ to be the map which sends $\pi(v)$ to $(\pi_0 \oplus \pi_1)(J(v))$ and $K:\hl(TG) \ra TG_0 \oplus
TG_1$ be the map which sends $\pi'(x)$ to the pair $(\pi_0(i^*_0(x)),\pi_1(i^*_1(x)))$.  That $I$
and
$K$ are well defined follows from standard diagram chases in Diagram \ref{diag2}.  Note that Diagram
\ref{diag2} corresponds to Diagram \ref{diag1} of Lemma \ref{lemmaglue} and extends that diagram.  It is
commutative and all sets of collinear arrows represent exact sequences of homomorphisms.

\begin{Diagram} \label{diag2}
\[
\divide\dgARROWLENGTH by3
\begin{diagram}
\node{0} \arrow{se} \node{0} \arrow{s} \node[2]{0} \arrow{s} \node[1]{0} \arrow{s} \\
  \node{0} \arrow{e} \node{H_0 \oplus H_1} \arrow[2]{e,t}{i_0 \oplus i_1} 
      \arrow[2]{s,l}{\hl \circ (i_0 \oplus i_1)}  \arrow[2]{se,t,1}{\hl_0 \oplus -\hl_1}
        \node[2]{H} \arrow[1]{e,t}{\pi} \arrow[2]{s,r}{J}
           \arrow{sw,-} \node[1]{TG} \arrow[3]{s,r}{I} \arrow[1]{e} \node[1]{0} \arrow[2]{s}\\
\node[3]{} \arrow{sw,t}{\hl} \\    
\node{0} \arrow{0} \node{H^*} \arrow[2]{e,t}{i_0^* \oplus i_1^*} \arrow[1]{s,l}{\pi'} 
					\node[2]{H_0^* \oplus H_1^*} \arrow[1]{e,t,-}{\pi^{\wedge}} \arrow[1]{s,-} 
        \arrow[1]{se,t}{\pi_0 \oplus \pi_1} \node[1]{} \arrow{e} \node[1]{TG^{\wedge}} \arrow{e}
					 \arrow{s,r}{Id} \node{0} \\ 
\node{0} \arrow{e} \node{\hl(TG)} \arrow[1]{s} \arrow[3]{e,t}{K} \node[2]{} \arrow{s}   
       \node[1]{TG_0 \oplus TG_1} \arrow{e} \arrow{s} \arrow[1]{se} \node{TG^{\wedge}} \arrow{e}
\arrow{s} \node{0}\\
\node[2]{0} \arrow[2]{e} \node[2]{\hl(TG)^{\wedge}} \arrow{s} \arrow{e,t}{Id} \node{\hl(TG)^{\wedge}}
\arrow{s} \arrow{e} \node[1]{0}\\
\node[4]{0} \node{0} 
\end{diagram}
\]
\end{Diagram}

\noindent
Let ${\rm pr}_i: TG_0 \oplus TG_1 \ra TG_i$ is the projection onto the
indicated factor and for $v \in H$ suppose that ${\rm pr}_0(I(\pi(v))) = 0$.  By definition this means
that there is some $v_0 \in H_0$ such that
$\hl(v)|_{H_0} = \hl_0(v_0)$.  That is,
$$\lambda(v_0,w) = \lambda(v,w) ~\forall w \in H_0.$$
It follows that $\lambda(v-v_0,w) = 0$ for all $w \in H_0$ and hence that $(v-v_0) \in H_0^{\perp}
= H_1$.  Therefore $v = v_0 + (v-v_0) = (i_0 \oplus i_1)(v_0,v-v_0)$ and hence $\pi(v) = 0 \in
TG$.  Hence ${\rm pr_0} \circ I$ is injective and a completely similar argument shows the same for
${\rm pr_1} \circ I$.  Now by Lemma \ref{lemmacok} applied to 
$\hl_0 \oplus -\hl_1 = J \circ (i_0 \oplus i_1)$ we see that there is a short exact sequence
$$ 0 \ra TG \ra TG_0 \oplus TG_1 \ra \hl(TG)^{\wedge} \ra 0.$$
Now $|\hl(TG)^{\wedge}|=|TG|$ and so $|TG|^2 = |TG_0|\cdot|TG_1|$.  But since ${\rm
pr}_i \circ I$ is injective $|TG_i| \geq |TG|$ and so $|TG|=|TG_i|$.  We conclude that
that ${\rm pr}_i \circ I$ is an isomorphism.  

We show next that $\delta^c(\kappa_0) \oplus -\delta^c(\kappa_1)$ is identically zero on 
$I(TG) \subset TG_0 \oplus TG_1$.  Given $v \in H$ there is a non-zero integer $r$ and elements $v_i
\in H_i$ for $i=0,1$ such that $v = r(v_0+v_1)$.  Thus, for any $w \in H_0$,
$\lambda(v,w) = \lambda(rv_0,w) = r\lambda_0(v_0,w)$ and hence $\hl(v)|_{H_0} = r\hl_0(v_0)$.  Where as,
for any $w \in H_1$, $\lambda(v,w) = \lambda(rv_1,w) = -r\lambda_1(v_1,w)$ and $\hl(v)|_{H_1} =
r\hl_1(-v_1)$.  We now calculate
$$\begin{array}{cl}
(\delta^c(\kappa_0) \oplus -\delta^c(\kappa_1))[I(\pi(v))] 
&= [\hl_0^{-1}(\hl(v)|_{H_0},\hl(v)|_{H_0} + \alpha_0) \\
& \hskip 4cm - \hl_1^{-1}(\hl(v)|_{H_1},\hl(v)|_{H_1} + \alpha_1)]/2 \\
&= [\lambda(v,v_0) + \alpha_0(v_0) - \lambda(v,-v_1) + \alpha_1(-v_1)]/2r \\
&= [\lambda(v,v_0 + v_1) + \alpha(v_0 + v_1)]/2r \\
&= [\lambda(v,v) + \alpha(v)]/2. 
\end{array}$$
Now this last quantity is integral as $\kappa \in \FF^c$ and $\alpha$ is characteristic for
$\lambda$.  Thus $$(\delta^c(\kappa_0) \oplus -\delta^c(\kappa_1))[I(\pi(v))] = 0 \in \Q/\Z$$ for
every $\pi(v) \in TG$.  So it follows that
$$\theta_{\kappa} := (I \circ {\rm pr}_1) \circ (I\circ {\rm pr}_0)^{-1}: TG_0 \ra TG_1$$
is an isometry from $\delta^c(\kappa_0)$ to $\delta^c(\kappa_1)$.  In fact $I(TG) = \{
([x_0],\theta_{\kappa}([x_0])) | [x_0] \in TG_0\} \subset TG_0 \oplus TG_1$ and since
$\delta^c(\kappa_0)
\oplus -\delta^c(\kappa_1)$ vanishes identically on $I(TG)$ we see that for all $[x_0] \in TG_0$
$$\delta^c(\kappa_0)([x_0]) = \delta^c(\kappa_1)(\theta_{\kappa}([x_0])).$$

\noindent
{\em (4.)}  Observe that 
$J(H) = \{ (x_0,x_1) \in H_0^* \oplus H_1^* | \theta_{\kappa}([x_0]) = [x_1] \}$ which is the underlying
abelian group of
$\kappa_0 \cup_{\theta_{\kappa}} \kappa_1^-$.  As above, for any $v \in H$ we obtain a
non-zero integer $r$ and $v_i \in H_i$, $i=0,1$, such that $ v = r(v_0+v_1)$,
$\hl(v)|_{H_0} = r\hl_0(v_0)$ and $r\hl(v)|_{H_1} = r\hl_1(-v_1)$.  A similar calculation to preceding
one yields 
$$\begin{array}{cl}
(\kappa_0 \cup_{\theta_{\kappa}} \kappa_1^-)(J(v)) & = (\kappa_0 \cup_{\theta_{\kappa}} \kappa_1^-)
((\hl(v)|_{H_0},\hl(v)|_{H_1})) \\
&= \hl_0^{-1}(\hl(v)|_{H_0},\hl(v)|_{H_0} + \alpha_0) 
- \hl_1^{-1}(\hl(v)|_{H_1},\hl(v)|_{H_1} +\alpha_1) \\  
& = [\lambda(v_0,v) + \alpha_0(v_0)]/r - [\lambda(-v_1,v) + \alpha(-v_1)]/r  \\  
& = [\lambda(v_0,v) + \alpha(v_0) + \lambda(v_1,v) + \alpha(v_1)]/r \\  
& = \lambda(v_0+v_1,v)/r + \alpha(v_0+v_1)/r \\  
& = \lambda(v,v) + \alpha(v) \\
& = \kappa(v).
\end{array}$$
So $J$ is an isometry from $\kappa$ to $\kappa_0 \cup_{\theta_{\kappa}} \kappa_1^-$.
\end{proof}

\section{Stable equivalence of quadratic functions} \label{secseqf}
In this section we shall establish that the induced family of quadratic linking functions is a
complete stable invariant of quadratic functions.  Our method of proof is inspired by Wilkens' proof
of the analogous topological result (see Section \ref{chaptop}.\ref{secextdiff}).

\begin{Lemma} \label{lemmastabclass}
Let $\kappa_0(H_0,\lambda_0,\alpha_0)$ and $\kappa(H_1,\lambda_1,\alpha_1)$ be nondegenerate quadratic
functions in the family $\QF^*$.  Suppose that $\theta$ is an isometry of the induced quadratic 
linking families, $\theta:\delta^*(\kappa_0) \cong \delta^*(\kappa_1)$.  Then, there exists a nonsingular
quadratic function $\mu_0$ and an isometry 
$\Theta:\kappa_0 \oplus \mu_0 \cong \kappa_1 \oplus (\kap)$ such that $\delta^*(\Theta) = \theta$.
\end{Lemma}
  
\begin{proof}
Consider the four isometric embeddings:
$$\begin{array}{cccc}
j_0:& \kappa_0 & \hra & (\kappa_0  \cup \kappa_0^{-})\\
& v & \mapsto & (\hat{\lambda}_0(v),0)\\
\end{array} \hskip 2cm
\begin{array}{cccc}
j_1:& \kappa_0^{-} & \hra & (\kappa_0  \cup \kappa_0^{-})\\
& v & \mapsto & (0,-\hat{\lambda}_0(v))\\
\end{array}$$
$$\begin{array}{cccc}
i_0:& \kappa_0 & \hra & (\kap)\\
& v & \mapsto & (\hat{\lambda}_0(v),0)\\
\end{array}
\hskip 2cm
\begin{array}{cccc}
i_1:&\kappa_1^{-} & \hra & (\kap)\\
& v & \mapsto & (0,-\hat{\lambda}_1(v)).\\
\end{array}
$$
We shall prove that there is a further isometric embedding $I:(\kappa_0  \cup \kappa_0^{-}) \hra (\kap)
\oplus  (\kap)^-$, a nonsingular quadratic function $\mu_0$ and an isometry $\Phi:(\kappa_0  \cup
\kappa_0^{-}) \oplus \mu_0^- \cong (\kap) \cup (\kap)^-$ such that the following diagram
commutes.\\
\[
\begin{diagram}
\node{\kappa_0} \arrow{e,t}{j_0} \arrow{se,t}{i_0} 
 \node{(\kappa_0  \cup \kappa_0^{-})} \arrow{e,t}{{\rm Id} \oplus 0} \arrow{se,t}{I}
  \node{(\kappa_0  \cup \kappa_0^{-}) \oplus \mu_0^-} \arrow{s,r}{\Phi} \\ 
\node[2]{(\kappa_0  \cup_{\theta} \kappa_1)} \arrow{e,t}{{\rm Id} \oplus 0}
 \node{(\kap) \oplus (\kap)^-}\\
\end{diagram}
\]
\noindent 
Let $\theta':H_0^* \ra H_1^*$ be any homomorphism extending $\theta$.  Recall that
$$(\kappa_0  \cup \kappa_0^{-})(x_0,y_0) = \lambda_0^{-1}(x_0,x_0 + \alpha_0) -
\lambda_0^{-1}(y_0,y_0+\alpha_0)$$ 
on $\{ (x_0,y_0) \in H_0^* \oplus H_0^* : \pi_0({x_0}) = \pi_0({y_0}) \}$.  We define
$$ I((x_0,y_0)) = [(x_0,\theta'(y_0)),(y_0,\theta'(y_0))]$$
and verify that $I: \kappa_0 \cup \kappa_0^- \hra (\kap) \oplus (\kap)^-$ is an isometric embedding.
$$
\begin{array}{l}
(\kap) \oplus (\kap)^-[I(x_0,y_0)] \\  
\hskip 2cm = \lambda_0^{-1}(x_0,x_0+\alpha_0)-\lambda_0^{-1}(y_0,y_0+\alpha_0) \\
\hskip 5cm +\lambda_1^{-1}(\theta'(y_0),\theta'(y_0)+\alpha_1)
-\lambda_1^{-1}(\theta'(y_0),\theta'(y_0)+\alpha_1)  \\
\hskip 2cm = (\kappa_0 \cup \kappa_0^-)([(x_0,y_0)]
\end{array}
$$
\noindent
Since $I(j_0(v)) = [(\hat{\lambda}_0(v),0),(0,0)] = ({\rm Id} \oplus 0)(i_0(v))$ the left side of
the diagram commutes.  Now $I(\kappa_0 \cup \kappa_0^-)$ is a nonsingular subform of 
$(\kap) \oplus (\kap)^-$ which is also nonsingular and hence there is a nonsingular form $\mu_0$ and an
isometry 
$$\Phi:(\kappa_0 \cup \kappa_0^-) \oplus \mu_0^- \cong (\kap) \oplus (\kap)^-$$
which extends $I$.  Thus the right side of the diagram commutes.  

Now as $\Phi$ is an isometry $\Phi(j_0(\kappa_0)^{\perp}) = (\Phi(j_0(\kappa_0)))^{\perp}$. 
But $j_0(\kappa_0)^{\perp} = j_1(\kappa_0^{-}) \oplus \mu_0^{-}$ and
$(\Phi(j_0(\kappa_0)))^{\perp} = (i_0(\kappa_0))^{\perp} = i_1(\kappa_1^{-})
\oplus (\kap)^{-}$.  We therefore define $\Theta_1:=\Phi|_{j_0(\kappa_0)^{\perp}}$ and then
$$
\begin{array}{c}
\Theta_1: j_1(\kappa_0^-) \oplus \mu_0^- \cong i_1(\kappa_1^-) \oplus (\kap)^- \\
\end{array}
$$
is an isometry.  It follows that $\Theta:= (i_1^{-1} \oplus {\rm Id}_{\kap}) \circ \Theta_1 \circ
(j_1 \oplus {\rm Id}_{\mu_0})$ is an isometry, $\Theta: \kappa_0 \oplus \mu_0 \cong \kappa_1 \oplus
(\kap)$.  It remains
to show that $\Theta$ induces $\theta$.  Let the underlying modules and forms
of $\mu_0$ and $\kap$ be $(H_{\mu_0},\lambda_{\mu_0})$ and $(H_{\kappa},\lambda_{\kappa})$ respectively. 
We must determine the right most vertical homomorphism in the following commutative diagram.
\[
\begin{diagram}
\divide\dgARROWLENGTH	by2
\node{0} \arrow{e} 
 \node{H_0 \oplus H_{\mu_0}} \arrow{e,t}{\hat{\lambda}_0 \oplus
             \hat{\lambda}_{\mu_0}} \arrow{s,r}{\Theta} 
   \node{H_0^* \oplus H^*_{\mu_0}} \arrow{e,t}{\pi_0} \arrow{s,r}{(\Theta^{*})^{-1}}
    \node{TG_0} \arrow{e} \arrow{s,r,..}{?}\node{0}\\
\node{0} \arrow{e} 
 \node{H_1 \oplus H_{\kappa}} \arrow{e,t}{\hat{\lambda}_1 \oplus \hat{\lambda}_{\kappa}} 
  \node{H_1^* \oplus H_{\kappa}^*} \arrow{e,t}{\pi_1}
   \node{TG_1} \arrow{e} \node{0.} \\
\end{diagram}
\]
As both $\mu_0$ and $\kap$ are nonsingular $\pi_0|_{H_{\mu_0}^*}$ and $\pi_1|_{H_{\kappa}^*}$ are
identically zero and so it suffices to consider $\Theta|_{H_0 \oplus 0} \circ {\rm pr}_{H_1}$.
But $\Theta|_{H_0 \oplus 0} = (i_1^{-1} \oplus Id_{\kap}) \circ I \circ j_1$ and hence
$\Theta|_{H_0 \oplus 0} \circ {\rm pr}_{H_1} = \hl_1^{-1} \circ \theta' \circ \hl_0$.  But now one
verifies that $(\hl_1^{-1} \circ \theta' \circ \hl_0)^{*-1} = \theta'$.  Hence
$\Theta$ induces the same map on $TG_0$ as $\theta'$ which is
$\theta$ as $\theta'$ was chosen to extend $\theta$.  We summarize this argument with the
following subdiagram of the previous diagram
\[
\begin{diagram}
\node{0} \arrow{e} 
 \node{H_0} \arrow{e,t}{\hat{\lambda}_0} \arrow{s,l}{\Theta|_{H_0 \oplus 0}}
   \node{H_0^*} \arrow{e,t}{\pi_0} \arrow{s,r}{\theta'}
    \node{TG_0} \arrow{e} \arrow{s,r}{\theta} \node{0}\\
\node{0} \arrow{e} 
 \node{H_1} \arrow{e,t}{\hat{\lambda}_1} 
  \node{H_1^*} \arrow{e,t}{\pi_1} 
   \node{TG_1} \arrow{e} \node{0.} \\
\end{diagram}
\]
\end{proof}

\begin{Corollary} \label{corstabclass}
Let $\kappa_0(H_0,\lambda_0,\alpha_0)$ and $\kappa(H_1,\lambda_1,\alpha_1)$ be quadratic
functions in the family $\QF^*$.  Suppose that $\theta$ is an isometry of the induced quadratic 
linking families, $\theta:\delta^*(\kappa_0) \cong \delta^*(\kappa_1)$.  Then, there exist nonsingular
quadratic functions $\mu_0$ and $\mu_1$ and an isometry 
$\Theta:\kappa_0 \oplus \mu_0 \cong \kappa_1 \oplus \mu_1$ which induces $\theta$.
\end{Corollary}

\begin{proof}
Choose a section $\Phi_0 \in \Sec(\tau_{G_0})$ and let $\Psi_i \in \Sec(\rho_i)$, $i=0,1$, be sections
such that $\Phi(\Psi_0) = \Phi_0$ and $\Phi(\Psi_1) = \Phi_1 := \tau_{\theta}(\Phi_0)$.  Then since
$\theta$ is an isometry of quadratic linking families, $\theta|_{TG_0}$ is an isometry from
$q^*(\kappa_0)(\Phi_0) = q^*(\kappa_0(\Psi_0))$ to 
$q^*(\kappa_1)(\Phi_1) = q^*(\kappa_1(\Psi_1))$.  Applying Lemma \ref{lemmastabclass} to the
$\theta|_{TG_0}$ and setting
$\mu_1:=(\kappa_0(\Psi_0)
\cup_{\theta|_{TG_0}} \kappa_1(\Psi_1)^-)$ we conclude that there is
a nonsingular quadratic function $\mu_0$ and an isometry 
$$\Theta_N: \kappa_0(\Psi_0) \oplus \mu_0 \ra  \kappa_1(\Psi_1) \oplus \mu_1$$
such that $\delta^*(\Theta_N) = \theta|_{TG_0}$.  
Now $\kappa_0 = \alpha|_{F_0} \oplus \kappa(\Psi_0)$ and $\kappa_1 = \alpha_1|_{F_1} \oplus
\kappa(\Psi_1)$ and the isometry $\Theta_N$ is defined on $H(\Psi_0) \oplus H_{\mu_0}$ where $H_{\mu_0}$
is the domain of $\mu_0$. Thus it remains to use $\theta$ to define an isometry from $\alpha_0|_{F_0}$ to
$\alpha_1|_{F_1}$.  By the definition of $\Phi_i$, $\pi_i|_{F_i^*(\Psi_i)}: F_i^*(\Psi_i)
\ra \Phi_i(G_i^{**})$ is an isomorphism for $i=0,1$.  Also, by the definition of
$\tau_{\theta}$ and the fact that $\Phi_1 = \tau_{\theta}(\Phi_0)$, 
$\theta(\Phi_0(G^{**}_0)) = \Phi_1(G^{**}_1)$.  Moreover, there are isomorphisms, $R_{F_i}: F_i^*(\Psi_i) \ra F_i^*$, $i=0,1$  defined by restricting
the domain of homomorphisms in $F_i^*(\Psi_i)$
$$\begin{array}{cccc}
R_{F_i}: & F_i^*(\Psi_i) & \ra & F_i^* \\
& z & \mapsto & z|_{F_i}.
\end{array}$$
We therefore define $\Theta_D: F_0 \ra F_1$ so that $\Theta_D^{*-1}$ makes the following diagram
(in which all morphisms are isomorphisms) commute.
\[
\begin{diagram}
\node{F_0^*} \arrow{s,r,..}{\Theta_D^{*-1}} \node[2]{F_0^*(\Psi_0)} 
\arrow[2]{w,t}{R_{F_0}} \arrow[2]{e,t}{\pi_0|_{F_0^*(\Psi_0)}} \node[2]{\Phi_0(G^{**}_0)}
\arrow{s,r}{\theta|_{\Phi_0(G^{**}_0)}} \\
\node{F_1^*} \node[2]{F_1^*(\Psi_1)} \arrow[2]{w,t}{R_{F_1}} \arrow[2]{e,t}{\pi_1|_{F_1(\Psi_1)^*}}
\node[2]{\Phi_1(G^{**}_1)} \\
\end{diagram}
\]
Under the splittings $G_i = \Phi_i(G^{**}_i) \oplus TG_i$ and
$H_i^* = F_i^*(\Psi_i) \oplus H_i^*(\Psi_i)$ the elements $\beta_i$ and $\alpha_i$ are mapped to
$(\Phi_i(\tau_{G_i}(\beta_i)),\beta_i - \Phi_i(\tau_{G_i}(\beta_i)))$ and
$(p_{F_i}(\Psi)^*(\alpha|_{F_i}),\alpha(\Psi_i))$ respectively (see Sections
\ref{chapprelim}.\ref{secdoqf} and \ref{chapprelim}.\ref{secqlfam}).  
Now $\pi_i(\alpha_i) = \beta_i$, and so
$\pi_i|_{F_i^*(\Psi_i)}(p_{F_i}(\Psi)^*(\alpha|_{F_i})) = \Phi_i(\tau_{G_i}(\beta_i))$.
Since $\theta(\beta_0) = \beta_1$ and $\theta$ respects the splittings of
$G_0$ and $G_1$, $\theta(\Phi_0(\tau_{G_0}(\beta_0))) = \Phi_1(\tau_{G_1}(\beta_1))$.  Finally, as
$R_{F_i}(p_{F_i}(\Psi)^*(\alpha|_{F_i})) = \alpha|_{F_i}$ it follows by chasing the diagram defining
$\Theta_D$ that $\Theta_D^*(\alpha_1|_{F_1}) = \alpha_0|_{F_0}$.  That is, $\Theta_D$ defines an isometry
from $\alpha_0|_{F_0}$ to $\alpha_1|_{F_1}$.  By definition $\delta^*(\Theta_D) =
\theta|_{\Phi_0(G^{**}_0)}$.  We therefore define the isometry
$$\Theta := \Theta_D \oplus \Theta_N: \kappa_0 = [\alpha_0|_{F_0} \oplus \kappa_0(\Psi_0) \oplus \mu_0]
\ra [\alpha_1|_{F_1} \oplus \kappa_1(\Psi_1) \oplus \mu_1] = \kappa_1$$
and note that 
$$\delta^*(\Theta) = \delta^*(\Theta_D) \oplus \delta^*(\Theta_N) = \theta|_{\Phi_0(G^{**}_0)} \oplus
\theta|_{TG_0} = \theta.$$

\end{proof}

Having assembled the necessary results we conclude this chapter with the proof of Theorem
\ref{thmstabequivs} which we restate for the reader's convenience.

\renewcommand{\thetheorem}{\ref{thmstabequivs}}

\begin{theorem} \label{thmstabequivs2}
Let $\kappa_0$ and $\kappa_1$ be two quadratic functions in the same family $\FF^*$.  Then the
following are equivalent.
\begin{enumerate}
\item{$\kappa_0$ is stably $\FF^*$-equivalent to $\kappa_1$, $\kappa_0 \sim_{\FF^*} \kappa_1$.}
\item{The quadratic linking families induced by $\kappa_0$ and $\kappa_1$ are
isometric, $\delta^*(\kappa_0) \cong \delta^*(\kappa_1)$.}
\item{$\kappa_0$ and $\kappa_1^-$ are $\FF^*$-orthogonal, $\kappa_0 \perp_{\FF^*} \kappa_1^-$.}
\end{enumerate}
\end{theorem}

\begin{proof} 
\noindent
{\em (1.) $\Rightarrow$ (2.)}  The argument should be familiar by now.  Let $\mu_0$ and $\mu_1 \in \FF^*$
be nonsingular and $\Theta: \kappa_0 \oplus \mu_0 \ra \kappa_1 \oplus \mu_1$ be an isometry.  Then by
Lemma \ref{lemmaqlfam2}, $\delta^*(\Theta)$ induces an isometry from $\delta^*(\kappa_0 \oplus \mu_0)$ to
$\delta^*(\kappa_1 \oplus \mu_1)$.  But for $i=0,1$, $\delta^*(\kappa_i \oplus \mu_i) =
\delta^*(\kappa_i) \oplus \delta^*(\mu_i) = \delta^*(\kappa_i)$.  The last equality holding because
$\mu_i$ induces the trivial quadratic linking family because it is nonsingular.  Hence
$\delta^*(\Theta)$ defines an isometry from $\delta^*(\kappa_0)$ to $\delta^*(\kappa_1)$.

\noindent
{\em (2.) $\Rightarrow$ (3.)}
If $\theta:\delta^*(\kappa_0) \cong \delta^*(\kappa_1)$ is an isometry then by the proof of Corollary
\ref{corstabclass} there are sections $\Psi_0 \in \Sec(\rho_0)$ and $\Psi_1 \in \Sec(\rho_1)$  
such that $\theta_{TG_0}:\kappa_0(\Psi_0) \cong \kappa_1(\Psi_1)$ and such that there is an 
isometry $\Theta_D: \kappa_0|_{F_0} \cong \kappa_1|_{F_1}$.  We apply
Lemma \ref{lemmaglue} to $\theta_{TG_0}$ and observe that $\kappa_0(\Psi_0)$ and $\kappa_1(\Psi_1)^-$
are orthogonal inside $\kappa_0(\Psi_0) \cup_{\theta|_{TG_0}} \kappa_1(\Psi_1)^-$.  Hence $\kappa_0$ and
$\kappa_1$ are orthogonal.

\noindent
{\em (3.) $\Rightarrow$ (1.)}
Suppose that $\kappa_0 \perp_{\FF^*} \kappa_1^-$.  This means firstly there is an isometry $\Theta_D:
\kappa_0|_{F_0} \cong \kappa_1|_{F_1}$.  Secondly it means that there are
sections $\Psi_0 \in \Sec(\rho_0)$ and $\Psi_1 \in \Sec(\rho_1)$, a nonsingular quadratic function
$\kappa$ and embeddings $i_0:\kappa_0(\Psi_0) \hra \kappa$ and $i_1:\kappa_1(\Psi_1)^- \hra \kappa$ 
such that $i_1(\kappa_1(\Psi_1)^-) = i_0(\kappa_0(\Psi_0))^{\perp}$.  We calculate that
$$\begin{array}{cll}
\delta^*(\kappa_0(\Psi_0)) & \cong \delta^*(i_0(\kappa_0(\Psi_0))) \\
& \cong \delta^*(i_0(\kappa_0(\Psi_0))^{\perp -}) &\text{By Lemma \ref{lemmasplit} {\em 3}} \\
& \cong \delta^*(i_1(\kappa_1(\Psi_1)^-)^-) \\
& \cong -\delta^*(i_1(\kappa_1(\Psi_1)^-)) \\
& \cong -\delta^*(\kappa(\Psi_1)^-)\\
& \cong \delta^*(\kappa(\Psi_1)).
\end{array}$$
So by Lemma \ref{lemmastabclass} $\kappa_0(\Psi_0)$ and $\kappa_1(\Psi_1)$ are stably
$\FF^*$-equivalent.  Since $\kappa_0 = \kappa_0|_{F_0} \oplus \kappa_0(\Psi_0)$, $\kappa_1 =
\kappa_1|_{F_1} \oplus \kappa_1(\Psi_1)$ and $\kappa_0|_{F_0} \cong \kappa_1|_{F_0}$ we deduce
that $\kappa_0$ and $\kappa_1$ are stably $\FF^*$-equivalent.
\end{proof}

\chapter{Classification of the manifolds P} \label{chaptop}
It is an immediate consequence of Corollary \ref{corwilkwall} and Theorem \ref{thmstabequivs} that
two highly connected manifolds $P_0$ and $P_1$ are almost diffeomorphic if and only if the quadratic
linking families $Q(P_0)$ and $Q(P_1)$ are isometric.  Corollary \ref{corwilkwall} asserts that $P_0$
and $P_1$ are almost diffeomorphic if and only if the stable quadratic functions $[\kappa(P_0)]$ and
$[\kappa(P_1)]$ are equal and Theorem \ref{thmstabequivs} asserts (in part) that $[\kappa(P_0)] =
[\kappa(P_1)]$ if and only if $Q(P_0) \cong Q(P_1)$.  In Section \ref{secthmb} of this chapter we prove
the sharper result that every isometry $\theta: Q(P_0) \ra Q(P_1)$ is realized as $f_{!}$ for some almost
diffeomorphism $f: P_0 \ra P_1$.  Section \ref{seckap} is an aside which explains the topology behind
the key algebraic Lemmas \ref{lemmasplit} and \ref{lemmaglue}.  In Section \ref{secsmooth} we consider
the smooth classification of highly connected $(8j-1)$-manifolds.

\section{Extensions of diffeomorphisms} \label{secextdiff}
In Section \ref{secwilkc} we reported the following

\renewcommand{\thetheorem}{\ref{thmwilk3.2}}

\begin{theorem}[Wilkens' Thesis 3.2]
Let $L_0$ and $L_1$ be $2k$-dimensional handlebodies $(k>2)$ with boundaries $\del L_0 = P_0$ and
$\del L_1 = P_1$.  Suppose that $f:P_0 \ra P_1$ is a diffeomorphism.  Then there
are handlebodies $V_0, V_1$ with boundaries the standard sphere and a diffeomorphism $g$
which makes the following diagram commute.
\[
\divide\dgARROWLENGTH by 2
\begin{diagram}
\node{P_0} \arrow{e,t}{f} \arrow{s,j} 
 \node{P_1} \arrow{s,J}  \\
\node{L_0 \natural V_0} \arrow{e,t}{g}
 \node{L_1 \natural V_1}
\end{diagram}
\]
\end{theorem}

\noindent
We give Wilkens' proof of this result both for completeness and because the steps in the proof of
 Lemma \ref{lemmastabclass} are algebraic analogues inspired by the geometric steps in Wilkens' proof.
\begin{proof}
To begin form the manifold $M = L_0 \cup_f -L_1$ by identifying points in $\del L_0 = P_0$ and
$\del L_1 = P_1$ via $f$.  Using Van Kampen's theorem and the Mayer-Vietoris sequence for
$M=L_0 \cup -L_1$ we see that $M$ is a closed highly connected $2k$-manifold.  We wish to show
that $M \# -M$ is the connected sum of the double of $L_0 = L_0 \cup_{Id} -L_0$ and some other
highly connected manifold which will be called $M_0$.  To do this we embed a $2^k$-disc $D^{2k}$ 
in $L_1$ which itself is embedded in $M$ and let $N = (M-D^{2k}) \times [0,1]$.  Then there is
a diffeomorphism $F_1: \del N \cong M \# -M$.  Now if we consider the pair $(N,L_0)$ where $L_0$ is
embedded as $L_0 \times 0 \hra (M \times {0})$ then
$$\begin{array}{cl}
H_l(N,L_0) & \cong H_l(M-D^{2k}, L_0) \text {~as~} N \simeq M - D^{2k}\\
& \cong H_l(L_1 - D^{2k},P_1) \text{~by excision}\\
& \cong \left\{ \begin{array}{cr}
H_{k}(L_1,P_1) & {~l = k}\\
0 & \text{~else} \end{array} \right\}.
\end{array}$$
Thus $N$ has a handlebody decomposition based on $L_0 \times [0,1]$ with $k$-handles attached,
one for each element of a basis of $H_k(M-D^{2k}, L_0)$.  The attaching maps for the cores of
these handles are embeddings $S^{k-1} \hra L_0 \times 1$ and as $L_0 \times {1}$ is $(k-1)$-connected we
may assume after isotopy that the images of the attaching maps all lie in another $2k$-disc $D_1^{2k}
\hra L_0 \times 1$.  If follows that there is diffeomorphism 
$F_0:\del N \cong (L_0 \cup_{Id} -L_0) \# M_0$ where $M_0$ is the boundary of the manifold
obtained by attaching the $j$-handles to $D_1^{2k} \subset D^{2k+1} $.  Notice that since all the
necessary isotopies occur in $L_0 \times {1}$ that we have a fixed embedding of $L_0 = (L_0 \times
0) \hra \del N$.  Hence if we define $G$ to be the diffeomorphism $G:=F_0^{-1} \circ F_1$ we
have a commutative diagram of embeddings and diffeomorphisms 
\[
\begin{diagram}
\divide\dgARROWLENGTH by 2
\node[2]{((L_0 \cup_{Id} -L_0)-D^{2j})} \arrow{e} 
	\node{(L_0 \cup_{Id} -L_0) \# M_0} \arrow[2]{s,r}{G}\\
\node{L_0} \arrow{ne,t}{i_0} \arrow{se,t}{i_1}\\
\node[2]{((L_0 \cup_{f} -L_1)-D^{2j})} \arrow{e}
 \node{(L_0 \cup_{f} -L_1)\# -M}
\end{diagram}
\]
\noindent
such that $G|_{i_0(L_0)} = Id$.  If follows that $G$ defines a diffeomorphism 
$$G|_{-L_0 \# M_0} :-L_0 \# M_0 \cong -L_1 \# -M$$
of the complements of $i_0(L_0)$ and $i_1(L_0)$. 
Therefore $G|_{-L_0 \# M_0}$ induces a diffeomorphism on the boundary which must be $f$ since we
have the following commutative square.
\[
\begin{diagram}
\divide\dgARROWLENGTH by 2
\node{\del L_0} \arrow{s,l}{Id} \arrow{e,t}{Id}
 \node{-\del L_0} \arrow{s,r}{G|_{-\del L_0}} \\
\node{\del L_0} \arrow {e,t}{f} 
 \node{-\del L_1}
\end{diagram}
\]
\end{proof}

\begin{remark}
Wilkens' theorem is, in particular, a result about bordism of diffeomorphisms which is a topic that
has since received significant attention.  In particular see Kreck \cite{Kr1} and Quinn \cite{Q}. 
One might hope to prove this theorem from a more general point of view by arguing first that the bordism
of diffeomorphisms over $P$ are trivial and then performing surgery on the cobounding diffeomorphism to
make it a diffeomorphism of highly connected manifolds.
\end{remark}
  
\section{Almost smooth classification of the manifolds $P$} 
\label{secthmb}
\renewcommand{\thetheorem}{\ref{thmb},\ref{themwallqlfam} }

The following theorem combines the statements of Theorem \ref{thmb} and \ref{thmwallqlfam}.

\begin{Theorem} \label{thmB}
For every finitely generated abelian group $G$ and every characteristic (resp. even)
quadratic linking family $Q$ defined on $G$, there is a highly connected $(8j-1)$-manifold $P$, $j=1,2$
(resp $j>2$), such that $Q(P) \cong Q$.  The quadratic linking family $Q(P)$ of a highly connected
manifold $P$ is a complete  invariant of almost diffeomorphisms.  That is, if $P_0$ and $P_1$ are two highly connected 
$(8j-1)$-manifolds then there is an almost diffeomorphism $f: P_0 \ra P_1$ such that $f_! = \theta:
H^{4j}(P_0) \cong H^{4j}(P_1)$ if and only if $\theta$ defines an isometry of quadratic linking families
from $Q(P_0)$ to $Q(P_1)$.
\end{Theorem}

\begin{proof}
Let $Q=(G,q^*,\beta)$ be a quadratic linking family and let $\Phi$ be any section in $\Sec(\tau_G)$. 
Then $G = \Phi(G^{**}) \oplus TG$ where $TG$ is the torsion subgroup of $G$.  In the characteristic 
case Corollary \ref{chapprelim}\ref{corcharrealise} ensures that there is a nondegenerate quadratic
function $\kappa' = \kappa(H',\lambda',\alpha') \in \FF^{c}$ and an isomorphism $\theta': \Cok(\hl') \ra
TG$ such that $q^{c}(\kappa') = q^{c}(\Phi) \circ \theta$.  We define
$\kappa:=\kappa' \oplus \tau_G(\beta): H' \oplus G^* \ra \Z$ and leave the reader to verify that 
$\theta \oplus \Phi: \Cok(\hl) \oplus G^{**} \ra G$ defines an isometry from $\delta^{c}(\kappa)$
to $Q$.  The argument in the even case is similar but applies instead Theorem \ref{thmwallrealise} to
deduce the existence of an nondegenerate even quadratic function such $\kappa$ such that
$q^{ev}(\kappa) \cong q^{ev}(\Phi)$.

Now let $L_i$ be handlebodies such that $\del L_i = P_i \# \Sigma_{P_i}$.  Suppose that there is an
isometry
$$\theta: Q(P_0) = \delta^*(\kappa(L_0)) \ra \delta^*(\kappa(L_1)) = Q(P_1).$$  By Corollary
\ref{corstabclass} there is a nonsingular quadratic function $\mu_0$ and an isometry
$$\Theta : \kappa(L_0) \oplus \mu_0 \ra \kappa(L_1) \oplus \kap$$
which induces $\theta$, $\delta^*(\Theta) = \theta$.  Since $\mu_0$ is nonsingular the
boundary of $L(\mu_0)$ is a homotopy sphere which we shall denote by $\Sigma_0$.  Now, by Wall's
classification of handlebodies $\Theta$ is realized by a diffeomorphism 
$$g(\Theta): L_0 \natural L(\mu_0) \ra L_1 \natural L(\kap).$$
Restricting $g$ to the boundary of $L_0 \natural L(\mu_0)$ we obtain a diffeomorphism
$$f: P_0 \# \Sigma_{P_0} \# \Sigma_0 \ra P_1 \# \Sigma_{P_1}$$
which defines an almost diffeomorphism $f': P_0 \ra P_1$ such that $f'_{!} = f_! = \theta$.

Conversely, suppose that $f: P_0 \ra P_1$ be an almost diffeomorphism.  By Proposition \ref{propad}
there is a homotopy sphere $\Sigma$ and a diffeomorphism $f': P_0 \#\Sigma \ra P_1$ which agrees with
$f$ outside a disc.  Then
$\Sigma_{P_1} = \Sigma_{P_0 \# \Sigma} = \Sigma_{P_0} \# \Sigma_{\Sigma}$ and there is a diffeomorphism
$$f^{''}:P_0 \# \Sigma \# \Sigma_{P_1} = P_0 \# \Sigma_{P_0} \# \Sigma \# \Sigma_{\Sigma} \ra P_1 \#
\Sigma_{P_1}.$$
Let $L_{\Sigma}$ be a handlebody bounding $\Sigma \# \Sigma_{\Sigma}$.  Applying Theorem
\ref{thmwilk3.2} to $f^{''}: \del (L_0 \natural L_{\Sigma}) \ra \del L_1$ we deduce
that there are handlebodies
$V_0$ and $V_1$ cobounding standard spheres and a diffeomorphism 
$$h: L_0 \natural L_{\Sigma} \natural V_0 \ra L_1  \natural V_1$$
which extends $f^{''}$.  The map induced by $h$, 
$\Theta(h) : \kappa(L_0) \oplus \kappa(L_{\Sigma} \natural V_0)  \ra \kappa(L_1) \oplus \kappa(V_1)$, is
an isometry.  By Lemma \ref{lemmaqlfam2} it follows that 
$$f^{''}_{!}: \delta^*(\kappa(L_0)) \oplus \delta^*(\kappa(L_{\Sigma} \natural V_0))  
\ra \delta^*(\kappa(L_1)) \oplus \delta^*(\kappa(V_1))
$$
is an isometry.  However as $V_0$ and $L_{\Sigma} \natural V_0$ are both bounded by homotopy spheres,
$\delta^*(\kappa(V_0))$ and $\delta^*(\kappa(L_{\Sigma} \natural V_1))$ are both the trivial quadratic
linking family and therefore $f^{''}_{!}$ defines an isometry from $\delta^*(\kappa(L_0)) = Q(P_0)$ to
$\delta^*(\kappa(L_1)) = Q(P_1)$.  Finally, $f^{''}_{!} = f_{!}$ since we may assume that the two maps
differ only on the interior of a disc.
\end{proof}

\section{The manifolds  $L(\kap)$  and $L_0 \cup_{f^o} -L_1$}\label{seckap}
Let $P_0$ and $P_1$ be the boundaries of handlebodies and let $f:P_0 \ra P_1$ be an almost
diffeomorphism.  Recall that by Proposition \ref{propad} that there is a homotopy sphere $\Sigma$, a disc
$D^{8j-1} \hra P_0$ and a diffeomorphism $f':P_0 \# \Sigma \ra P_1$ such that $f'|_{P_0 - D^{8j-1}} =
f|_{P_0-D^{8j-1}}$ (where the connected sum of $\Sigma$ and $P_0$ is understood to occur inside
$D^{8j-1}$).  We define $f^o$ to be the restriction $f^o: = f|_{P_0 - D^{8j-1}}$ which is a
diffeomorphism onto its image.

\begin{Lemma} \label{lemmaL0ufL1}
Let $L_0=L(\kappa_0)$ and $L_1=L(\kappa_1)$ be handlebodies with associated nondegenerate
quadratic functions $\kappa_0 = \kappa(H_0,\lambda_0)$ and $\kappa_1 = \kappa(H_1,\lambda_1)$.  Let
$P_0$ and $P_1$ be the boundaries of $L_0$ and $L_1$, let $f: \del P_0 \ra \del P_1$ be an almost
diffeomorphism which induces  $f_{!} = \theta: H^{4j}(P_0) \ra H^{4j}(P_1)$ and let $\Sigma$ and $f^o$
be as above.  Then  $\theta: Q(P_0) \cong Q(P_1)$ is an isometry of quadratic linking families and
$L_0 \cup_{f^o} -L_1$ is  diffeomorphic to $\kappa_0 \cup_{\theta} \kappa_1^-$.
\end{Lemma}

\begin{proof}
The map $f^o$ is a diffeomorphism onto its image and, after perhaps rounding corners, the adjunction
space $L:=L_0 \cup_{f^o} -L_1$ is a smooth manifold with boundary.  The boundary of $L$ is obtained by 
gluing the disc $D^{8j-1} \subset P_0$ to a disc $D_1^{8j-1}$ which is the complement of the range of
$f^0$.  Thus $\del L$ is a homotopy sphere.  Now let $P^o$ be the subset $P_0-D^{8j-1} =
f^o(P_0-D^{8j-1}) \subset P_1$ contained in $L:=L_0 \cup_{f^o} -L_1$.  Applying Van-Kampen's theorem and
the Mayer-Vietoris sequence to the decomposition of $L$ into the triple $(L_0,L_1,P^o)$ we conclude that
$L$ is $(4j-1)$ connected.  It follows that $L$ falls under Wall's classification of handlebodies and we
let $\kappa=\kappa(H,\lambda,\alpha)$ be the quadratic function of $L$.  We will show that $\kappa$ is
isometric to $\kap$.  There are embeddings of
$i_0:H_0 \hra H$ and $i_1: H_1 \hra H$  embedding $\kappa_0$ and $\kappa_1^-$ respectively into $\kappa$ 
corresponding to the obvious embeddings of $L_0$ and $-L_1$ into $L_0 \cup_{f^o} -L_1$.  We claim that
that $\kappa_0^{\perp} = \kappa_1^-$ because firstly it is clear that spheres in $L_0$ representing
homology classes do not intersect spheres in $L_1$ and vice versa.  Secondly, suppose that
$\bar{v}:S^{4j} \hra L$ represents $v \in H_{4j}(L)$ and that $\lambda(v,v_0) = 0$ for every $v_0 \in
i_0(H_0)$.  It follows via the Whitney trick that $\bar{v}(S^{4j})$ may be isotoped away from all the
handles of $L_0$ and thence into $-L_1$.  Therefore $v \in i_1(H_1)$.  Applying Lemma \ref{lemmasplit} we
conclude that $\kappa \cong \kappa_0 \cup_{\theta_{\kappa}} \kappa_1^{-}$ and that 
$$\theta_{\kappa} : Q(P_0) = \delta^*(\kappa_0) \ra \delta^*(\kappa_1) = Q(P_1)$$
is an isometry of quadratic linking families.  It remains therefore to demonstrate
that $\theta_{\kappa} = \theta$.  

From the proof of Lemma \ref{lemmasplit} recall the projections $\pi_i: H_i^* \ra \Cok(\hl_0) TG_0 =
H^{4j}(P_0)$, $i=0,1$ and recall that $\theta_{\kappa}: TG_0 \ra TG_1$ is the unique isomorphism such
that $(\pi_0 \oplus \pi_1)(i_0^* \oplus i_1^*)(H^*) \subset TG_0 \oplus TG_1$ is the graph of
$\theta_{\kappa}$.  Now let $j_0: P^o \hra L_0$ and $j_1: P^o \hra -L_1$ denote the inclusions.  The
usual Mayer-Vietoris arguments show that there is a long exact sequence 
$$ 0  \lra   H^{4j}(L)  \stackrel{i_0^* \oplus i_1^*}{\lra}  H^{4j}(L_0) \oplus H^{4j}(-L_1) 
\stackrel{j_0^* -j_1^*}{\lra}  H_{4j-1}(P^o)  \lra  0. $$
Identifying $H^{4j}(P^o) = H^{4j}(P_0-D^{8j-1})$ we see that 
$(j_0^* - j_1^*)(x_0,x_1) = [x_0] - f^*[x_1]$.  Hence 
$$\begin{array}{cl}
(\pi_0 \oplus \pi_1)(i_0^* \oplus i_1^*)(H^*) & = \{ ([x_0],[x_1]) | f^{*-1}([x_0]) = [x_1] \} \\
& = \{ ([x_0],[x_1]) | \theta([x_0]) = [x_1] \}
\end{array}$$
and $\theta = \theta_{\kappa}$ as required.
\end{proof}

Given the manifold $L(\kap)$ with boundary $\Sigma$ we let $D^{8j}$ denote the almost smooth disc
obtained by taking the cone of $\Sigma$.  Hence $D^{8j}$ is a $PL$ manifold equipped with a smooth
structure except at the cone point which we denote by $p_0$.  Now form 
$M(\kap) := L(\kap) \cup_{\Sigma} D^{8j}$ which is an almost smooth manifold with singular point
$p_0$.  

\begin{Lemma} \label{lemmahcbd}
Let $L_0 = L(\kappa_0)$ and $L_1=L(\kappa_1)$ be handlebodies associated to nondegenerate quadratic
functions $\kappa_0 = \kappa(H_0,\lambda_0,\alpha_0)$ and $\kappa_1 = \kappa(H_1,\lambda_1,\alpha_1)$. 
If $\theta: \delta^*(\kappa_0)\cong \delta^*(\kappa_1)$ is an isomorphism of the induced quadratic
linking families then
\begin{enumerate}
\item{there is an embedding of the disjoint union $L(\kappa_0) \sqcup -L(\kappa_1) \hra M(\kap)$
 with complement an almost smooth h-cobordism from $\del L_0$ to $- \del L_1$,} 
\item{there is an almost diffeomorphism $f(\theta):\del L_0 \cong \del L_1$ realizing $\theta$ on
homology,}
\item{the manifold $M(\kap)$ is almost diffeomorphic to the manifold $L_0 \cup_{f(\theta)} -L_1$ formed
by gluing $L_0$ to $-L_1$ with $f(\theta)$. }
\end{enumerate}
\end{Lemma}

\begin{proof}
{\em 1.} Recall that the group $H = \{(x_0,x_1) \in H_0^* \oplus H_1^* |\theta([x_0]) = [x_1] \}$
is the domain $\kap$ and that there are isometric embeddings
$$\begin{array}{cccl}
i_0: & H_0 & \hra & H\\
& v_0 & \mapsto & (\hat{\lambda}_0(v_0),0)
\end{array} \hskip 2cm
\begin{array}{cccl}
i_1: & H_1 & \hra & H\\
& v_1 & \mapsto & (0,-\hat{\lambda}_1(v_1)).
\end{array}
$$
from Lemma \ref{lemmaglue}.  For $j=0,1$, the group $i_j(H_j)$ is a summand of $H$ for it is the
kernel of the homomorphism $H \ra H_j^*$ which sends a pair $(x_0,x_1) \in H$ to $x_j \in H_j^*$. 
Thus if $\{ v_{01}, \dots, v_{0m} \}$ and $\{ v_{11}, \dots, v_{1n} \}$ are bases for $H_0$ and
$H_1$ respectively, then there are bases $\{y_1, \dots, y_N \}$ and $\{z_1, \dots, z_N \}$ for $H$
such that $i_0(v_{0j}) = y_j$ for $1 \leq j \leq l$ and $i_1(v_{1j}) = z_j$ for $1 \leq j \leq
m$.  According to Wall's Theorem \ref{thm1wall1} the handlebody $L:=L(\kap)$ has a presentation for
each of the bases $\{ y_i \}$ and $\{ z_i \}$.  In fact, Wall proved the following.  Let $D^{8j}
\hra L$ be a fixed disc in the interior of $L$ and let $K$ be the closure of $L - D^{8j}$.  Then
$$H_{4j}(L) \cong H_{4j}(L,D^{8j}) \cong H_{4j}(K,\del D^{8j}) \cong \pi_{4j}(K,\del D^{8j}).$$
If $\{z'_1,\dots,z'_n\}$ is any base for $H_{4j}(L)$ then for each $z'_i$ there are
pairwise disjoint embeddings of $\bar{h}_i: (D^{4j},\del D^{4j}) \hra (K,\del D^{8j})$ which are
transverse to $\del D^{8j}$ along
$\bar{h}_i(\del D^{4j})$ and which represent the element in $\pi_{4j}(K,\del D^{8j})$ corresponding
to $z'_i$.  These embeddings can be thickened to handles $h_i:(D^{4j} \times D^{4j}, \del D^{4j}
\times D^{4j}) \hra (K,\del D^{8j})$ and the resulting handle body $L' \subset L$ obtained from the
union of $D^{8j}$ and the union of the images of the $h_i$ is diffeomorphic to $L$.  In fact
$L - L'$ is an $h$-cobordism from $\del L$ to $\del L'$.  We now choose disjoint discs $D^{8j}_0
\hra L$ and $D^{8j}_1 \hra L$ one for each basis $\{ y_i \}$ and $\{ z_i \}$ and attach
handles $h_{0i}$ and $h_{1j}$ to these discs as above but only for basis elements
$\hl_0(v_{0j})$ and $\hl_1(v_{1j})$.  The manifolds $D^{8j}_0
\cup(\cup_{i=1}^l h_{0i})$ and $D^{8j} \cup (\cup_{j=1}^m h_{1j})$ thus define embeddings of
$L_0$ and $-L_1$ into $L$.  Recall that $\lambda$
denotes the intersection form of $L$ and that
$$\lambda[(\hl_0(v_0),0),(0,-\hl_1(v_1)] = 0$$
for every $v_0 \in H_0$ and $v_1 \in H_1$.  It follows that the handles $h_{0i}$ and $h_{1j}$ have
algebraic intersection zero and so by the Whitney trick they may be isotoped to have empty
geometric intersection.  We may therefore assume that $L(\kappa_0)$ and $-L(\kappa_1)$ are
disjointly embedded into $L(\kap) \subset M(\kap)$ and we use the same symbols to denote the 
disjointly embedded manifolds.  

Let $N := M(\kap) - {\rm Int}(L_0 \sqcup -L_1)$ be the complement of the interior
of $L_0 \sqcup -L_1$ inside $M(\kap)$ which is a codimension zero submanifold of $M$.  
We orient $N$ compatibly with $M$ and then $\del N = \del L_0 \sqcup -\del L_1$.  Let $J:L_0 \sqcup -L_1
\hra M$, $i_N: N \hra M~$ and $i_{\del N}: \del N \hra L_0 \sqcup -L_1$ denote the obvious inclusions and
$E^*: H^*(M,L_0 \sqcup -L_1) = H^*(N,\del N)$ the excision isomorphisms.  The long exact cohomology
sequences of the pairs $(N,\del N)$ and $(M,L_0 \sqcup -L_1)$ are linked in the following commutative
diagram (where we abbreviate $L_0 \sqcup -L_1$ to $L_{0\sqcup 1}$).
\[
\begin{diagram} 
\divide\dgARROWLENGTH by2
\node{H^0(N,\del N)} \arrow{e} \node{H^0(N)} \arrow{e} \node{H^0(\del N)} \arrow{e} \node{H^1(N,\del N)}
\arrow{e} \node{H^1(N)} \\
\node{H^0(M,L_{0\sqcup 1})} \arrow{e} \arrow{n,l}{E^0}
\node{H^0(M)} \arrow{e,t}{J^0} \arrow{n,l}{i_N^0} 
\node{H^0(L_{0\sqcup 1})} \arrow{e} \arrow{n,l}{i_{\del N}^0}
\node{H^1(M,L_{0\sqcup 1})} \arrow{e} \arrow{n,l}{E^1} 
\node{H^1(M)} \arrow{n,l}{i^0_N}\\
\node{0} \arrow{e} \arrow{n,l}{=}
\node{\Z} \arrow{e} \arrow{n,l}{=}
\node{\Z \oplus \Z} \arrow{e} \arrow{n,l}{=}
\node{\Z} \arrow{e} \arrow{n,l}{=} \node{0} \arrow{n,l}{=} \\
\end{diagram}
\]
As $E^0$, $I^0_{\del N}$ and $E^1$ are all isomorphisms we conclude that $i_N^0$ is an isomorphism and
hence that $N$ is connected.  A simple application of Van Kampen's theorem shows that $N$ is also simply
connected.  Applying Poincar\'{e} duality to the manifold $N$ we have that $H_{8j-*}(N) \cong
H^*(N,\del N) \cong H^*(M,L_0 \sqcup -L_1)$.  We may calculate the groups $H^*(M,L_0 \sqcup -L_1)$
again using the long exact cohomology sequence of the pair $(M,L_0 \sqcup -L_1)$.  This sequence is zero
below dimension $(8j-1)$ except in dimensions $4j$ and $4j+1$.  The following diagram shows the
relative cohomology sequences for the pairs $(M,L_0 \sqcup -L_1)$ and $(N,\del N)$ on the left.  On
the right it shows an isomorphic sequence extending the right hand side of Diagram \ref{diag1} in Lemma
\ref{lemmaglue} from where we recall the projections $\pi_0: H_0^* \ra TG_0$, $\pi_1:H_1^* \ra TG_1$,
$\pi: H \ra TG = \{([x_0],\theta[x_0] \} \subset TG_0 \oplus TG_1$ and $\Pi: H_0^* \oplus H_1^* \ra
TG^{\wedge}$.

\begin{Diagram} \label{diag3}
\[
\divide\dgARROWLENGTH by2
\begin{diagram}
\node{0} \arrow{s} \node[1]{0} \arrow{s} \node{0} \arrow{s} \node{0} \arrow{s} \\
\node[1]{H^{4j}(M)} \arrow[1]{e,t}{i_N^*} \arrow[1]{s,r}{J^*}
          \node[1]{H^{4j}(N)} \arrow[1]{s,r}{I^*} 
\node{H} \arrow{e,t}{\pi} \arrow{s,r}{i} \node{TG} \arrow{s} \\   
\node[1]{H^{4j}(L_0 \sqcup -L_1)} \arrow[1]{s,l}{\delta_1} 
 \arrow[1]{e,t}{J_0^* \oplus J_1^*} \node[1]{H^{4j}(\del N)} \arrow{s,l}{\delta_2} 
\node{H_0^* \oplus H_1^*} \arrow{e,t}{\pi_0 \oplus \pi_1} 
\arrow{s,r}{\Pi} \node{TG_0 \oplus TG_1} \arrow{s}\\
\node[1]{H^{4j+1}(M, L_0 \sqcup -L_1)} \arrow{s} \arrow{e,t}{E^{4j+1}} 
\node{H^{4j+1}(N, \del N)} \arrow{s} 
\node{TG^{\wedge}} \arrow{e} \arrow{s} \node{TG^{\wedge}} \arrow{s}\\
\node[1]{0} \node{0}
\node{0} \node{0}\\ 
\end{diagram}
\]
\end{Diagram}

Applying Poincar\'{e} duality and the universal coefficient theorem for torsion groups we conclude
that
$$
H^*(N) \cong \left\{ \begin{array}{ll}
\Z & *=0, 8j-1\\
TG & *=4j\\
0 & \text{else.}
\end{array} \right\}
$$
With the identifications $H^{4j}(\del L_0) = TG_0$, $H^{4j}(-\del L_1) = TG_1$, $H^{4j}(N) = TG$ and
$H^{4j}(\del N) TG_0 \oplus TG_1$ and the notations $I_0 : \del L_0 \hra N$ and $I_1: -\del L_1 \hra N$
for the standard inclusions we conclude also that for 
$x = ([x_0],\theta[x_0]) \in H^{4j}(N)$, $I_0^*(x) = [x_0]$ and $I_1^*(x) = \theta[x_0]$.  Given our
calculation of $H^*(N)$ it follows that both
$I_0$ and
$I_1$ are $I_0$ and $I_1$ are (co)homology equivalences
of simply connected manifolds.  So by Whitehead's theorem both $I_0$ and $I_1$ are homotopy
equivalences.   Thus $N$ is an $h$-cobordism from $\del L_0$ to  $- \del L_1$.  It is an almost smooth
manifold since the singular point of $M$ in contained in $N$.

\noindent 
{\em (2.)}  The almost smooth manifold $N$ has one singular point with associated exotic sphere
$\Sigma := \del L(\kap)$.  As described in Section \ref{secmaness} we may construct a smooth
$h$-cobordism from $\del L_0 \# \Sigma$ to $\del L_1$ by thickening a path from a point on $\Sigma$ to
$\del L_0$ and then cutting out this thickened path to produce $N' \subset N$.  We note that cutting
out this path has no effect on the $4j$th cohomology groups of the manifolds involved.  We denote the
inclusions of the boundary components of the modified manifold by $I_0': \del L_0 \# \Sigma \hra N'$ and
$I_1':-\del L_1
\hra N'$.  Now let  $F: N' \ra (\del L_0 \# \Sigma) \times [0,1]$ be a diffeomorphism which restricts to
the identity on $\del L_0 \# \Sigma$.  The diffeomorphism at the other end of $N'$ is a diffeomorphism 
$F_1:=F \circ I_1': -\del L_1 \ra -(L_0 \# \Sigma)$.  Given $x \in H^{4j}(N')$ we have from above that 
$I_1^{'*}(x) = \theta I_0^{'*} (x)$.   But since
$I_0' \circ F|_{\del L_0 \# \Sigma} = Id_{\del L_0 \# \Sigma}$ we calculate that
$$\begin{array}{cl}
F_1^*(I_0^{'*}(x)) & = (F \circ I_1')^*((I_0')^*(x))\\
& = I_1^{'*}(F^*(I_0^{'*}(x)))\\
& = I_1^{'*}((I_0' \circ F_{\del L_0 \# \Sigma})^*(x))\\
& = I_1^{'*}(x) \\
& = \theta(I_0^{*'}(x)).
\end{array}$$
Thus $F_1$ induces $\theta$ on cohomology.  Therefore we take $f(\theta)$ to be the map 
$F_1^{-1}: \del L_0 \# \Sigma \ra \del L_1$ which we may consider as an almost diffeomorphism from 
$\del L_0$ to $\del L_1$ with one singular point whose associated exotic sphere is
$\Sigma$.  The isomorphism $f(\theta)_{!} = (F_1^{-1})^{*-1} = F_1^* = \theta$.

\noindent
{\em (3.)}   
The almost diffeomorphism $f(\theta)$ has one singular point so let $f^o=f(\theta)|_{\del L_0 -
D^{8j-1}}$.  By Lemma \ref{lemmaL0ufL1} the adjunction space $L_0 \cup_{f^o} -L_1$ space is diffeomorphic
to $L(\kap)$.  If we now form the adjunction space $L_0 \cup_{f(\theta)} -L_1$ then the effect is to
identify the two discs from which the homotopy sphere $\Sigma = \del (L_0 \cup_{f^o} -L_1)$ is formed. 
But this has the same effect as coning off $\Sigma$  to produce $M(\kap)$.  Hence $M(\kap)$ and $L_0
\cup_{f(\theta)} -L_1$ are diffeomorphic except at one singular point.
\end{proof}

\begin{remark} \label{remaltproof}
Note that for highly connected rational homology spheres bounding handlebodies
Lemmas \ref{lemmaL0ufL1} and \ref{lemmahcbd} together give a proof that the quadratic
linking function is a well defined and complete invariant of almost diffeomorphisms.  This second
proof is independent of Theorem \ref{thmwilk3.2}.  Its engine is the $h$-cobordism theorem which is the
heart of Lemma \ref{lemmaL0ufL1}.  We note also that with a little care this alternate proof could be
extended to all highly connected $(8j-1)$-manifolds and thus give a second proof of Theorem \ref{thmB}
which is independent of Theorem \ref{thmwilk3.2}.
\end{remark}


\section{Smooth classification of the manifolds $P$} \label{secsmooth}
Recall Corollary \ref{corwilkwall} which states that a highly connected $(8j-1)$-manifold $P$ is
classified up to almost diffeomorphism by its stable quadratic function $[\kappa(P)]$.  We begin this
section by defining a smooth stable quadratic function for $P$ and showing that it is a complete
diffeomorphism invariant of those $P$ which bound handlebodies.  
We shall call a quadratic function $\mu \in \FF^{8j}$  {\em spherical} if and only if the boundary of
$L^{8j}(\mu)$ is diffeomorphic to the standard sphere.  Now we have  already noted in the proof of Lemma
\ref{lemmahcb2} that every homotopy sphere
$\Sigma$ bounding a handlebody lies in $bP_{8j}$ and by Theorem \ref{thmstolz}, $\Sigma$ is standard if
and only if $s_1(\Sigma) = 0 \in \Q/\Z$.  We can therefore characterize spherical quadratic functions
completely algebraically as follows.

\begin{Lemma} \label{lemmaspherical}
A quadratic function $\mu(H,\lambda,\alpha) \in \FF^{8j}$ is spherical if and only if it is nonsingular
and  
$$s_1(\del L(\mu)) = (1/|bP_{8j}|) \left(\tilde{S}_{(j,2j)}\lambda^{-1}(\alpha,\alpha) -
\frac{1}{8}\sigma(\mu) \right) = 0 \in \Q/\Z.$$
\end{Lemma}

\noindent
We shall call two quadratic functions $\kappa_0$ and $\kappa_1$ in $\FF^{8j}$ {\em smoothly}
equivalent if there are spherical quadratic functions $\mu_0$ and $\mu_1$ in $\FF^{8j}$ and an isometry
$\Theta: \kappa_0 \oplus \mu_0 \cong \kappa_1 \oplus \mu_1$.  It is not difficult to see that smooth
equivalence is in fact an equivalence relation.  So given a highly connected $(8j-1)$-manifold we
define the smooth quadratic function of $P$, $[\kappa(P)]^s_{\FF^{8j}}$, to be the smooth equivalence
class of $\kappa(L)$ where $L$ is any handlebody such that $\del L = P \# \Sigma_P$. 

\begin{Theorem} \label{thmsphericalclass}
Two highly connected $(8j-1)$-manifolds $P_0$ and $P_1$ are diffeomorphic if and only if 
$$\Sigma_{P_0} \cong \Sigma_{P_1} \hskip 0.5cm \text{and} \hskip 0.5cm [\kappa(P_0)]^s_{\FF^{8j}} =
[\kappa(P_1)]^s_{\FF^{8j}}.$$
\end{Theorem}

\begin{proof}
That the manifolds $P_0 \# \Sigma_{P_0}$ and $P_1 \# \Sigma_{P_1}$ are diffeomorphic if and only if
$[\kappa(P_0)]^s_{\FF^{8j}} = [\kappa(P_1)]^s_{\FF^{8j}}$ follows from Wilkens' Theorem \ref{thmwilkc}
and Wall's classification of handlebodies in a manner completely analogous to Corollary
\ref{corwilkwall}.  The theorem now follows immediately.
\end{proof}

\noindent
The determination of smooth equivalence for nondegenerate quadratic functions is a straightforward
matter and this allows us to give a more explicit diffeomorphism classification for rational homology
spheres.

\begin{Theorem} \label{thmsmoothclass}
Highly connected $(8j-1)$-dimensional rational homology spheres $P_0$ and $P_1$ are diffeomorphic 
if and only if 
$$Q(P_0) \cong Q(P_1), \hskip 0.5cm s_1(P_0) = s_1(P_1) \hskip 0.5cm \text{and} \hskip 0.5cm 
\Sigma_{P_0} \cong \Sigma_{P_1}.$$
Moreover, if $P_0$ and $P_1$ are diffeomorphic then every isometry $\theta:Q(P_0) \ra Q(P_1)$ is
realized by a diffeomorphism.

\end{Theorem}

\begin{proof}
By Theorem \ref{thmb}, $P_0$ and $P_1$ are almost diffeomorphic if and only if $Q(P_0) \cong Q(P_1)$. 
So we may suppose that there is a homotopy sphere $\Sigma$ such that $P_0$ is diffeomorphic to 
$P_1 \# \Sigma$.  But then 
$$s_1(\Sigma) = s_1(P_0) - s_1(P_1)  \text{~~~and~~~} \Sigma_{\Sigma} \cong \Sigma_{P_0} - \Sigma_{P_1}$$
and the first assertion follows since $\Sigma$ is diffeomorphic to the standard sphere if and only if
$s_1(\Sigma) = 0$ and $\Sigma_{\Sigma} \cong S^{8j-1}$.

Now suppose that $P_0$ and $P_1$ are diffeomorphic and that $\theta:Q(P_0) \ra Q(P_1)$ is an isometry. 
By Lemma \ref{lemmastabclass} there is a nonsingular quadratic function $\mu_0 \in \FF^{8j}$ and an
isometry
$$\Theta: \kappa_0 \oplus \mu_0 \ra \kappa_1 \oplus \kap$$
where $\kappa_i = \kappa(L_i)$ and $L_i$ is a handlebodies with $\del L_i = P_i \# \Sigma_{P_i}$,
$i=0,1$.  Of course by Wall's classification of handlebodies, $\Theta$ is realized by a diffeomorphism
$g: L_0 \natural L(\mu_0) \cong L_1 \natural L(\kap)$ and so it suffices to show that $\mu_0$ and $\kap$
are spherical.  By Lemma \ref{lemmaspherical} this means showing that $s_1(\del L(\mu_0)) = 0 =
s_1(\del L(\kap))$.  But now if $\kap = \kappa(H,\lambda,\alpha)$ and $\kappa_i =
\kappa(H_i,\lambda_i,\alpha_i)$ then 
$$\sigma(\kap) = \sigma(\kappa_0) - \sigma(\kappa_1) \text{~~~and~~~}
\lambda^{-1}(\alpha,\alpha) = \lambda_0^{-1}(\alpha_0,\alpha_0) -  \lambda_1^{-1}(\alpha_1,\alpha_1).$$
It follows immediately that
$$s_1(\del L(\kap)) = s_1(\del L_0) - s_1(\del L_1) = s_1(P_0) - s_1(P_1) = 0$$
and then we have
$$s_1(\del L(\mu_0)) = s_1(\del L_1) + s_1(\del(L(\kap)) -s_1(\del L_0) = 0.$$
\end{proof}

The question of the smooth equivalence of degenerate quadratic functions is a difficult algebraic
problem.  It is equivalent to calculating the the inertia group of the boundaries of the corresponding
handlebodies.  The inertia group of a manifold $M$, $I(M)$ is the subgroup of the group of
homotopy spheres for which $\Sigma \# M$ is diffeomorphic to $M$.  If $P$ is a highly connected
$(8j-1)$-manifold then $I(P) \subset bP_{8j}$ for if $[\Sigma] \neq 0 \in \Omega_{8j-1}^{HC}$ then $[P]
\neq [P \# \Sigma]$ and so $P$ and $P \# \Sigma$ cannot be diffeomorphic.  If in addition $P$ is a
rational homology sphere then $I(P) = 0$ for in this case $s_1$ is defined for $P$ and detects the
addition of any exotic sphere in $bP_{8j}$.  In cite \cite{Wi1} and \cite{Wi3} Wilkens investigated the
inertia groups of highly connected manifolds of dimension $7$ and $15$ and discovered them to be quite
subtle.  In particular, he found examples where $I(P_0) = 0 =I(P_1)$ but $I(P_0 \# P_1) \neq 0$. 
Wilkens concluded \cite{Wi3} with a conjecture concerning the inertia groups of highly connected
manifolds of dimension $7$ and $15$ which we have been unable to prove or disprove.

\chapter{Linking forms and linking functions} \label{chaplink}
We have now classified highly connected manifolds in dimension $7$ and $15$ up to almost
diffeomorphism by their quadratic linking functions which need not be homogeneous.  Thus the
existence of bundles of Hopf invariant one leads us to a somewhat more general class of quadratic
object in order to classify manifolds.  In this chapter we review the classifications of symmetric
bilinear forms on finite abelian groups started by Wall and finished by Kawauchi and Kojima.  We
next recall Nikulin's classification of quadratic linking forms refining a given bilinear form and
then generalize Nikulin's result to all quadratic linking functions.  Throughout this chapter $G$
shall denote a finite abelian group.

\section{Bilinear forms on finite abelian groups} \label{secbil}
Recall that a nondegenerate symmetric bilinear form $b$ on a
finite abelian group $G$ is called a {\em linking form} $(G,b)$.  Two linking forms $(G_0,b_0)$ and
$(G_1,b_1)$ are isometric if there is an isomorphism of abelian groups $\theta: G_0 \cong G_1$
such that $b_0 = b_1 \circ (\theta \times \theta)$.  The block sum of two linking forms $(G_1,b_1)$
and $(G_2,b_2)$ is defined in the obvious way to be $(G_1 \oplus G_2,b_1 \oplus b_2)$.  If $H$ is a
sub-group of $G$ we let $H^{\perp}$ denote the annihilator of $H$, 
$$H^{\perp} = \{g \in G | b(g,h) = 0 ~\forall h \in H\}.$$  In Section \ref{secblf} we promised
a proof of the following fundamental lemma which (according to Wall \cite{Wa3}) is well know
and may be found in \cite{Wa3}.
  
\begin{Lemma} \label{lemmabsplitt}
Let $(G,b)$ be a linking form and $H$ a subgroup of $G$ such that $b|_{H\times H}:H \times H \ra
\Q/\Z$ is nonsingular.  Then $(H^{\perp},b|_{H^{\perp}\times H^{\perp}})$ is nonsingular and
$(G,b)$  splits as the block sum
$(H,b|_{H\times H})
\oplus (H^{\perp},b|_{H^{\perp}\times H^{\perp}})$.
\end{Lemma}

\begin{proof}
Let $R_H:G^{\wedge} \ra H^{\wedge}$ be the map which restricts the domain of an element in the
torsion dual of $G$ to the subgroup $H$.  Now consider the composition of the adjoint map of $b$ with
$R_H$, 
$$R_H \circ \hat{b}:G \ra H^{\wedge}.$$  By definition $H^{\perp}$ is the kernel of
$R_H \circ \hat{b}$ and by hypothesis $(R_H \circ \hat{b})|_H$ is an isomorphism.  Thus
$H \cap H^{\perp} = \{e\}$ and $H$ and $H^{\perp}$ generate $G$ which is therefore the direct sum of
$H$ and $H^{\perp}$.  As $b(h,h') = 0$ for all $h \in H$ and $h' \in H^{\perp}$ the form $b$ splits
over $H \oplus H^{\perp}$ and it follows that $b|_{H^{\perp}}$ must be nonsingular for if it were
singular then $b$ would also be singular.
\end{proof}

Let $p$ be a prime.  Every linking form splits as a block sum over the $p$-primary components of
its domain.  We therefore let ${\mathcal{N}}_p$ denote the semi-group of linking forms on
$p$-groups under block sum and $\NN$ denote abelian semi-group of all linking forms under block
sum.  Hence
${\mathcal{N}} \cong~{\bigoplus}{_p}{\mathcal{N}}_p$.  The structure of linking forms for
$p$ an odd prime is relatively simple.  In this case every linking form splits as the block sum of
linking forms over cyclic groups $\Z/p^k$.  The following notation comes from Wall \cite{Wa3} and
Kawauchi and Kojima \cite{KK}.  The form, $A_p^k(\theta)$ denotes the unique (up
to isometry) linking form $(\Z/p^k,b)$ with $b(e,e)=\theta$ where $e$ is a generator of $\Z/p^k$. 
When $p$ is understood we write simply $A^k(\theta)$.  If $p$ is an odd prime there are two
isometry classes of linking form over $\Z/p^k$ depending on whether $\theta$ is a quadratic residue
$( \modu ~p)$ or not.  The semi-group ${\mathcal{N}}_p$ is generated by the linking forms
$A^k(1)$ and  $A^k(\theta), \theta$ a quadratic non-residue where $k$ ranges over the natural
numbers.  All relations amongst the linking forms are generated by the following isometry
$$A^k(1) \oplus A^k(1) \cong A^k(\theta) \oplus A^k(\theta).$$  Thus for odd primes the
Grothendieck group of $\NN_p$, ${\mathcal{G}}(\NN_p)$, is isomorphic to the group
$\Z^{\infty} \oplus (\Z/2)^{\infty}$,
where the infinite summands are generated by $[A^k(1)]$ and the summands of order two are generated
by $[A^k(1)]-[A^k(\theta)]$.

The main result of \cite{KK} was the computation of the abelian semi-group of $2$-primary linking
forms, $\NN_2$.  The generators had already been discovered by Wall \cite{Wa3} however the precise
relations between them were unknown.  We now review the generators of $\NN_2$.  There are the
cyclic linking forms $A^k(n)$ for $n$ odd.  In this case $A^1(n) = A^1(1)$ for all $n$,
$A^2(n) = A^2((-1)^{(n-1)/2})$ and $A^k(n) \cong A^k(n')$ if and only if $n \equiv n' (8)$ for all
$k
\geq 3$.  There are also the hyperbolic and pseudo-hyperbolic linking forms with matrices
$$E^{k,0} = \left(
\begin{array}{cc}
0 & 2^{-k}\\
2^{-k} & 0\\
\end{array} \right) ~(k \geq 1)
\hskip 1 cm  E^{k,1} = \left(
\begin{array}{cc}
2^{1-k} & 2^{-k}\\
2^{-k} & 2^{1-k}
\end{array} \right) ~(k \geq 2).
$$

\begin{Theorem}[Kawauchi and Kojima \cite{KK} Theorem 0.1] \label{thmkk01}
$\NN_2$ has a presentation with generators $A^k(n)$, where $n$ is $1 (k=1), \pm1(k=2),\pm1
\text{~or~} \pm5 (k
\geq 3)),E^{k,0} (k\geq1)$ and $E^{k,1} (k\geq2)$ and relations\\
$$\begin{array}{lcrr}
A^k(n_1) \oplus A^k(n_2) = A^k(n_1+4) \oplus A^k(n_2+4) & (k\geq3) & & (0.1)\\
A^k(n) \oplus 2A^k(-n) = A^k(-n) \oplus E^{k,0} & (k \geq 1) & &(0.2) \\
3A^k(n) = A^k(-n+4) \oplus E^{k,0} & (k\geq2) & & (0.3)\\
2E^{k,0} = 2E^{k,1} & (k\geq2) & & (0.4)\\
A^k(n_1) \oplus A^{k+1}(n_2) = A^k(n_1 + 2n_2) \oplus A^{k+1}(n_2 + 2n_1) & (k\geq1) & & (1.1)\\
A^k(n) \oplus E^{k+1,1} = A^k(n+4) \oplus E^{k+1,0} & (k\geq1) & & (1.2)\\
E^{k,1} \oplus A^{k+1}(n) = E^{k,0} \oplus A^{k+1}(n+4) & (k\geq2) & & (1.3)\\
A^k(n_1) \oplus A^{k+2}(n_2) = A^k(n_1+4) \oplus A^{k+2}(n_2+4) & (k\geq1) & & (2.1)\\
\end{array}$$
\end{Theorem}

\noindent
We call a finite abelian group $G$, as opposed to a linking form, {\em
homogeneous} if $G$ is isomorphic to $(\Z/p^k)^l$ for some integer $l$.  It is implicit in the
statement of the statement of Theorem \ref{thmkk01} that every linking form $b$ is isomorphic to a
direct sum of linkings over homogeneous groups which is fact that was already proven in
\cite[Lemma 8]{Wa3}.  A complete list of homogeneous linking forms is given in \cite[34-35]{KK}. 
Kawauchi and Kojima also provide a complete set of invariants which classify linking forms and we
review these now.  First let $\bar{G}_p^k$ be the subgroup of $G$ generated by elements of order
$p^s$, $s \leq k$.  Then let $\widetilde{G}_p^k := \bar{G}_p^k/(\bar{G}_p^{k-1} +
p\bar{G}_p^{k+1})$.  The group $\widetilde{G}_p^k$ is a vector space over $\Z/p$ and we define 
$r^k_p(G) = \dim_{\Z/p}(\widetilde{G}_p^k)$.  Thus each linking form $(G,b)$ is isomorphic to the
direct sum $\bigoplus_{k\geq 1} (G_p^k,b^k_p)$ with each $G_p^k \cong (\Z/p^k)^{r_p^k(G)}$.  The
function

$$\begin{array}{cccc}
\widetilde{b}^k_p : & \widetilde{G}_p^k \times \widetilde{G}_p^k & \ra & \Q/\Z \\
& ([x],[y]) & \mapsto & p^{k-1}b(x,y).\\
\end{array}$$
is well defined and indeed defines a nonsingular linking form on $\widetilde{G}_p^k$.
For each linking form $(G,b)$ we make the following

\begin{Definition}[\cite{KK} Definition 1.1] \label{defKK}
Let $c^k(b) \in \widetilde{G}^k_2$ be the characteristic element of $\widetilde{b}^k_2$ which is
uniquely defined by the identity
$$\widetilde{b}_2^k(c^k(b),[x]) = \widetilde{b}^k_2([x],[x]).$$
\end{Definition}
\noindent
Observe that $c^k(b) \neq 0$ if and only if $(G,b)$ contains an indecomposable cyclic
linking form $A^k(n_k)$ as a summand.  We now wrap up a lemma, corollary and definition of
Kawauchi and Kojima into a single proposition.

\begin{Proposition}[\cite{KK} \S 2]
Let $(G,b)$ be a linking form with $c^k(b) = 0$ and let $q_k$ be the (well defined) function 
$$\begin{array}{cccc}
q_k: & G_2/\bar{G}^k_2 & \ra & \Q/\Z \\
& [x] & \mapsto & 2^{k-1}b(x,x). \\
\end{array}$$
The complex number
$$GS_k(b) : = \sum_{[x]\in G_2/\bar{G}^k_2} {\rm exp} (2 \pi i q_k([x])) $$
(where $\expo$ denotes the standard exponential and $i = \sqrt{-1}$) 
is non-zero and there is an integer $\sigma$ defined uniquely ${\rm mod}~8$ such that $GS_k(b) =
|GS_k(b)|\cdot \exp(2\pi i \sigma/8)$.  The invariant $\sigma_k(b)$ is defined to be $\sigma_k(b)
\in \Z/8$ if $c^k(b) = 0$ and $\infty$ if $c^k(b) \neq 0$.  
\end{Proposition}

\begin{Theorem}[\cite{KK} 4.1] \label{thmkk}
Two linking forms $(G_2,b)$ and $(G_2^{'},b')$ are isometric if and only if $r^k_2(G_2) =
r^k(G'_2)$ and $\sigma_k(b) = \sigma_k(b')$ for all $k \geq 1$.
\end{Theorem}

\begin{Corollary}
Let $r^k_p:{\mathcal{G}}(\NN_2) \lra \Z$ be the $k$th-rank homomorphism defined on the
Grothendieck group of $\NN_2$ by
$$
\begin{array}{cccc}
r_2^k:{\mathcal{G}}(\NN_2) & \lra & \Z[k] \\
([b_0],[b_1]) & \mapsto & r_2^k(b_0) - r_2^k(b_1).
\end{array}
$$
Then the homomorphism $r_2$ which is the sum of $r_2^k$ for all positive integers $k$ gives an
isomorphism
$$r_2: \bigoplus_{k=1}^{\infty}r_2^k : {\mathcal{G}}(\NN_2) \cong~ \bigoplus_{k=1}^{\infty}
\Z[k].$$
\end{Corollary}

\begin{proof}
The homomorphism $r_2$ is clearly surjective.  Let $b_0$ and $b_1$ be linking forms for which
$([b_0],[b_1])$ lies in the kernel of $r_2$.  If $b$ is
the linking form
$$b=\bigoplus_{\{ k:r_2^k(b_0) \neq 0\} } A^k(1)$$
then $\sigma_k(b_0 \oplus b) = \sigma_k(b_1 \oplus b) = \infty$ for every positive integer $k$. 
By Theorem \ref{thmkk}, $b_0 \oplus b \cong b_1 \oplus b$ and hence the kernel of $r_2$ is trivial.
\end{proof}

\begin{remark}
Notice that cancellation fails in $\NN_2$ and hence information is lost in passing to
$\mcal{G}(\NN_2) \cong (\NN_2 \oplus \NN_2)/\Delta$, (where $\Delta = \{([q],[q]) \in \NN_2 \times
\NN_2 \}$).  In particular, if 
$$\bar{\Delta}= \{([q_0],[q_1]):(j[q_0],j[q_1]) \in \Delta \text{~~for some $j \in \Z^{+}$} \},\\$$
then there is an isomorphism of {\em sets}
$$\bar{\Delta}/\Delta \cong \Z/2^{\infty} \cup \Z/4^{\infty}.$$ 
\end{remark}

\section{Quadratic forms on finite abelian groups} \label{secquad}
Recall from Section \ref{secqlf} that a homogeneous quadratic refinement of a linking form $(G,b)$
is a function $q:G \ra \Q/\Z$ such that $q(-x) = q(x)$ for all $x \in G$ and such that 
$$b(x,y) = q(x+y) - q(x) - q(y) ~~ \in \Q/\Z$$
for all $x,y \in G$.  We note that many authors including Wall consider quadratic refinements in
$\Q/2\Z$ in which case the defining equation is
$$\bar{q}(x+y) = \bar{q}(x) + \bar{q}(y) + 2b(x,y)~~ \in \Q/2\Z$$
where multiplication by $2$ induces an isomorphism $\cdot 2:\Q/\Z \cong \Q/2\Z$.  The two
definitions are equivalent with $\bar{q} = (\cdot 2) \circ q$.  We call any such $q$ a
quadratic linking form.  Two quadratic linking forms $(G_0,q_0)$ and
$(G_1,q_1)$ are isometric if there exists and isomorphism of abelian groups $\theta: G_0 \cong
G_1$ such that $q_0 = q_1 \circ \theta$.  Note that any isometry of quadratic refinements is
necessarily an isometry of the associated bilinear forms.  If $G$ contains no $2$-torsion
then every bilinear form on $G$ has a unique quadratic refinement and so we shall be concerned with
bilinear forms on $2$-groups.  Given that Kawauchi and Kojima have classified bilinear forms on
$2$-groups we shall focus on the question of when two homogeneous quadratic refinements of the same
linking form are isometric.  Recall from Lemma \ref{lemmaqlfu2} {\em 4}, that if $q$ is a fixed
quadratic linking form refining $b$, then all other quadratic linking forms refining $b$ are of
the form $q_a$, where $a \in G$ and $2a=0$.  Thus our goal can be restated as
follows: given $a, c \in \GE$ when is $q_a \cong q_c$?  We shall discover a most satisfactory
answer to this question in Corollary \ref{corhqf}.  We turn now to the definition of the
invariants which resolve this question.  Recall Wall's theorem that every quadratic form $(G,q)$
is the boundary of some even quadratic form
$\kappa(H,\lambda,0)$, 
$$0 \lra H \stackrel{\hat{\lambda}}{\lra} H^* \ra G \lra 0.$$ 

\begin{Lemma} \label{lammasigmawd}
Let $q$ be a quadratic linking form and $\kappa(\lambda,0)$ any nondegenerate even quadratic
function such that $q^{ev}(\kappa) = q$.  The signature of $\lambda~(\modu~ 8)$ is independent of
the choice of $\kappa(\lambda,0)$ and is thus an invariant of $q$: 
$$\sigma(q) := \sigma(\lambda) ~ (\modu ~ 8)~ \in \Z/8.$$
\end{Lemma}

\begin{proof}
We use the technique of algebraic gluing from Lemma \ref{lemmaglue}.  Let $\kappa_0 =
\kappa(\lambda_0)$ and $\kappa_1 = \kappa(\lambda_1)$ be two nondegenerate quadratic forms
inducing $q$, then the quadratic form $(\kappa_0 \cup_q \kappa_1^{-})$ defined in Lemma
\ref{lemmaglue} is nonsingular and even and hence has signature divisible by $8$ according to
\cite[Theorem 5.1]{MH}.  If
$\lambda$ is the bilinear form associated to $(\kappa_0 \cup_q \kappa_1^{-})$ then $\lambda
\tensor Id_{\Q} = (\lambda_0 \tensor Id_{\Q}) \oplus (-\lambda_1 \tensor Id_{\Q})$ and hence
$\sigma(\lambda) = \sigma(\lambda_0) - \sigma(\lambda_1)$ and $\sigma(q)$ is well-defined.
\end{proof}

\begin{Theorem}[Nikulin \cite{Ni}] \label{thmnik}
Two quadratic linking forms $q$ and $q'$ are isometric if and only if $b(q) \cong b(q')$ and
$\sigma(q) \equiv \sigma(q') ( \modu ~ 8)$.
\end{Theorem}

\noindent
The invariant $\sigma(q)$ may also be computed by a Gauss-sum formula as follows.

\begin{Definition}
Let $(G,q)$ be a quadratic linking form on a finite abelian group.  Define the Gauss sum invariant
of $q$ by
$$GS(q) = 1/\sqrt{|G|}\cdot\sum_{x \in G} {\rm exp} (2\pi iq(x)) \in \C.$$ 
\end{Definition}

\begin{Theorem}[Milgram \cite{Milg}] \label{thmmilg}
Let $(G,q)$ be a quadratic linking form and let $\kappa(\lambda,0)$ be a
nondegenerate quadratic function such that $q^{ev}(\kappa) = q$.  Then
$$GS(q) =  {\rm exp}(2\pi i \sigma(\lambda)/8).$$
\end{Theorem}
\noindent
Theorem \ref{thmmilg} gives another proof that $\sigma(q): = \sigma(\lambda) ~(\modu ~ 8)$ is well
defined.  It follows in particular from Theorem \ref{thmmilg}, that $GS(q)$ is never zero. 
Moreover, the important part of this complex number for our purposes is its argument, so we make
the following

\begin{Definition} \label{defKA}
Let $(G,q)$ be a quadratic linking form.  Define the Kervaire-Arf
invariant of $q$, $K(q) \in \Q/\Z$, to be the argument of $GS(q)$ divided by $2\pi$.
\end{Definition} 

\begin{Corollary} \label{corhqf}
Two quadratic linking forms $q$ and $q'$ are isometric if and only if $b(q) \cong b(q')$ and $K(q)
= K(q')$.  For a fixed quadratic linking form $q$ and for any $a,c \in \bar G^1_2$, $q_a$ and $q_c$
are isometric if and only if $q(a) = q(c)$.
\end{Corollary}

\begin{proof}
The first statement is a trivial consequence of Theorems \ref{thmnik} and 
\ref{thmmilg}.  The second statement follows from the following simple calculation.
$$\begin{array}{cl}
GS(q_a) & = \sum_{x \in G} \expo(2\pi i q_a(x))\\
& = \sum_{x \in G} {\rm exp}(2\pi i (q(x+a) - q(a)))\\
& = \left[\sum_{x \in G} {\rm exp}(2\pi i q(x))\right].{\rm exp}(-2\pi i q(a))\\
& = GS(q).{\rm exp}(-2\pi i q(a)).
\end{array}
$$
Thus $K(q_a) = K(q) - q(a)$ and so $q_a \cong q_c$ if and only if and only if $q(a)
= q(c)$.
\end{proof}

\noindent
We shall not give Nikulin's proof of Theorem \ref{thmnik} here and mention only that it is achieved
by induction over the number of generators of $G$ and the classification of bilinear forms on
$2$-groups (which Nikulin seems to have obtained independently from Kawauchi and Kojima).  To
conclude this section we mention another result of Wall's which reduces the problem of classifying
quadratic forms to the classification of bilinear forms.  If $G$ is a finite $2$-group choose a
basis for $G$, that is a set of generators
$\{x_i\}$ for cyclic groups of which $G$ is a direct sum.  Let $G_{\times 2}$ be the group with
generators $\{x_i'\}$ where the order of $x_i'$ is twice the order of $x_i$.  There is a short
exact sequence of groups 
$$0 \ra \bar{G}_{\times 2}^1 \ra G_{\times 2} \ra G \ra 0.$$

\begin{Theorem}[Wall \cite{Wa3} Theorem 5] \label{thmwallG2}
A quadratic linking form $q:G \ra \Q/\Z$ determines and is determined by a linking form
$B(q):G_{\times 2} \times G_{\times 2} \ra \Q/\Z$.
\end{Theorem}

\noindent
With this result of Wall it is not difficult to use the classification of Kawauchi
and Kojima to prove Nikulin's theorem.  In fact, let $b(q): G \times G \ra \Q/\Z$ be the bilinear
form which $q$ refines.  By definition $GS_k(B(q)) = GS_{k-1}(b(q))$ for $k>1$ and $GS_1(B(q)) =
GS(q)$.  (Note that  since $G_{\times 2}$ has no elements of order two $c^1(B(q)) = 0$ and the
results of \cite{KK} give an alternate proof that $GS(q)$ is non-zero.)  Theorem \ref{thmnik} now
follows from Theorems \ref{thmwallG2} and \ref{thmkk} applied to $B(q)$ and the fact that linking
forms on groups of odd order have unique quadratic refinements.  In fact we may now state the
following extension of Theorem \ref{thmkk} to quadratic linking forms.

\begin{Corollary} \label{corkkq}
Two quadratic linking forms $(G_2,q)$ and $(G_2^{'},q')$ are isometric if and only if $r^k_2(G_2) =
r^k(G'_2)$ and $\sigma_k(B(q)) = \sigma_k(B(q'))$ for all $k \geq 1$.  
\end{Corollary}

\begin{proof}
By Theorem \ref{thmwallG2} $q$ and $q'$ are isometric if and only if $B(q)$ and $B(q')$ are
isometric and this occurs according to Theorem \ref{thmkk} if and only if $r^k_2(G_2) =
r^k(G'_2)$ and $\sigma_k(B(q)) = \sigma_k(B(q'))$ for all $k \geq 1$.
\end{proof}

\begin{remark}
Let $q: G \ra \Q/\Z$ be a quadratic linking form defined on the $2$-group $G$ and let the range of
$q$, $R(q)$, be the images of $q$ in $\Q/\Z$ counted with multiplicity.  It follows from Corollary
\ref{corkkq} that if $(G,q')$ is another quadratic liking form then $q$ is isometric to $q'$ if
and only if $R(q) = R(q')$.  
\end{remark}

\section{Quadratic functions on finite abelian groups} \label{secquadf}
A quadratic linking functions is a function $q:G \ra \Q/\Z$ which refines a bilinear form $b$ but
which need not be homogeneous.  Just as in the
homogeneous case we may define the Gauss-sum invariant and the associated Kervaire-Arf invariant.  
$$GS(q) = \sum_{x \in G} {\rm exp}(2\pi iq(x)) \in \C \hskip 2cm K(q) = Arg(GS(q))/2\pi$$
However the signature invariant of a general quadratic linking function must be somewhat
modified as follows.  Recall that for every quadratic linking function $q$, Corollary
\ref{corcharrealise} guarantees a characteristic quadratic function $\kappa$ such that
$q^c(\kappa) \cong q$.

\begin{Proposition} \label{props_1wd}
Let $(G,q)$ be a quadratic linking function and let $\kappa_0 = \kappa(\lambda_0,\alpha_0)$ be a
nondegenerate quadratic function such that $\delta^c(\kappa_0) = q$.  Then the quantity
$$\bar{s}_1(q) := (\lambda_0^{-1}(\alpha_0,\alpha_0) - \sigma(\lambda_0))/8~ \in \Q/\Z.$$
is independent of the choice of $\kappa_0$.
\end{Proposition}

\begin{proof}
Suppose that $\kappa_0$ and $\kappa_1$ are two quadratic functions such that $\delta^c(\kappa_0)
= q = \delta^c(\kappa_1)$.  Then by Lemma \ref{lemmaglue} we may form the nonsingular
characteristic quadratic function $\kappa_0~\cup_{Id}~\kappa_1^-$.  We let $\lambda$ and $\alpha$
denote respectively the bilinear and linear parts of $\kappa_0\cup_{Id}\kappa_1^-$.
Since $\alpha$ is characteristic for $\lambda$ it is well known that
$\lambda(\alpha,\alpha) \equiv \sigma(\lambda)$ (mod $8$) (see for example \cite {MH}).  However
$$\begin{array}{cl}
[\lambda(\alpha,\alpha) - \sigma(\lambda)]/8 ~ \modu ~ \Z ~ = [\lambda_0^{-1}(\alpha_0,\alpha_0)
- \lambda_1^{-1}(\alpha_1,\alpha_1) -(\sigma(\lambda_0) - \sigma(\lambda_1))]/8 ~ \modu ~ \Z\\
\hskip 1cm 0 \hskip 1cm  \modu ~ \Z ~ = [\lambda_0^{-1}(\alpha_0,\alpha_0) -\sigma(\lambda_0)]/8
-[\lambda_1^{-1}(\alpha_1,\alpha_1) -\sigma(\lambda_1)]/8 ~ \modu ~ \Z.
\end{array}$$
As the right hand side is the difference of the possible values for $\bar{s}_1(q)$ we have proven
the proposition.
\end{proof}

\begin{remark} \label{rems1}
We have called this invariant of a quadratic linking functions $\bs$ because when $q =
\delta^c(\kappa_0)$ then $\bs(q) = \bs(\del L(\kappa_0))$ where $\bs$ is the almost diffeomorphism
invariant of $\del L_0$ which was defined in Section \ref{chapprelim}.\ref{secs_1}.
\end{remark}

\begin{Proposition} \label{props1}
Let $q:G \ra \Q/\Z$ be a quadratic linking function on a finite abelian group.  Then
$$\bar{s}_1(q) = - K(q)$$
\end{Proposition}

\begin{proof}
By the proof of Corollary \ref{corcharrealise} there is a nondegenerate characteristic quadratic
function
$\kappa = \kappa(H,\lambda,\alpha)$ with $\lambda$ an even bilinear form such that $q =
\delta^c(\kappa)$.  It follows that $2$ divides $\alpha$ and so we let $a =[\alpha/2] \in G =
H^*/\hl(H)$.  Moreover, if we let $q^o = \delta^{c}(\kappa(\lambda,0))$ be the quadratic
linking form defined by $\lambda$ alone then $q = q^o_a$ by Lemma \ref{lemmaqlfpert}.  By
Proposition \ref{props_1wd}
$$\begin{array}{cl}
\bar{s}_1(q) & = [\lambda^{-1}(\alpha,\alpha) - \sigma(\lambda)]/8 ~ \modu ~ \Z \\
& = \lambda^{-1}(\alpha/2,\alpha/2)/2 - K(q^o) \\
& = q'(a) - K(q^o)\\
& = -K(q).
\end{array}$$
The last step follows by the arithmetic of Corollary \ref{corhqf} which showed that
$K(q^o_a) = K(q^o) - q^o(a)$ (and which made no use of the fact operating then that $2a=0$).  
\end{proof}

\begin{Corollary}
If $q_0$ and $q_0$ are two quadratic linking functions, then 
$K(q_0 \oplus q_1) = K(q_0) + K(q_1)$.
\end{Corollary}

\begin{proof}
This will follow from the simple fact that $\bar{s}_1$ is additive.  For let
$\kappa_0(\lambda_0,\alpha_0)$ and $\kappa_1(\lambda_1,\alpha_1)$ be two quadratic functions such
that $q_i = \delta^c(\kappa_i)$.  Then
$$\begin{array}{cl}
\bar{s}_1(q_0 \oplus q_1) &= [(\lambda_0 \oplus \lambda_1)^{-1}(\alpha_0,\alpha_1) -
\sigma(\lambda_0 \oplus \lambda_1)]/8 ~ \modu ~ \Z \\
&= [\lambda_0^{-1}(\alpha_0,\alpha_0) - \sigma(\lambda_0)]/8 + [\lambda_1^{-1}(\alpha_1,\alpha_1) -
\sigma(\lambda_1)]/8 ~ \modu ~ \Z \\
&= \bar{s}_1(q_0) + \bar{s}_1(q_1))
\end{array}$$
\end{proof}

Now recall the Wilkens data $(b(q),\beta(q))$ of a quadratic linking function $q$ which is the
pair composed of the linking form refined by $q$ and the even element $\beta(q)$ defined in
\ref{lemmaqlfu2}, {\em 5}.  The element $\beta(q)$ is twice the linear part of $q$.

\begin{Lemma} \label{lemmaGodd}
Let $q:G \ra \Q/\Z$ be a quadratic linking function.  If the order of $G$ is
odd then $(b(q),\beta(q))$ is a complete invariant of $q$.
\end{Lemma}

\begin{proof}
In a group of odd order there is a unique quadratic linking form $q^o$ refining $b(q)$.  As there
is a unique $\gamma \in G$ such that $q=q^o_{\gamma}$ and $\gamma = (1/2)\beta(q)$ is determined by
$\beta(q)$, we see that $(b(q),\beta(q))$ determines $q$.
\end{proof}

Before stating and proving the classification of
quadratic linking functions we establish notation for quadratic linking functions parallel to the
notation for linking forms in Section \ref{secbil}.  Let $q^k(n)$, $q^{k,0}$ and $q^{k,1}$ denote
homogeneous quadratic refinements of $A^k(n)$, $E^{k,0}$ and $E^{k,1}$ respectively.  We define 
homogeneous quadratic refinements $q^k(n)$, $q^{k,0}$, and $q^{k,1}$ of $A^k(n)$, $E^{k,0}$, and
$E^{k,1}$ as follows.  For $r,s \in \Z$,
$$\begin{array}{c}
(k \geq 1)~A^k(n) = (\Z/2^k[e],A^k_n(e,e) = n2^{-k}) \\ q^k(n)(re) = r^2 n 2^{-k-1}\\
\end{array}$$
$$\begin{array}{c}
(k\geq 1)~E^{k,0} = (\Z/2^k[e_1] \oplus \Z/2^k[e_2], ~~E^{k,0}(e_i,e_i) = 0, ~~E^{k,0}(e_1,e_2) =
2^{-k}) \\ 
q^{k,0}(r e_1 + s  e_2) = rs2^{-k}
\end{array}$$
$$\begin{array}{c}
(k\geq 2)~E^{k,1} = (\Z/2^k[e_1] \oplus \Z/2^k[e_2], ~~E^{k,1}(e_i,e_i) = 2^{1-k},
~~E^{k,1}(e_1,e_2) = 2^{-k})
\\ q^{k,1}(r e_1 + s e_2) = (r^2 + rs +s^2)2^{-k}.\\
\end{array}$$

\noindent
We shall continue to with our convention of using a subscript to denote a linear perturbation of a
quadratic function, $q_a(x) = q(x) + b(x,a)$.  So for example 
$$[q^k(n)_{e}](re) = q^k(n)(re) + A^k(n)(re,e) = r^2 n 2^{-k-1} + rn2^{-k}.$$

\begin{Theorem}\label{thmqf}
Let $q_0$ and $q_1$ be quadratic linking functions over a finite abelian group $G$.  
Then $q_0$ is isometric to $q_1$ if and only if 
$(b(q_0),\beta(q_0)) \cong (b(q_1),\beta(q_1))$ and $K(q_0) = K(q_1)$.
\end{Theorem}

\begin{proof}
By Lemma \ref{lemmaGodd}, $(b(q),\beta(q))$ is an invariant of the isometry class of any
quadratic linking function $q$.  Moreover, $K(q)$ is manifestly an invariant.  Therefore,
it remains to show that $q_0$ and $q_1$ are isometric if $(b(q_0),\beta(q_0)) \cong
(b(q_1),\beta(q_1))$ and $K(q_0) = K(q_1)$.  By Lemma
\ref{lemmaGodd} (or by the results of Wilkens' thesis) the theorem is true if $b(q_0)$ and
$b(q_1)$ are linking forms on a group of odd order.  We therefore assume that the linking forms
are defined on
$2$-primary groups.  Let $q$ have Wilkens data $(b(q),\beta(q))$ and a fixed Kervaire-Arf invariant
$K(q) \in \Q/\Z$.  The strategy of the proof shall be to show that $q$ is uniquely determined up to
isometry by these data using induction on the number of indecomposable summands in $b(q)$.

We first observe that if $\theta: (b(q_0),\beta_0)) \cong (b(q_1),\beta(q_1))$ then we may
replace $(b(q_1),\beta(q_1))$ with $(b(q_1)\circ(\theta\times \theta),\theta^{-1}(\beta(q_1))$.  We
shall therefore always assume that $(b(q_0),\beta(q_0)) = (b(q_1),\beta(q_1))$.  In order to prove
the existence of the isometries of quadratic linking functions predicted by the theorem we shall
constantly use the following lemma.

\begin{Lemma} \label{lemmaKWisom}
Two quadratic linking functions $(G_0,q_0)$ and $(G_1,q_1)$ are isometric if and only if there are
bases $\{e_1, \dots, e_n \}$ and $\{f_1, \dots, f_n \}$ for $G_0$ and $G_1$ respectively such that
and $q_0(e_i) = q_1(f_i)$ and $b(q_0)(e_i,e_j) = b(q_1)(f_i,f_j)$ for all $i$ and $j$.
\end{Lemma}

\begin{proof}[Proof of Lemma]
Clearly the condition given is necessary for $q_0$ and $q_1$ to be isometric.  It is sufficient
because the relation $q_i(x+y) = q_i(x) + q_i(y) + b_i(x,y)$, $i=0,1$, means that $q_i$ is
determined by
$b(q_i)$ and the values $q_i$ takes on a basis.
\end{proof}

The base case of our induction is when $b(q)$ is indecomposable.  In this case Wilkens' thesis
implies that there are at most two quadratic linking functions with invariants $(G,b,\beta)$ and
one may check using Wilkens' calculations that in the cases where there are two quadratic linking
functions then they have distinct Kerviare-Arf invariant.  We now give a proof of the theorem in
the indecomposable case which does not rely on Wilkens' results.  Consider the case where
$(b(q),\beta(q)) = (A^k(n),\beta)$ is a linking form on the cyclic group $\Z/2^k$.  We must
therefore consider $q^k(n)_a$ and $q^k(n)_{a+2^{k-1}e}$ where $2a = \beta$.  We have already noted
that the proof of Corollary \ref{corhqf} shows that $K(q_a) = K(q) - q(a)$ for any
quadratic linking function.  Hence
$$\begin{array}{cl}
K(q^k(n)_{a+2^{k-1}e})-K(q^k(n)_a) & = q^k(n)(a+2^{k-1}e) - q^k(n)(a) \\
& = (na^2 + 2^kna + n2^{2k-2} - na^2).2^{-k-1} \\
& = na2^{-1} + n2^{k-3}.
\end{array}$$
We see that $K(q^k(n)_{a+2^{k-1}e}) \neq K(q^k(n)_a$ when $k=1$, when $k=2$ and $2|a$ and when
$k\geq 3$ and $2 \not  | \; a$ (recall that $n$ is odd).  In these cases the two quadratic
refinements are not isometric, this is detected by the Kervaire-Arf invariant and there is nothing
further to prove.  In the remaining cases ($k=2$ and $2 \not \! | \; a$ or $k\geq 3$ and $2 | a$)
we must show that
$q^k(n)_{a+2^{k-1}e}$ and
$q^k(n)_a$ are isometric.  We let $\theta:\Z/2^k \ra \Z/2^k$ be multiplication by $(1+2^{k-1})$
and verify that $\theta$ is the desired isometry.  Firstly let $k=2$ and $2 \nd a$.
$$\begin{array}{cl}
q^2(n)_a(3e) & = q^2(n)(e+2e) + A^2(n)(e+2e,a) \\
& = q^2(n)(e) + q^2(2e) + A^2(n)(e,2e) +A^2(n)(2e,a) + A^2(n)(e,a) \\
& = q^2(n)(e) + 2^{-1} +  A^2(n)(e,2e) + 2^{-1}+ A^2(n)(e,a)\\
& = q^2(n)(e) + A^2(n)(e,a+2e)\\
& = q^2(n)_{a+2e}(e)\\ 
\end{array}$$
Hence by Lemma \ref{lemmaKWisom} applied to $(q^k(n)_{a+2^{k-1}e},e)$ and $(q^k(n)_a,3e)$ we have
the desired isometry.  In the second case $k \geq 3$ and $2|a$.

$$\begin{array}{cl}
q^k(n)_a((1+2^{k-1})e) & = q^k(n)(e+2^{k-1}e) + A^k(n)(e+2^{k-1}e,a) \\
& = q^k(n)(e) + q^k(2^{k-1}e) + A^k(n)(e,2^{k-1}e) + A^k(n)(a,e + 2^{k-1}e) \\
& = q^k(n)(e) + 0  + A^k(n)(e,2^{k-1}e)  + A^k(n)(a,e)  \\ 
& = q^k(n)(e) + A^k(n)(e,a+2^{k-1}e)\\
& = q^k(n)_{a+2^{k-1}e}(e)\\
\end{array}$$
Again by Lemma \ref{lemmaKWisom} applied to $(q^k(n)_{a+2^{k-1}e},e)$ and $(q^k(n)_a,(1+2^k)e)$ we
have the desired isometry.  We now consider the hyperbolic and pseudo-hyperbolic linking forms. 
When $(b(q),\beta(q)) = (E^1_0,0)$ the classification of quadratic refinements via the Kervaire-Arf
invariant is well known: $K(q^1) = K(q^1_{e_1}) = K(q^1_{e_2}) = 0$ and $K(q^1_{e_1 + e_2}) =
1/2$.  It is perhaps most instructive simply to list the table of values for all four quadratic
refinements of $E^1$.
$$\begin{array}{ccccc}
& q^1 & q^1_{e_1} & q^1_{e_2} & q^1_{e_1+e_2}\\
0 & 0 & 0 & 0 & 0 \\
e_1 & 0 & 0 & 1/2 & 1/2 \\
e_2 & 0 & 1/2 & 0 & 1/2 \\
e_1 + e_2 & 1/2 & 0 & 0 & 1/2 \\
\end{array}$$
The first three quadratic linking functions are evidently isometric via permutations of $e_1, e_2$
and $e_1 + e_2$.  We now consider Wilkens data $(b(q),\beta(q)) = (E^{k,\epsilon},\beta)$ where
$\epsilon = 0,1$, $k \geq 2$, $\beta = 2a$ and $a = a_1e_1 + a_2e_2$ is any element in
$\Z/2^k \oplus \Z/2^k$.  We let $c=c_1e_1+c_2e_2$ denote any of the elements
$e_1, e_2$ or $e_1 + e_2$ so that $c_i = 0$ or $1$.  Again we apply Corollary \ref{corhqf} to
calculate that
$$\begin{array}{cl}
K(q^{k,\epsilon}_{a+2^{k-1}c}) - K(q^{k,\epsilon}_a) & = q^{k,\epsilon}(a+2^{k-1}c) -
q^{k,\epsilon}(a) \\
& = q^{k,\epsilon}(a) + q^{k,\epsilon}(2^{k-1}c) + E^{k,\epsilon}(a,2^{k-1}c) - q^{k,\epsilon}(a)
\\
& = E^{k,\epsilon}(a,2^{k-1}c) \\
& = (a_1c_1 + a_2c_2 )2^{-1}\\
& = \left\{ \begin{array}{rl}
0 & \text{~if~}(a_1c_1 + a_2c_2)\text{is even} \\
1/2 & \text{~if~} (a_1c_1 + a_2c_2)\text{is odd} \\
\end{array}  \right\}.
\end{array} $$
Hence we see Wilkens' criteria that when the order of the group is greater than $4, (k \geq 2)$,
that $(E^{k,\epsilon},\beta)$ is ambiguous if and only if $4$ does not divide $\beta = 2a_1e_1 +
2a_2e_1$.  More importantly we must exhibit an isometry $\theta$ from $q^{k,\epsilon}_a$ to
$q^{k,\epsilon}_{a+2^{k-1}f}$ when $a_1c_1 + a_2c_2$ is odd.  Now either both $a_i$ are odd in
which case we assume (as we may interchange bases elements at will) that $c_2 = 0$ or precisely one
$a_i$ is odd in which case we assume it to be $a_1$ and hence $c_1 = 1$.  We describe the necessary
isometry
$\theta$ in these two cases and present the matrices which by Lemma \ref{lemmaKWisom} ensures
that $\theta$ is an isometry.

$$(a_i \text {~odd}, c_1 = 1, c_2 = 0 ) \left\{ \begin{array}{ccc}
q^{k,\epsilon}_a & \stackrel{\theta}{\ra} & q^{k,\epsilon}_{a + 2^{k-1}f} \\
e_1 & \mapsto & e_1 \\
e_2 & \mapsto & e_2 + 2^{k-1}e_1.
\end{array} \right\}
$$
$$ \left[
\begin{array}{cc}
q^{k,\epsilon}_a(e_1) & E^{k,\epsilon}(e_1,e_2) \\
* & q^{k,\epsilon}_a(e_2)
\end{array}
\right] 
 =
\left[
\begin{array}{cc}
q^{k,\epsilon}_a(e_1) & 2^{-k} \\
* &  q^{k,\epsilon}_a(e_2)
\end{array}
\right]$$

$$ \left[
\begin{array}{cc}
q^{k,\epsilon}_{a+2^{k-1}c}(e_1) & E^{k,\epsilon}(e_1,e_2+2^{k-1}e_1) \\
* & q^{k,\epsilon}_{a+2^{k-1}c}(e_2+2^{k-1}e_1)
\end{array}
\right] 
 =
\left[
\begin{array}{cc}
q^{k,\epsilon}_a(e_1) + c_2 & 2^{-k} \\
* &  q^{k,\epsilon}_a(e_2)+2^{-1}(c_1+a_2)
\end{array}
\right]$$

$$(a_1 \text{~odd}, a_2 \text{~even}, c_1 = 1)
\left\{ \begin{array}{ccc}
q^{k,\epsilon}_a & \stackrel{\theta}{\ra} & q^{k,\epsilon}_{a + 2^{k-1}f} \\
e_1 & \mapsto & e_1 + 2^{k-1}c_2e_2 \\
e_2 & \mapsto & e_2 + 2^{k-1}e_2.
\end{array} \right\} $$

$$ \left[
\begin{array}{cc}
q^{k,\epsilon}_{a+2^{k-1}c}(e_1+2^{k-1}c_2e_2) & E^{k,\epsilon}(e_1,e_2+2^{k-1}c_2e_1) \\
* & q^{k,\epsilon}_{a+2^{k-1}c}(e_2+2^{k-1}e_1)
\end{array}
\right] 
 = $$
$$\left[
\begin{array}{cc}
q^{k,\epsilon}_a(e_1) + 2^{-1}(c_2 + c_2a_1)& 2^{-k} \\
* &  q^k_a(e_2)+2^{-1}(c_1+a_1)
\end{array}
\right]$$

This completes the classification of quadratic linking forms refining indecomposable pairs
$(b(q),\beta)$.  We call the Wilkens data $(b,\beta)$ ambiguous if there is more than one 
isometry class of quadratic linking function with Wilkens data $(b,\beta)$.  
We call the Wilkens data $(b,\beta)$ unambiguous if there is a unique quadratic linking function
with invariants $(b,\beta)$.   We have just reproven Wilkens' result that an indecomposable
pair $(b,\beta)$ is ambiguous if and only if the order of the group is less than or equal to $4$
and $\beta = 0$ or the order of the group is greater than $4$ and $\beta$ is not divisible by
$4$.  Moreover we note that in all ambiguous cases except for $(A^1(1),0)$ that the $K(q_0) -
K(q_1) = 1/2$ where $q_0$ and $q_1$ are the two possible refinements.  When $(b(q),\beta(q)) =
(A^1(1),0)$ then there are two possible quadratic refinements $q^1(1)$ and $q^1(1)_e$ and these
have distinct Kervaire-Arf invariants: 
$$K(q^1(1)) = \frac{1}{8}
\hskip 0.5cm \text{and} \hskip 0.5cm K(q^1(1)_e)=\frac{7}{8}. $$ 

The inductive step requires us to verify the theorem in a large number of special cases.  We start
with a simple case by way of illustration.  Suppose that $(b(q),\beta(q)) = (b_0,\beta_0) \oplus
(b',\beta')$ where $(b_0,\beta_0)$ is an unambiguous pair with quadratic refinement $q_0$.  It
follows that $q \cong q_0 \oplus q'$ where $q'$ refines $(b',\beta')$.  But since the Kervaire-Arf
invariant is additive over block sum $K(q') = K(q) - K(q_0)$ is determined and hence by induction
$q'$ is uniquely determined.  If suffices therefore to restrict our attention to pairs
$(b(q),\beta(q))$ consisting entirely of ambiguous irreducible components.  If these all have
$\beta(q) = 0$ then all the quadratic linking functions are homogeneous and hence quadratic
linking forms in which case Nikulin has proven the result.  We thus check the cases in which at
least one of the ambiguous indecomposable pairs has $\beta(q) \neq 0$ and therefore $4 \nd
\beta(q)$.  We label that the ambiguous indecomposable inhomogeneous Wilkens data $(I,\beta)$
and recall that they are as follows:
$$\begin{array}{ccc}
(A^k(n),\beta)~k>2,~4\nd \beta, & (E^{k,\epsilon},\beta)~k>1,~4\nd \beta. 
\end{array}$$
Suppose that $(b(q),\beta(q)) = (A^1(1),0) \oplus (I,\beta)$ then there are four possible
quadratic refinements and four different values of $K(q)$ so the result is proved.  In the
remaining cases we consider Wilkens data of the form $(I,\beta_0) \oplus (J,\beta_1) \oplus
(b',\beta')$ where
$(I,\beta_0)$ is ambiguous, indecomposable and inhomogeneous and $(J,\beta_1)$ is ambiguous and
either indecomposable or $(A^1(1) \oplus A^1(1), 0)$.  Let $q^I$
and $*q^I$ denote the two isometry classes of quadratic refinements of $(I,\beta_0)$ and let
$q^J$ and $*q^J$ be quadratic refinements of $(J,\beta_1)$ such that $K(q^J) - K(*q^J) = 1/2$.  We
shall show that $*q^I \oplus *q^J \cong q^I \oplus q^J$.  It follows that every quadratic
refinement of $(I,\beta_0) \oplus (J,\beta_1) \oplus (b',\beta')$ is isometric to $q^I \oplus
q^{''}$ for some $q^{''}$ which refines $(J,\beta_1) \oplus (b',\beta')$ and hence by induction
that the value of $K(q)$ determines the isometry class of any quadratic refinement of $(I,\beta_0)
\oplus (J,\beta_1) \oplus (b',\beta')$.  

We now proceed to the necessary computations.  Let $a_i$ be elements of the domain of $I$ or
$J$ such that $2a_i = \beta_i$.   From our computations in the indecomposable case we know that if
$I$ or $J$ is isomorphic to $E^{k,\epsilon}$ with $k>1$ then, by swapping $e_1$ and $e_2$ if
necessary, we may assume that $K(q^{k,\epsilon}_{a_i}) - K(q^{k,\epsilon}_{a_i+2^{k-1}e_1}) =
1/2$.  Finally, when considering pairs of similar linking forms we shall distinguish
generators by the place they hold.  So, for example, we shall denote generators of $A^k(n_0)
\oplus A^l(n_1)$ by
$(e,0)$ and
$(0,e)$.

$$\begin{array}{lr}
(b(q),\beta(q)) = (A^k(n_0),\beta_1) \oplus (A^l(n_1),\beta_2)  &
(k \geq 3 , l \geq 2)\\
\end{array} $$
$$\begin{array}{ccccccc}
q^k(n_0)_{a_0} & \oplus & q^l(n_1)_{a_1} & \ra & q^k(n_0)_{a_0+2^{k-1}e} & \oplus &
q^l(n_1)_{a_1+2^{l-1}e}\\
(e & , & 0) & \mapsto & (e & , & 2^{l-1}e)\\
(0 & , & e) & \mapsto & (2^{k-1}e & , & e)
\end{array}$$ \vskip 0.5cm

$$\begin{array}{lr}
(b(q),\beta(q)) = (A^k(n_0),\beta_1) \oplus (E^{l,\epsilon},\beta_2)  & (k,l \geq 2).\\
\end{array}$$ 
$$\begin{array}{ccccccc}
q^k(n_0)_{a_0} & \oplus & q^{l,\epsilon}_{a_1} & \ra & q^k(n_0)_{a_0+2^{k-1}e} & \oplus &
q^{l,\epsilon}_{a_1+2^{l-1}e_1}\\
(e & , & 0) & \mapsto & (e & , & 2^{l-1}e_1)\\
(0 & , & e_1) & \mapsto & (0 & , & e_1)\\
(0 & , & e_2) & \mapsto & (2^{k-1}e & , & e_2)
\end{array}$$ \vskip 0.5cm
\pagebreak
$$\begin{array}{lr}
(b(q),\beta(q)) = (E^{k,\epsilon},\beta_1) \oplus (E^{l,\epsilon},\beta_2))  & (k, l \geq 2).\\
\end{array}$$
$$\begin{array}{ccccccc}
q^{k,\epsilon}_{a_0} & \oplus & q^{l,\epsilon}_{a_1} & \ra & q^{k,\epsilon}_{a_0+2^{k-1}e} & \oplus
& q^{l,\epsilon}_{a_1+2^{l-1}e_1}\\
(e_1 & , & 0) & \mapsto & (e_1 & , & 0)\\
(e_2 & , & 0) & \mapsto & (e_2 & , & 2^{l-1}e_1)\\
(0 & , & e_1) & \mapsto & (0 & , & e_1)\\
(0 & , & e_2) & \mapsto & (2^{k-1}e_1 & , & e_2)
\end{array}$$ \vskip 0.5cm


$$\begin{array}{lr}
(b(q),\beta(q)) = (A^k(n_0),\beta_1) \oplus (E^{1,0},0)  & (k
\geq  3, 4\not| ~\beta_1) .\\
\end{array}$$
$$\begin{array}{ccccccc}
q^k(n_0)_{a_0} & \oplus & q^1_{e_1 + e_2} & \ra & q^k(n_0)_{a_0+2^{k-1}e} & \oplus &
q^1\\
(e & , & 0) & \mapsto & (e & , & e_1+e_2)\\
(0 & , & e_1) & \mapsto & (2^{k-1}e & , & e_1)\\
(0 & , & e_2) & \mapsto & (2^{k-1}e & , & e_2)
\end{array}$$ \vskip 0.5cm

$$\begin{array}{lr}
(b(q),\beta(q)) = (A^k(n_0),\beta_1) \oplus (A^1(1)\oplus A^1(1),0))  & (k
\geq 3, 4\not| ~\beta_1) .\\
\end{array}$$
$$\begin{array}{ccccccccccc}
q^k(n_0)_{a_0} & \oplus & q^1(1)_{a_1} & \oplus & q^1(1)_{a_2} & \ra & q^k(n_0)_{a_0+2^{k-1}e} &
\oplus & q^1(1)_{a_1+e} & \oplus & q^1(1)_{a_2+e}\\
(e & , & 0 & , & 0) & \mapsto & (e & , & e & , & e)\\
(0 & , & e & , & 0) & \mapsto & (2^{k-1}e & , & e & , & 0)\\
(0 & , & 0 & , & e) & \mapsto & (2^{k-1}e & , & 0 & , & e)
\end{array}$$ \vskip 0.5cm


\pagebreak
$$\begin{array}{lr}
(b(q),\beta(q)) = (E^{k,\epsilon},\beta_1) \oplus (E^{1,0},0))  & (k
\geq  2, 4\not| ~\beta_1) .\\
\end{array}$$
$$\begin{array}{ccccccc}
q^{k,\epsilon}_{a_0} & \oplus & q^1_{e_1 + e_2} & \ra & q^{k,\epsilon}_{a_0+2^{k-1}e_1} & \oplus
& q^1\\
(e_1 & , & 0) & \mapsto & (e_1 & , & 0)\\
(e_2 & , & 0) & \mapsto & (e_2 & , & e_1 + e_2)\\
(0 & , & e_1) & \mapsto & (2^{k-1}e_1 & , & e_1)\\
(0 & , & e_2) & \mapsto & (2^{k-1}e_1 & , & e_2)
\end{array}$$ \vskip 0.5cm

$$\begin{array}{lr}
(b(q),\beta(q)) = (E^{k,\epsilon},\beta_1) \oplus (A^1(1) \oplus A^1(1) ,0))  &
(k \geq  2, 4\not| ~\beta_1) .\\
\end{array}$$
$$\begin{array}{ccccccc}
q^{k,\epsilon}_{a_0} & \oplus & q^1_{a_1} & \ra & q^{k,\epsilon}_{a_0+2^{k-1}e} & \oplus
& q^1_{a_1+2^{l-1}e_1}\\
(e_1 & , & 0) & \mapsto & (e_1 & , & 0)\\
(e_2 & , & 0) & \mapsto & (e_2 & , & 2^{l-1}e_1)\\
(0 & , & e_1) & \mapsto & (0 & , & e_1)\\
(0 & , & e_2) & \mapsto & (2^{k-1}e_1 & , & e_2)
\end{array}$$ \vskip 0.5cm

\noindent
This concludes the proof of Theorem \ref{thmqf}.
\end{proof}

\begin{Corollary} \label{corwilkdata}
Let $P_0$ and $P_1$ be highly connected rational homology spheres of dimension $7$ or $15$ with
isometric Wilkens invariants $(G_0,b_0,\beta_0)$ and $(G_1,b_1,\beta_1)$.  Then $P_0$ and
$P_1$ are almost diffeomorphic (resp. diffeomorphic) if and only if $\bar{s}_1(P_0) =
\bar{s}_1(P_1)$ (resp. $s_1(P_0) = s_1(P_1)$ and $\Sigma_{P_0} = \Sigma_{P_1}$).
\end{Corollary}

\begin{proof}
By Theorem \ref{thmb}, $P_0$ and $P_1$ are almost diffeomorphic if and only if $Q(P_0) =
(G_0,q^c_0,\beta_0)$ is isometric to $Q(P_1) = (G_1,q^c_1,\beta_1)$.  As $P_0$ and $P_1$ are 
rational homology spheres $q_0^c=q_0$ and $q_1^c=q_1$ are quadratic linking functions which by
hypothesis define isometric Wilkens data.  Now by Theorem \ref{thmqf}, $q_0$ and $q_1$ are
isometric if and only if
$K(q_0) = K(q_1)$ but by Proposition \ref{props1}, $K(q_i) = -\bs(q_i)$, $i=0,1$, and by Remark
\ref{rems1}, $\bs(q_i) = \bs(P_i)$.  Hence $P_0$ and $P_1$ are almost diffeomorphic if
and only if $\bs(P_0) = \bs(P_1)$.  By Theorem \ref{thmsmoothclass}, $P_0$ and $P_1$ are
diffeomorphic if and only if $q_0$ is isometric to $q_1$, $s_1(P_0) = s_1(P_1)$ and 
$\Sigma_{P_0} = \Sigma_{P_1}$.  As $\bs$ is a multiple of $s_1$, we conclude that $P_0$ and $P_1$
are diffeomorphic if and only if $s_1(P_0) = s_1(P_1)$ and $\Sigma_{P_0} = \Sigma_{P_1}$.
\end{proof}

\chapter{Classification up to homeomorphism and homotopy} \label{chaphom}
In this final chapter we shall restrict our attention to highly connected manifolds of
dimension $7$ and thus the pair $(L,P)$ shall always denote an $8$-dimensional handlebody $L$ with 
boundary $P$.  We remind the reader that we consider only smooth manifolds $P$.  We shall extend the
almost diffeomorphism classification obtained in Chapter
\ref{chaptop} to a classification up to $PL$-homeomorphism, homeomorphism and homotopy.  The
classification of manifolds is often presented in the opposite direction, moving from the homotopy
classification to the homeomorphism classification via the surgery exact sequence and then upwards via
smoothing theory to the $PL$-classification and the smooth classification.  However, in the case of
highly connected manifolds it has proven easier to classify up to almost diffeomorphism first and then
descend to the other categories using smoothing theory and the surgery exact sequence in the opposite
direction.

\section{Smoothing theory} \label{secsmooththeory}
We begin by briefly recalling the theory of smoothing manifolds.  Suppose that $M$ is a $PL$
manifold of dimension greater than $4$.  If we embed $M$ in Euclidean space of sufficiently highly
dimension then the embedding we choose is unique up to isotopy and hence $M$ has a well defined
stable normal bundle which we denote by $\nu_M^{PL}: M \ra BPL$.  Now there is a forgetful map from
$\pi_{O,PL}:BO \ra BPL$ which classifies the universal smooth bundle $\gamma: EO \ra BO$ as a
$PL$-bundle. It is a remarkable fact (proven in \cite{HM}) that the question of whether $M$ may be
equipped with a smooth structure reduces to the question of whether the map stable $PL$-normal
bundle of $M$ lifts through $\pi_{O,PL}$.
\[
\begin{diagram}
\node[2]{BO} \arrow{s,r}{\pi_{O,PL}} \\
\node{M} \arrow{e,b}{\nu_M^{PL}} \arrow{ne,..}
 \node{BPL}
\end{diagram}
\]
In fact, if $\nu$ is a lift of $\nu_M^{PL}$ then there is a smooth structure on $M_{\nu}$
for which $\nu$ classifies the stable normal bundle.  We would now like to determine the number of
distinct smooth structures which the $PL$ manifold $M$ admits.  It amounts to the same thing,
but gives us a little more freedom, to determine the smoothings of $M$ where we define a smoothing
of $M$ to be a $PL$-homeomorphism $f:M_{\nu} \ra M$ where $M_{\nu}$ is a smooth manifold.  Such a
smoothing defines the smooth structure $f^{-1*}M_{\nu}$ on $M$.  Two smoothings $f_0:M_{\nu_0} \ra
M$ and $f_1:M_{\nu_1} \ra M$ are called concordant if there is a smooth structure on $M \times I$ which
restricts to be $f_0^{-1*}M_{\nu_0}$ on $M \times {0}$ and $f_1^{-1*}M_{\nu_1}$ on $M \times {1}$. 
Another theorem of \cite{HM} states that the concordance classes of smoothings of $M$ are in
bijective correspondence with $[M,PL/O]$ where $PL/O$ denotes the fiber of $\pi_{O,PL}$.  We shall
call two smoothings $f_i: M_{\nu_i} \ra M$ diffeomorphic if there is a diffeomorphism $g:M_{\nu_0}
\ra M_{\nu_1}$ such that $f_0 = f_1 \circ g$.  Of course one
often would rather know the classification of the smoothings of $M$ up to diffeomorphism and these
are the quotient of the concordance classes of smoothings of $M$ under the action of the
$PL$-homeomorphisms of $M$.  A $PL$-homeomorphism $h: M \ra M$ acts on a smoothing $M_{\nu}: \ra M$
by composition, $h \cdot f = h \circ f$.  If $h$ is $PL$-isotopic to the identity then it acts
trivially on any smoothing.  In the case where $M \cong (S^n)^{PL}$ it is known that all
$PL$-homeomorphisms which preserve orientation are $PL$ isotopic to the identity and hence
$[S^n,PL/O] = \pi_n(PL/O)$ corresponds bijectively with the set of orientation preserving diffeomorphism
classes of smoothings of $S^n$, otherwise known as $\Theta_n$ (for $n \neq 3,4$).  For $i\leq 6$ it
is known that $\pi_i(PL/O) = 0$.  In general the computation of $[M,PL/O]$ for an arbitrary
manifold is a (not necessarily easy) task in homotopy theory.   

The problem of classifying the $PL$ structures on a given topological manifold was solved by Kirby
and Siebenmann \cite{KiSi} who showed that the problem again reduced to the lifting question for the
classifying maps of the stable tangent bundles.  Specifically every compact topological manifold $M$ 
has a stable topological normal bundle which classified by a map $\nu_M^{Top}$ from $M$ to the
classifying space $BTop$.  There is a universal forgetful map $\pi_{PL,Top}: BPL \ra BTop$ and a
topological manifold $M$ admits  a $PL$-structure if and only if the following lifting problem can be
solved. 
\[
\begin{diagram}
\node[2]{BPL} \arrow{s,r}{\pi_{PL,Top}} \\
\node{M} \arrow{e,b}{\nu_M^{Top}} \arrow{ne,..}
 \node{BTop}
\end{diagram}
\]
Moreover, the homotopy fiber of $\pi_{PL,Top}$, $Top/PL$, has but one non-trivial homotopy group and
this is $\Z/2$ in dimension $3$.  There is thus a single obstruction to placing a $PL$ structure on a
topological manifold $M$ which lies in $H^4(M;\Z/2)$ and it is called the Kirby-Siebenmann invariant. 
If $M$ is a topological manifold which admits a $PL$ structure then the concordance classes of $PL$
structures on $M$ are in bijective correspondence with $[M,Top/PL] \cong H^3(M;\Z/2)$.

Having reviewed the fundamental results concerning the smoothing of manifolds we now consider the $PL$
and topological classification of highly connected manifolds of dimension $7$.  We shall write $P^{PL}$
and $PL^{Top}$ for the $PL$ and topological manifolds underlying the smooth manifold $P$.  

\begin{Theorem} \label{thm7homeo}
Let $P_0$ and $P_1$ be highly connected, closed, smooth $7$-manifolds with
quadratic linking families $Q(P_0)$ and $Q(P_1)$.
The following are equivalent:
\begin{enumerate}
\item{$P_0$ is homeomorphic to $P_1$,}
\item{$P_0$ is $PL$-homeomorphic to $P_1$,}
\item{$P_0$ is almost diffeomorphic to $P_1$,}
\item{$Q(P_0)$ is isometric to $Q(P_1)$.}
\end{enumerate}
\end{Theorem}

\begin{proof}
Evidently {\em 3} implies {\em 2} and {\em 2} implies {\em 1}.  The main result of Chapter
\ref{chaptop} was the equivalence of {\em 3} and {\em 4}.  In dimension $7$ we see that
$[P,PL/O] \cong H^7(P;\pi_7(PL/O)) \cong \Theta_7$.  It follows that all possible smooth structures on
$P_0^{PL}$ are obtained by taking the connected sum of $P_0$ with a homotopy sphere and thus
{\em 2} implies {\em 3}.  Finally we show that {\em 1} implies {\em 2}.

Let $f:P_0 \ra P_1$ be a homeomorphism.  Suppose, to begin with, that $P_0$ and $P_1$ have torsion free
third homology.  Then $Q(P_i) \cong (\Z^l,0,\beta_i)$ where $l$ is the rank of $H^4(P_i)$, 
$0$ indicates the trivial quadratic linking function and $\beta_i$ is the
tangential invariant of $P_i$.  By Novikov's theorem on the topological invariance of the rational
Pontryagin classes, $f_!(\beta_0) = \beta_1$ and hence $f$ induces an isometry of the quadratic linking
families.  Hence by Theorem \ref{thmb}, $P_0$ and $P_1$ are almost diffeomorphic and in particular
$PL$-homeomorphic. If $P_0$ and $P_1$ have finite fourth cohomology then by Corollary
\ref{corwilkdata} the quadratic linking functions $Q(P_i) = (TG,q_i,\beta_i)$ are determined by
their Wilkens data $(b(q_i),\beta_i)$ and Kervaire-Arf invariants $K(q_i)$.  Now the
homeomorphism $f$ preserves the linking form since $b$ and it also preserves the tangential invariant
since by \cite{KS} $\beta(P_i) = \halfp(P_i)$ is a topological invariant.  But by Lemma \ref{props1}
$K(q_i) = -\bar{s}_1(q_i) = -\bar{s}_1(P_i)$ and again by \cite{KS} the invariant $\bar{s}_1$ is a
topological invariant.  Thus $P_0$ and $P_1$ have isometric quadratic linking functions and are
therefore almost diffeomorphic and in particular $PL$-homeomorphic.

Thus if $P$ is a highly connected $7$-manifold with either finite homology or torsion free homology, then
every element of $H^3(P;\pi_3(PL/O))$ is realized by a self-smoothing $f:P^{PL} \ra P^{Top}$.  Finally
suppose that $H^4(P_0) \cong F \oplus TG$ where $F$ is a free group and $TG$ is a finite group.
We choose a splitting of $P_0$ as the connected sum of manifolds $P_0 = P_0(F) \# P_0(\text(TG)$
where $F = H^4(P(F))$ and $TG=H^4(P(TG))$.  Given any $x \in
H^3(P_0;\pi_3(PL/O))$ we may write $x = y \oplus z$ where $y \in
H^3(P_0(F);\pi_3(PL/O))$ and $z \in H^3(P_0(TG)^{Top};\pi_3(PL/O))$.  Then
by the previous arguments $y$ and $z$ are realized by self $PL$ smoothings
$f_y:P_0(TG)^{PL} \ra P_0(TG)^{Top}$ and $f_z:P_0(F)^{PL} \ra
P_0(F)^{Top}$.  If we now form the connected sum $f_y \# f_z:P_0^{PL} \ra P_0^{Top}$ then it realizes
$x$ and is manifestly a self smoothing.  Thus every element of $H^3(P_0;\pi_3(PL/O))$ is realized
by a self $PL$ smoothing and hence every $PL$ manifold homeomorphic to $P_0$ is
$PL$-homeomorphic to $P_0$ and in particular $P^{PL}_1$ is $PL$-homeomorphic to $P_0$.
\end{proof}

\section{The surgery exact sequence} \label{secsurg}
We present now a rapid review of the surgery exact sequence which is an important tool for the
classification of manifolds. We begin with a pair of finite CW-complexes
$(Y,X)$ which we assume to be simply connected and where we allow $X$ to be empty.  Let $\Cat$ denote
either of the smooth ($O$), piecewise linear ($PL$) or topological ($\Top$) categories.  We wish to know
if  there is a homotopy equivalence $f:(W,M) \ra (Y,X)$ from a compact pair of $\Cat$-manifolds
$(W,M)$ to $(Y,X)$. A necessary condition for such a homotopy equivalence is of course that the
pair $(Y,X)$ satisfy Poincar\'{e} duality.  That is, for some $n$ there is a class $[Y,X] \in
H_n(Y,X;\Z)$ such that $\cap[Y,X]:H^k(Y,X) \cong H_{n-k}(Y)$.  As the pair $(Y,X)$ is a finite complex
there a large $N$ and an embedding $(Y,X) \hra D^N,S^{N-1})$ where $(D^N,S^{N-1})$ is $N$-disc with
boundary the $(N-1)$-sphere $S^{N-1}$.  Inside $D$, the complex $Y$ has a regular
neighborhood $N$.  The manifold $N$ is a smooth manifold and retracts onto
$Y$ via a retraction we denote by $r:N \ra Y$.  It is a consequence of Poincar\'{e} duality for the pair
$(Y,X)$ that the map $r|_{\del N}$, when considered as a fibration, has homotopy fiber the homotopy type
of a sphere.  The stable fiber homotopy type of the fibration defined by
$r|_{\del N}$ is independent of the embedding and is a homotopy invariant of the pair
$(Y,X)$ (\cite{Sp}).  It is called the Spivak normal bundle of $(Y,X)$ and is classified by a map $\nu: Y
\ra BG$ where $BG$ is the classifying space for spherical fibrations.  Now if $Y$ were a $\Cat$ manifold
then $r$ would be a stable normal bundle of $Y$ and $r|_{\del N}$ would be the sphere bundle of the
stable normal bundle.  Hence a further necessary condition for $(Y,X)$ to have the homotopy type of a
compact manifold with boundary is that its Spivak normal bundle have a $\Cat$ bundle reduction. 
Assuming that the dimension of $(Y,X)$ as a Poincar\'{e} space is greater than $5$ then it is a theorem
of Browder  \cite{Bro} and Novikov \cite{No} (in the smooth and $PL$ categories) and Kirby and
Siebenmann \cite{KiSi} (in the $Top$ category)\footnote{By their extension of transversality methods to
$Top$ \cite[III]{KiSi})} that a $Cat$ reduction of the Spivak normal bundle of a simply connected
Poincar\'{e} space is in fact a sufficient condition for it to have the homotopy type of a $Cat$
manifold.

With an answer to the question of when a homotopy type $(Y,X)$ contains the homotopy type of a
manifold, surgery proceeds to the task of classifying $Cat$-manifolds homotopy equivalent to
$(Y,X)$ up to $Cat$-isomorphism.  However, the primary objects which surgery theory seeks to
classify are in fact smoothings of a finite $CW$-pair $(Y,X)$.  The  $\Cat$-structure set of $Y$,
${\mathcal S}^{\Cat}(Y)$, consists of equivalence classes of $\Cat$-smoothings in $Y$ where a
$\Cat$ smoothing in $Y$ is an orientation preserving homotopy equivalences of pairs $f:(M,
\del M)  \rightarrow (Y, X)$ from a compact $\Cat$-manifold with boundary to
$(Y,X)$.  Two homotopy equivalences
$(M_{1}, f_{1})$ and $(M_{2}, f_{2})$ are equivalent in
${\mathcal S}^{\Cat}(Y)$ if there is a $\Cat$-isomorphism $g:
M_1 \rightarrow M_2$ such that $f_1$ is homotopy equivalent to $f_2 \circ
g$.  In order to compute the set of $\Cat$ smoothings in $Y$ one maps ${\mathcal
S}^{\Cat}(Y)$ to the set  ${\mathcal{N}}^{Cat}(Y)$ of degree one normal maps in $Y$ modulo normal
bordism.  We now define this set.  A degree one normal map to $Y$ is a
bundle map $(f, \hat f)$ 
$$
\begin{array}{ccc}
\nu_{M}  & \stackrel{\hat{f}}{\ra} & \nu \\
\downarrow & & \downarrow \\
(M,\del M)  & \stackrel{f}{\ra} & (Y,X)          
\end{array}
$$
from the stable normal bundle $\nu_{M}$ of a Cat-manifold $M$ to a Cat-bundle over
$Y$ such that map induced on homology by $f$ carries the fundamental class
of $(M,\del M)$ to $[Y,X]$.  The equivalence between degree one normal maps is
simply that of bordism with boundary over the bundle $\nu \rightarrow Y$. 
Specifically, two degree one normal maps
$$
\begin{array}{ccc}
\nu_{M_{0}} &  \stackrel{\hat{f_0}}{\ra} & \nu_{0} \\  
\downarrow & & \downarrow \\
(M_{0},\partial M_0) & \stackrel{f_0}{\ra} & (Y,X)   
\end{array} \hskip 2cm
\begin{array}{ccc}
\nu_{M_1} & \stackrel{\hat{f_1}}{\ra} & \nu_{1} \\  
\downarrow & & \downarrow \\
(M_{1},\partial M_1) & \stackrel{f_1}{\ra} & (Y,X)   
\end{array}
$$
are normally bordant if there exists a compact manifold $W$ and a degree
one normal map $(F, \hat{F})$
$$
\begin{array}{ccc}
\nu_{W}  & \stackrel{\hat{F}}{\ra} & \nu \times I \\
\downarrow & & \downarrow        \\
W & \stackrel{F}{\ra} & Y \times I    
\end{array}
$$
with the following two properties: {\em (i}) $\partial
W=M_{1}~\cup~M_{2}~\cup~N$ where $N$ is a compact manifold with 
$\partial N = N \cap (M_1 \cup M_2) = \partial (M_1 \cup M_2)$,
and {\em (ii)} $(F|_{M_i}, \hat F|_{M_i}) = (f_i, \hat f_i)$.

To report the essential results of simply connected surgery theory we
first define the following abelian groups.
$$ L_{n}(\Z) = \left\{ \begin{array}{ll}
	   {\mathbf Z}, & \mbox{ if } n = 4\,k;\\
	   {\mathbf Z}_{2} , & \mbox{ if } n = 4\,k + 2;\\
	    0 , & \mbox{otherwise.} 
	\end{array} 
	\right\}, \hskip 2cm L_{n}(\Z \rightarrow \Z) = 0. $$

\noindent
The following two facts are well-known, see \cite{Wa5} or \cite{MM}.
\begin{Fact} \label{factses}
(1)  Let $M$ be a simply connected manifold and let $G/\Cat$ denote the fiber of the map $B\Cat
\rightarrow BG$, then ${\mathcal N}^{Cat}(M) \cong [M, G/\Cat]$.\\ 
(2)  If $k>5$ there is a pair of surgery exact sequences of sets
$$
\begin{array}{ccccccc}
L_{k+1}(\Z \rightarrow \Z) &
\longrightarrow &
{\mathcal S}^{\Cat}(W^{k}) &
\stackrel{{\mathcal{\eta}^{\Cat}}}{\longrightarrow} &
{\mathcal N}^{\Cat}(W^{k}) &
\longrightarrow & 
L_{k}(\Z \rightarrow \Z)\\
\downarrow & & i^* \downarrow & & i^* \downarrow & & \downarrow \\ 
L_k(\Z) & \longrightarrow &
{\mathcal S}^{\Cat}(\partial W^{k-1}) &
\stackrel{{\mathcal{\eta}^{\Cat}}}{\longrightarrow} &
{\mathcal N}^{\Cat}(\partial W^{k-1}) & \longrightarrow & L_{k-1}(\Z)
\end{array}
$$
where exactness means that each set has a preferred element $*$ and that at every stage the
pre-image of $*$ is the image of the previous map.  The base point of both $\SCat$ and $\NCat$ is
the the class of the identity map $Id: W \ra W$ or $Id: \del W \ra \del W$.  The base point of the
$L$-groups is the identity element.  In fact there is an action of the left-most $L$-group on
$\SCat$ and two smoothings in $\SCat$ are in the same orbit if and only if they have
the same image under $\eta^{Cat}$.  In the smooth category this action
is connect sum in the domain of a smoothing with a homotopy sphere bounding a parallelizable
manifold and in either the $PL$ or $Top$ categories the action is trivial and $\eta$ is
injective.  Finally, the vertical maps $i^*$ are given by restricting a smoothing or degree one
normal map to its boundary.
\end{Fact}

We now illustrate this theory in the $PL$ case for the manifolds
$(W,\del W)=(L,P)$.  In this case $n=8$,
$L_{7}(e) = 0$ and $\eta^{PL}$ is a bijection.  Thus both $PL$-structure
sets and $PL$-normal invariant sets coincide.  Moreover the normal
invariants are extremely easy to calculate since $[L,G/PL]$ and
$[P,G/PL]$ can be identified, via the the primary obstruction to null
homotopy, with $H^4(L;\pi_4(G/PL)) \cong H^4(L;\Z)$ and
$H^4(P;\pi_4(G/PL)) \cong H^4(P;\Z)$ respectively.   The surgery exact
sequences become the following commutative squares.\\
$$
\left\{
\begin{array}{ccc}
\SPL (L) & \stackrel{\eta}{\cong} & \NPL(L) \\
i^* \downarrow  &   & i^* \downarrow \\
\SPL (P) & \stackrel{\eta}{\cong} & \NPL(P) \\
\end{array} \right\}
\cong \left\{
\begin{array}{ccc}
H^{4}(L) & \cong & H^{4j}(L)\\
i^* \downarrow &  & i^*\downarrow\\
H^{4}(P) & \cong & H^{4j}(P)\\
\end{array} \right\}
\noindent
$$
We emphasize that the vertical maps $i^*$ are given by simply restricting a
$PL$-structure or $PL$-normal invariant to the structure or normal
invariant on the boundary.  Since the map $i^*:H^{4}(L) \ra H^{4}(P)$ is always onto it 
follows that every degree one normal map in $P$ is normally bordant to the restriction of a
degree one normal map in $L$ and thus that every equivalence class of $PL$-smoothings in $P$ has a
representative which is the restriction of a $PL$-smoothing in $L$.  

\section{Constructing smoothings in $L$ and $P$} \label{secsmoothings}
Let $SG(n)$ denote the $H$-space of orientation preserving self homotopy equivalences of $S^{n-1}$.
In what follows the groups $SG(4)$ and $SO(4)$, and in particular $\pi_3(SG(4))$ and
$\pi_3(SO(4))$, shall play a crucial role so we begin by digressing somewhat to describe these groups
and the relationships between them.

We start with $SO(4)$. Let $|~|$ denote the standard metric on $\R^4$ and identify $\R^4$ with the
quaternions, $\Ha$.  Then $S^3 = \{ x \in \Ha : |x| = 1 \}$ is the set of unit quaternions.  For $y \in
\Ha$ and $x \in S^3$ we define the function
$$\begin{array}{cccc}
a(m,n):& S^3 & \ra & SO(4) \\
& x & \mapsto & [y \mapsto x^{n+m}\cdot y \cdot x^{-m}]
\end{array}$$
where $\cdot$ indicates quaternionic multiplication.  The mapping $a(m,n)$ is manifestly linear and
orientation preserving and since $|x^{n+m}\cdot y \cdot x^{-m}| = |x|^{n+m}|y||x|^{-m|} = |y|$ and it
follows that $a(m,n)(x)$ is an element of $SO(4)$ for each $x \in S^3$.  We define $\ba (m,n)$ to be
the homotopy class of $a(m,n)$.  It is well known that the  mapping
$$\begin{array}{ccc}
\Z \oplus \Z & \ra & \pi_3(SO(4)) \\
(m,n) & \mapsto & \ba(m,n)
\end{array}$$
defines an isomorphism of groups.  
Given the function $a(m,n)$ we define the associated twist $A(m,n)$
as follows
$$\begin{array}{cccc}
A(m,n):& S^3 \times D^4 & \ra & S^3 \times D^4 \\
& (x,y) & \mapsto & (x,a(m,n)(y)).
\end{array}$$
Identifying $S^3 \times D^4$ with $\del D^4 \times D^4$ we may use $A(m,n)$ to glue two copies of $D^4
\times D^4$ together thereby form a $D^4$-bundle over $S^4$, $\xi(m,n):L(m,n) \ra S^4$.  We define the
Euler number of $\ba(m,n)$ to be the  Euler number of $\xi(m,n)$,
$E(\ba(m,n)) := E(\xi(m,n))$.  The following
facts are also well known and may be found in \cite{Hu} ($r=\ba(1,0)$ and $s=\ba(0,1)$).
\begin{enumerate}
\item{The Euler number of $\ba(m,n)$, $E(\ba(m,n))$, is $n$.}
\item{The stabilization of $\ba(m,n)$, $S(\ba(m,n))$, is $2m+n \in \Z = \pi_3(SO)$.}
\item{If $i:SO(3) \hra SO(4)$ denotes the standard inclusion with induced map $i_*:\pi_3(SO(3)) \hra
\pi_3(SO(4))$, and if $m \in Z = \pi_3(SO(3))$, then $i_*(m) =  \bar \alpha(m,0)$.}
\end{enumerate}

By restricting an element of $SO(n)$ to $S^{n-1}$ we obtain an orientation preserving homotopy
equivalence of $S^{n-1}$ and this defines an inclusion $i_n:SO(n) \hra SG(n)$.

\begin{Lemma} \label{lemmapi_3(SG(4))}
There is an isomorphism $\pi_3(SG(4)) \cong \Z/12 \oplus \Z$ such that $$i_{4*}(\ba(m,n)) = (m \,
(\modu~12),n).$$
\end{Lemma}

\begin{proof}
Let $SF(n)$ denote the sub-$H$-space of $SG(n+1)$ consisting of orientation
preserving homotopy equivalences of $S^n$ which fix a point $p$.   By the adjoint correspondence there
is an isomorphism $A_{n,j}:\pi_j(SF(n)) \cong \pi_{j+n}(S^n)$.  We relate the homotopy groups of $SF(n)$
to $SG(n+1)$ by noting that $SF(n)$ is the fiber of the fibration
$$\begin{array}{cccc}
ev: & SG(n+F) & \ra & S^n \\
& f & \mapsto & f(p).
\end{array}$$
Moreover, $ev|_{SO(n+1)}: SO(n+1) \ra S^n$ is the standard
fibration with fiber $SO(n) \subset SF(n)$.  When $n=3$, $ev|_{SO(4)}$ has a section $s:S^3 \ra SO(4)$
given by $s(x) = a(0,n)(x)$.  Applying the long exact homotopy sequence to this pair of fibrations at
$\pi_3$ we obtain the following commutative diagram with exact rows.
\[
\begin{diagram}
\divide\dgARROWLENGTH	by2
\node{0} \arrow{e} \node{\Z} \arrow{s,r}{\cong} \arrow{e} \node{\Z \oplus \Z}
\arrow{s,r}{\cong} \arrow{e} \node{\Z} \arrow{s,r}{\cong} \arrow{e} \node{0} \\
\node{0} \arrow{e} \node{\pi_3(SO(3))} \arrow{s,t}{j_*} \arrow{e} \node{\pi_3(SO(4))}
\arrow{s,r}{i_{4*}} \arrow{e,t}{ev|_{SO(4)*}} \node{\pi_3(S^3)} \arrow{s,r}{Id} \arrow{e} \node{0} \\
\node{0} \arrow{e} \node{\pi_3(SF(3))} \arrow{s,r}{\cong} \arrow{e} \node{\pi_3(SG(4))}
\arrow{s,r}{\cong} \arrow{e,t}{ev_*} \node{\pi_3(S^3)} \arrow{s,r}{\cong} \arrow{e} \node{0} \\
\node{0} \arrow{e} \node{\Z/12}  \arrow{e} \node{\Z/12 \oplus \Z} \arrow{e} \node{\Z}  \arrow{e}
\node{0} 
\end{diagram}
\]  
To conclude the proof we note the following.  Firstly, if $n \in \pi_3(S^3)$ then $s_*(n) =
\ba(0,n) \in \pi_3(SO(4))$.  Secondly, if $m \in \pi_3(SO(3))$ then $i_*(m) = \ba(m,0)$.
Finally, if $j_n:SO(n) \hra SF(n)$ denotes the inclusion of one fiber
into the other then the homomorphism
$$A_{n,j} \circ j_{n*}: \pi_j(SO(n)) \ra \pi_{n+j}(S^n)$$
is the usual $J$-homomorphism and the $J$-homomorphism $J_{3,3}:\pi_3(SO(3)) \ra
\pi_6(S^3)$ is known to be the surjection $\Z \ra \Z/12$ (\cite{JW2}).
\end{proof}

Now, let $L$ be a $8$-dimensional handlebody and let the quadratic function of $L$ be $\kappa(L) =
\kappa(H,\lambda,\alpha)$.  Recall from Section \ref{chapprelim}.\ref{secwall} that for every choice of
basis, $\{ v_1 \dots v_l \}$, for $H$ there is a corresponding presentation for the handlebody $L$,
$$\phi=\sqcup_{i=1}^l \phi_i : \sqcup_{i=1}^l (S^{3} \times D^4)_i \hra \del D^{8},$$
and that $L$ is diffeomorphic to the adjunction space $D^8 \cup_{\phi} \sqcup_{i=1}^l (D^4 \times
D^4)_i$.  The isotopy invariants of $\phi$ are the linking numbers $c_{ij}$ of the
attaching spheres $\bar{\phi}(S^{3}_i)$ where $S^3_i = (S^3 \times 0)_i$, and the framing
invariants $\ba(i) \in \pi_3(SO(4))$ which we denote now by $ba_i$.  They are related to $\lambda$ and
$\alpha$ as follows:
$$c_{ij} = \lambda(v_i,v_j)~(i \neq j), \hskip 1cm E(\ba_i) = \lambda(v_i,v_i), \hskip 1cm S(\ba_i) =
\alpha(v_i).$$

\begin{Lemma}
Suppose that the handlebody $L$ has a presentation 
$$\phi=\sqcup_{i=1}^l \phi_i : \sqcup_{i=1}^n (S^{3} \times D^4)_i \hra \del D^{8}$$  
with invariants $c_{ij}, (i\neq j)$, and $\ba_i$, and let $\bg_i = \ba(m_i,n_i) \in \pi_3(SO(4))$.  Then
$$\phi(\bg_1, \dots, \bg_l) : = \sqcup_{i=1}^l \phi_i \circ A(m_i,n_i)$$
is a presentation with invariants $c_{ij}, (i\neq j)$, and $\ba_i + \bg_i$. 
\end{Lemma}

\begin{proof}
As $A(m_i,n_i)(x,0) = (x,0)$ for every $x \in S^3$, $\phi_i \circ A(m_i,n_i)|_{S^3_i} =
\phi_i|_{S^3_i}$.  Hence the attaching spheres and linking numbers $c_{ij}$ are unchanged.  Recall
that the framing invariant $\ba_i$ is defined as the homotopy class of $a(\phi_i)$ where 
$$\psi_i:\phi_i((S^3 \times S^4)_i) \ra \bar{\phi}_i(S^3_i) \times D^4$$
is a framing of $\phi_i((S^3 \times D^4)_i)$ which extends to a framing of the normal bundle of a
$4$-disc embedded in $D^8$ bounding $\phi_i((S^3 \times 0)_i)$ and where
$$\begin{array}{cccc}
\psi_i \circ \phi_i:& (S^{3} \times D^4)_i & \ra & \bar{\phi}_i(S^{3}_i) \times D^4\\
& (x,y) & \mapsto & (\bar{\phi}_i(x),a(\phi_i)(x)(y)).
\end{array}$$ 
After isotopy we may assume that $a(\phi_i) = a(p_i,q_i)$ for some $p_i$ and $q_i$.  For the 
presentation $\phi(\bg_1, \dots, \bg_l)$ we consider instead the diffeomorphism 
$\psi_i \circ \phi_i \circ A(m_i,n_i)$ which is
$$\begin{array} {cccccc}
A(p_i,q_i) \circ A(m_i,n_i) : & (S^{3} \times D^4)_i & \ra & \bar{\phi}_i(S^{3}) \times D^4 \\  
& (x,y) & \mapsto & (\bar{\phi}_i(x),x^{m_i+n_i+p_i+q_i}.y.x^{-m_i-p_i}).\\
\end{array}$$
It follows that the framing functions for $a(\phi(\bg_1, \dots, \bg_l)_i)$ are $a(n_i+p_i,m_i+q_i)$ and
that the framing invariants for $\phi(\bg_1, \dots, \bg_l)$ are
$$\ba(n_i+p_i,m_i+q_i) = \ba(p_i,q_i) + \ba(n_i,m_i) = \ba_i + \bg_i.$$
\end{proof}

\noindent
We shall abbreviate $\bg_1,\dots, \bg_l$ to $\bg$ and denote the handlebody obtained from the
twisted presentation by  $L(\bg)$.  

\begin{Lemma} \label{lemmahes}
Let $L$, $\phi$, $L(\bg)$ and $\phi(\bg)$ be as in the previous lemma
and let  $\bg_i \in \Ker[i_{4*}:\pi_3(SO(4)) \ra \pi_3(SG(4))]$ for $i=1,\dots,l$.  Then there is a
homotopy equivalence
$$g(\bg): (L(\bg), \del L(\bg)) \ra (L, \del L).$$
\end{Lemma} 

\begin{proof}
We make the identifications 
$$L = D^8 \cup_{\phi} \sqcup_{i=1}^l (D^4 \times D^4)_i \text{~~~and~~~}
L(\bg) = D^8 \cup_{\phi(\bg)} \sqcup_{i=1}^l (D^4 \times D^4)_i.$$
The map $g(\bg)$ shall take each handle of $L$ to the corresponding handle of 
$L(\bg)$.  On the zero handle, $D^8$, we define $g(\bg)$ to be the
identity.  Now consider the $i$th handles of $L$ and $L(\bg)$ which we denote by $(D^4
\times D^4)_i$.  The map $g(\bg)$ has been defined only on $(\del D^4 \times D^4)_i$
where it is $\phi_i^{-1} \circ \phi(\bg)_i $.  Now $\phi(\bg)_i = \phi_i
\circ A(m_i,n_i)$ where $A(m_i,n_i)$ is the twisting corresponding to $a(m_i,n_i)$.  It follows that
$g(\bg)|_{(\del D^4 \times D^4)_i} = A(m_i,n_i)$.  As $a(m_i,n_i)$ represents $\bg_i$ and
$i_{4*}(\bg_i) = 0 \in \pi_3(SG(4))$ we conclude from Lemma
\ref{lemmapi_3(SG(4))} firstly that each $n_i = 0$ and that $12$ divides $m_i$.  Secondly,  for every
$j$, $i_4 \circ a(12j,0): S^3 \ra SG(4)$ is null homotopic so we let $b(12j,0): D^4 \ra SG(4)$ be such 
a null homotopy and we define an extension of
$A(12j,0)$ as follows.  For $x \in D^4$, $y \in S^3$, $t \in [0,1]$ and $ty \in D^4$ let 
$$\begin{array}{cccc}
B(12j,0): & D^4 \times D^4 & \ra & D^4 \times D^4 \\
& (x,ty) & \mapsto & (x,tb(12j,0)(x)(y)).
\end{array}$$
On each handle $(D^4 \times D^4)_i$ we define $g(\bg)$ to be $B(m_i,0)$.

We now show that $g(\bg)$ is a homotopy equivalence of pairs.  The boundary of $L$ is
obtained by performing surgeries on $\del D^8$ described by the embeddings $\phi$,
$$\del L = (\del D^8 -[\sqcup_{i=n}^l \phi_i((\del D^4 \times D^4)_i)] 
\cup_{\phi|_{\sqcup _{i=1}^l \del D^4 \times \del D^4}} \sqcup_{i=1}^l (D^4 \times \del D^4)_i.$$
A similar statement holds true for $\del L(\bg)$.  Since   
$g(\bg)|_{\del D^8} = Id_{\del D^8}$ and for each $i$, $B(m_i,0)((D^4 \times \del D^4)_i) = (D^4
\times \del D^4)_i$, we see that $g(\bg)$ maps the boundary of $L(\bg)$ to the boundary
of $L$.   The generators of $H_4(L)$ and $H_4(L(\bg))$ may be represented by embedded spheres
constructed by gluing the cores of each handle to $4$-discs embedded in $D^8$ and bounding the
corresponding attaching spheres.  Since $B(m_i,0)$ is the identity map when restricted to the core of
each handle and since
$g(\bg)$ is the identity on the zero handle we see that $g(\bg)$ carries the generators of $H_4(L)$ to
the generators of $H_4(L(\bg))$.   Since $L$ and $L(\bg)$ have non-zero homology only in dimensions
$0$ and $4$ we conclude that $g(\bg)_*:H_*(L(\bg)) \cong H_*(L)$.  Since all the groups involved are
free we conclude that $g(\bg)^*:H^*(L) \cong H^*(L(\bg))$.  We now apply the naturality of the
Poincar\'{e} duality isomorphisms $H^{8-*}(L) \cong H_*(L,\del L)$ and $H^{8-*}(L(\bg)) \cong
H_*(L(\bg),\del L(\bg))$ to conclude that $g(\bg)_*:H_*(L(\bg),\del L(\bg)) \cong H_*(L,\del L)$. 
Applying the five lemma to the long exact homology sequences of $(L,\del L)$ and $(L(\bg),\del L(\bg))$,
it follows that $(g(\bg)|_{\del L(\bg)})_*:H_*(\del L(\bg)) \cong H_*(\del L)$.  Thus $g(\bg)$ induces a
homology equivalence of simply connected pairs and so by Whitehead's Theorem $g(\bg)$ induces a homotopy
equivalence of pairs.
\end{proof}

\begin{Corollary} \label{corhes}
Let $L$ be a handlebody with quadratic function $\kappa(H,\lambda,\alpha)$ and for each
$\epsilon \in 24 \cdot H^*$ let $L_{\epsilon}$ be the handlebody with quadratic function
$\kappa(L_{\epsilon}) = \kappa(H,\lambda,\alpha + \epsilon)$.  Then there is a homotopy equivalence
$$g_{\epsilon}:  (L_{\epsilon},\del L_{\epsilon}) \ra (L, \del L)$$ 
such that $g_{\epsilon*} = Id: H \ra H$.
\end{Corollary}

\begin{proof}
Let $\{v_1, \dots, v_l \}$ be a basis for $H$ and let $\phi$ and $\phi_{\epsilon}$ be presentations
for $L$ and $L_{\epsilon}$ respectively which correspond to this basis.  If the invariants of $\phi$
are $c_{ij}$ and $\ba_i$ then the invariants for $\phi_{\epsilon}$ are $c_{ij}$ and $\ba_i + \bg_i$
where $\bg_i = \ba(\epsilon(v_i)/2,0)$.  As $24$ divides $\epsilon$, $12$ divides $\epsilon(v_i)/2$
for each $i$ and so by Lemma \ref{lemmapi_3(SG(4))}, $\bg_i \in \Ker(i_{4*})$.  Thus we apply Lemma
\ref{lemmahes} and set $f_\epsilon = g(\bg)$.  The induced map $g_{\epsilon_*}$ is the
identity because, as we noted in the proof of Lemma \ref{lemmahes}, $g_{\epsilon}$ maps the embedded
sphere representing $v_i$ in $L_{\epsilon}$ to the corresponding sphere in $L$.  Hence for each $i$,
$g_{\epsilon*}(v_i) = v_i$.  Since $\{ v_1, \dots, v_n \}$ is a basis for $H$, $g_{\epsilon*} = Id_H$.
\end{proof}

\begin{Lemma} \label{lemmanormalinvts}
Let $g_{\epsilon} : L(\kappa(H,\lambda,\alpha + \epsilon)) \ra L(\kappa(H,\lambda,\alpha))$ be the 
homotopy equivalence of handle bodies defined in Corollary \ref{corhes}.  Then
$g_{\epsilon}$ has $PL$-normal invariant $\eta^{PL}(g_{\epsilon}) = \epsilon/24 \in \NPL(L) \cong
H^{4j}(L) = H^*.$
\end{Lemma}

\begin{proof}  Having previously identified $\NPL(L)$ with $H^4(L;\pi_4(G/PL))$
we consider the change of coefficient map to $U:H^4(L;\pi_4(G/PL)) \ra H^4(L;\pi_4(BPL))$.  If
$[f] \in H^4(L;\pi_4(G/PL))$ corresponds to a degree one normal map 
$$
\begin{array}{ccc}
\nu_{M}  & \stackrel{\hat{f}}{\ra} & \nu \\
\downarrow & & \downarrow \\
(M,\del M)  & \stackrel{f}{\ra} & (L,P),          
\end{array}
$$
then $U([f])$ corresponds to the difference of bundles $\nu_L - \nu \in [L,BPL]$.  Now
the map $\pi_4(G/PL) \ra \pi_4(BPL)$ is $\cdot 24: \Z \ra \Z$ as it has cokernel $\pi_3^S \cong
\Z/24$ and $\pi_4(BPL) \cong \pi_4(BO) \cong \Z$.  Thus $U$ is injective and all normal invariants
in $L$ are detected by the bundle data.  We wish to calculate $U(\eta(g_{\epsilon}))$.  Now, as
$g_{\epsilon}$ is a homotopy equivalence $\eta( g_{\epsilon})$ is the degree on
normal map
$$
\begin{array}{ccc}
\nu_{L_{\epsilon}}  & \stackrel{\hat{g_{\epsilon}}}{\ra} & (g_{\epsilon}^{-1})^{*}\nu_{L_{\epsilon}} \\
\downarrow & & \downarrow \\
(L_{\epsilon},\del L_{\epsilon})  & \stackrel{g_{\epsilon}}{\ra} & (L,\del L)        
\end{array}
$$
where $L=L(\kappa(\lambda,\alpha))$ and $L_{\epsilon} = L(\kappa(\lambda,\alpha + \epsilon))$.
Note that $\nu_L = -\alpha$ is the inverse of the stable tangent bundle of $L$ and that
$g_{\epsilon}$ induces the identity on cohomology.
$$\begin{array}{cl}
U(\eta(g_{\epsilon})) & = \nu_{L} - (g_{\epsilon}^{-1})^*\nu_{L_{\epsilon}} \\  
& =  (-\alpha) - (g_{\epsilon}^{-1})^*(-\alpha - \epsilon) \\
& = \epsilon \\
\end{array}$$
Since $U$ is the coefficient homomorphism which is multiplication by $24$,
$\eta(g_{\epsilon}) = \epsilon/24$.
\end{proof}

\begin{Corollary} \label{corsslp}
Let $L$ be an $8$-dimensional handlebody with quadratic function
$\kappa = \kappa(H,\lambda,\alpha)$ and boundary $P = \del L$. If $g_{\epsilon}: L_{\epsilon} : =
L_(\lambda, \alpha + \epsilon) \ra L(\lambda,\alpha)$ is the homotopy equivalence of Corollary
\ref{corhes} and let $f_{\epsilon} = \del L_{\epsilon}$, then, 
\begin{enumerate}
\item{ $\SPL(L) = \{ [g_{\epsilon}: L_{ \epsilon} \ra L] |\, \epsilon \in  24\cdot H^* \}$, }

\item{ $\SPL(P) = \{ [f_{\epsilon} : \del L_{\epsilon} \ra \del L] |\, \epsilon \in 24 \cdot H^* \}$. }
\end{enumerate}
\end{Corollary}

\begin{proof}
The first assertions follows immediately from Lemma \ref{lemmanormalinvts} and the isomorphism 
$\SPL(L) \cong H^4(L)$ which was discussion in Section \ref{secsurg} following Fact \ref{factses}.  The
second assertion follows from the fact that the map  $i^*: \SPL(L) \ra \SPL(P)$ is onto, which was
demonstrated in the same discussion.
\end{proof}

It is a simple task to use Corollary \ref{corsslp} to descend from the $PL$ classification of
handlebodies $L$ to their homotopy classification.  We introduce some extra notation so that the theorem
we state is of the same form as Wall's corresponding result for closed manifolds (\cite{Wa1}
Lemma 8).  We let $J:\pi_3(SO) \ra \pi_3^S$ denote the stable $J$-homomorphism which is known to be the
surjection
$\Z \ra \Z/24$.  The $J$-homomorphism induces a coefficient homomorphism for any space $X$,
$J_*:H^*(X;\pi_3(SO(4)) \ra H^*(X;\pi_3^s)$.

\begin{Theorem} \label{thmhomotlclass}
Two handlebodies $L(\kappa_0)$ and $L(\kappa_1)$ with quadratic functions 
$\kappa_0= \kappa(H_0,\lambda_0,\alpha_0)$ and $\kappa_1=\kappa(H_1,\lambda_1,\alpha_1)$ are homotopy
equivalent if and only if there is an isometry of their intersection forms, $\Theta:H_0 \ra H_1,
\lambda_1 \circ (\Theta \times \Theta) = \lambda_0$, such that $J_*(\alpha_1) =
J_*(\Theta^*(\alpha_0))$.  Moreover, every such isomorphism
$\Theta$ is induced by a homotopy equivalence $f:L_0 \ra L_1$.
\end{Theorem}

\begin{proof}
From Corollary \ref{corsslp} we deduce that the handlebody $L_1$ is homotopic to $L_0$ if and only if
$L_1$ is $PL$-homeomorphic to $L(\kappa(H_0,\lambda_0,\alpha_0 + \epsilon))$ for some $\epsilon \in 24
\cdot H^*$.  Now by the classification of handlebodies (see Remark \ref{rempltophb}) this is the case if
and only if there is an isometry 
$$\Theta: \kappa(H_1,\lambda_1,\alpha_1) \cong \kappa(H_0,\lambda_0,\alpha_0 + \epsilon))$$
for some $\epsilon \in 24 \cdot H^*$.  But by Proposition \ref{propqfrecov} there is such an isometry if
and only if there is an isomorphism $\Theta:H_1 \ra H_0$ such that 
$\lambda_1 \circ (\Theta \times \Theta) = \lambda_0$ and
$\Theta^*(\alpha_0 + \epsilon) = \alpha_1$ for some $\epsilon \in 24 \cdot H^*$.  But 
$\Ker(J_*) = 24 \cdot H^*$ and thus there is an $\epsilon \in 24 \cdot H^*$ such that 
$\Theta^*(\alpha_0 + \epsilon) = \alpha_1$ if and only if $J_*(\Theta^*(\alpha_0)) = J_*(\alpha_1)$.

Finally, suppose $\Theta: H_1 \ra H_0$ is as in the statement of this theorem.  Then by the remarks
above there is an $\epsilon \in 24 \cdot H^*_0$ such that $\Theta$ defines an isometry $\kappa_1 \cong
\kappa(H_0,\lambda_0,\alpha_0 + \epsilon)$.  This isometry is implemented by a diffeomorphism $g$
which we compose with $g_{\epsilon}$ from Corollary \ref{corhes}.  Since $g_{\epsilon*} = Id$,
$(g_{\epsilon} \circ g)_* = g_* = \Theta$.
\end{proof}

In order to apply Corollary \ref{corsslp} to the homotopy classification of highly connected
$7$-manifolds we first remind the reader of the definition of the characteristic quadratic
linking family induced by a characteristic quadratic function and define the related concept of a
$J$-quadratic linking family.  Suppose that $\kappa(H,\lambda,\alpha)$ is a (characteristic) quadratic
function.  Recall that $G$ and $F$ denote respectively $\Cok(\hl)$ and $\Ker(\hl)$ ant that $\tau_G:G
\ra G^{**}$ and $\rho:H \ra H/F$ denote the canonical maps whose sets of sections are $\Sec(\rho)$ and
$\Sec(\tau_G)$.  Each section $\Psi \in \Sec(\rho)$ induces a splitting $H^* = F^*(\Psi) \oplus
H^*(\Psi)$ and every $x \in H^*$ splits as $x = (p_F(\Psi)^*(x|_F),x(\Psi))$ where $p_F(\Psi): H \ra 
F$.  In addition, each $\Psi$ defines a section $\Phi(\Psi)$ in $\Sec(\tau_G)$ and also a nondegenerate
quadratic function $\kappa(\Psi)$.  The quadratic linking function of $\kappa$,
$\delta^c(\kappa)$, is the triple $(G,q^c(\kappa),\beta)$ where $\beta = [\alpha]$ is the image of
$\alpha$ in $G$ and $q^c(\kappa): \Sec(\tau) \ra \QL(b(\lambda))$ is the function defined
by $q^c(\Phi(\Psi)) = q^c(\kappa(\Psi))$.  Here $\QL(b)$ denotes the set of quadratic linking functions
which refine the linking form induced by $\lambda$, $b(\lambda)$.  A crucial property of $q^c(\kappa)$
is that it is equivariant with respect to the actions of $\Hom(G^{**},TG)$ on $\Sec(\tau_G)$ and
$\QL(b(\lambda))$ which are defined, for $\phi \in \Hom(G^{**},TG)$, $\Phi \in \Sec(\tau_G)$ and $q \in
\QL(b(\lambda))$ by
$\phi\cdot\Phi = \Phi - \phi$ and $\phi \cdot q = q_{\phi(\tau_G(\beta)/2)}$.  We now make the following
definitions for an abelian group $G$ containing the even element $\beta$ and with torsion subgroup $TG$. 
Firstly, given a linking form $b$ on
$TG$, let $\QL_J(b)$ be the subset of the power set of
$\QL(b)$ consisting of sets of the form
$$S(q) := \{q_{12a} | a \in TG, q \in \QL(b) \}.$$
We let $\Hom(G^{**},TG)$ act on $\QL_J(b)$ by setting $\phi \cdot S(q) := S(q_{\phi(\tau(\beta))/2})$
where $\phi \in \Sec(\tau_G)$ and $S(q) \in \QL_J(b)$.

\begin{Definition}

\begin{enumerate}
\item{A {\em $J$-quadratic linking family} is a triple $(G,q^c_J,\beta_J)$ where $G$ is a finitely
generated abelian group, $\beta_J$ is an even element of $G \tensor \Z/24$ and $q^c_J:\Sec(\tau_G) \ra
\QL_J(b)$ is a function for some linking form $b$ on $TG$ and which satisfies the following properties.}
\begin{enumerate}
\item[(A)]{ $q^c_J$ is equivariant with respect to the actions of $\Hom(G^{**},TG)$ on 
$\Sec(\tau_G)$.}
\item[(B)]{If $q^c_J(\Phi) = S(q)$ then $\beta_J = (\Phi(\tau(\beta))+\beta(q))\tensor 1 \in G
\tensor \Z/24$.}
\end{enumerate}
\item{For any characteristic quadratic linking family, $Q = (G,q^c,\beta)$, the of $J$-quadratic linking
family of $Q$ is the triple $Q_J = (G,q_J^c,\beta_J)$, where $\beta_J = \beta \tensor 1 \in G \tensor
\Z/24$ and
$$ q^c_J(\Phi) := S(q^c(\Phi)) = \{ q^c(\Phi)_{12a}|\, a \in TG \}. $$}
\item{Let $P$ be a highly connected $7$-dimensional manifold.  The $J$-quadratic linking family of $P$,
$Q_J(P)$, is the $J$-quadratic linking family of $Q(P)$.}
\end{enumerate}
\end{Definition}

\noindent  
Recall from Section \ref{chapprelim}.\ref{secqlfam} that an isomorphism of abelian groups $\theta:G_0 \ra
G_1$ induces a bijection $\tau_{\theta}:\Sec(\tau_{G_0}) \ra \Sec(\tau_{G_1})$.  If $Q_{J0} =
(G_0,q^c_{J0},\beta_{J0})$ and $Q_{J1}=(G_1,q^c_{J1},\beta_{J1})$ are $J$-quadratic linking families
then an isomorphism
$\theta:G_0 \cong G_1$ shall be called a $J$-isometry if $(\theta \tensor Id)(\beta_{J0}) = \beta_{J1}$
and if for every $\Phi \in \Sec(\tau_{G_0})$,
$$q^c_{J1}(\tau(\theta(\Phi)) \circ \theta|_{TG_0} = q^c_{J0}(\Phi)$$
where, if $q^c_{J1}(\tau(\theta(\Phi)) = S(q)$, then $q^c_{J1}(\tau(\theta(\Phi)) \circ \theta|_{TG_0}
= S(q \circ \theta|_{TG_0})$.

\begin{Theorem} \label{thmhomotpclass}
Let $P_0$ and $P_1$ be a pair of $2$-connected $7$-manifolds with $J$-quadratic linking families
$Q_J(P_0) = (G_0,q^c_{J0},\beta_{J0})$ and $Q_J(P_1) = (G_1,q^c_{J1},\beta_{J1})$.  Then there is a
homotopy equivalence $f:P_0 \ra P_1$ such that $f_! = \theta : H^{4}(P_0) \cong H^{4}(P_1)$
if and only if $\theta$ defines a $J$-isometry of $J$-quadratic linking families from $Q_J(P_0)$ to
$Q_J(P_1)$.
\end{Theorem}

\begin{proof}
In dimension $7$ every highly connected manifold bounds a handlebody.  So for $i=0,1$, let $P_i = \del
L_i$ and let the handlebody $L_i$ have quadratic function $\kappa_i = \kappa(H_i,\lambda_i,\alpha_i)$. 
Suppose that $f: P_0 \ra P_1$ is a homotopy equivalence.  As there are no obstructions to smoothing
$7$-dimensional $PL$-manifolds, $\SPL(P_1) \cong \SO(P_1)/\Theta_7$, where $\Theta_7$ acts by connected
sum with a homotopy sphere in the domain of a homotopy equivalence.  Applying Corollary \ref{corsslp}
we conclude that for some $\epsilon \in 24 \cdot H^*_1$ there is an almost diffeomorphism $f':P_0 \ra
\del L_{\epsilon}$ such that $f = f_{\epsilon} \circ f'$, where 
$g_{\epsilon}:L_{\epsilon} := L(\kappa(H_1,\lambda_1,\alpha_1 + \epsilon)) \ra L_1$ is the homotopy
equivalence of Lemma
\ref{corhes} and $f_{\epsilon} =g_{\epsilon}|_{L_{\epsilon}}$.  By Theorem
\ref{thmb}, $f'_!$ defines an isometry of quadratic linking families
$f'_!:Q^c(P_0) \cong \delta^c(\kappa(\lambda_1,\alpha_1 + \epsilon))$.  It is clear from the definitions
that an isometry of quadratic linking families induces an isometry of the induced
$J$-quadratic linking families.  Moreover, $f_! = f_{\epsilon !} \circ f'_! = f'_!$ since 
$f_{\epsilon !} = Id_{G_1}$.  So it remains to check that $Id_{G_1}$ defines a $J$-isometry from
$Q_J(\del(L_{\epsilon}))$ to $Q_J(\del L_1) = Q_J(P_1)$.  Clearly since $24$ divides $\epsilon$,
$\beta_1
\tensor 1 = (\beta_1 + [\epsilon]) \tensor 1 \in G_1 \tensor \Z/24$.  For any $\Phi
\in \Sec(\tau_{G_1})$, let $\Psi \in \Sec(\rho_1)$ be such that $\Phi(\Psi) = \Phi$.  Recalling that
$\epsilon(\Psi)$ denotes the component of $\epsilon$ in $H^*_1(\Psi)$ we calculate that
$$\begin{array}{cll}
q^c_J(\kappa(\lambda_1,\alpha_1 + \epsilon))(\Phi) & = S[q^c(\kappa(\lambda_1,\alpha_1 +
\epsilon))(\Psi)]\\ 
& = S[q^c(\kappa(\lambda_1(\Psi),\alpha_1(\Psi) + \epsilon(\Psi))] \\
& = S[q^c(\kappa(\lambda_1(\Psi),\alpha_1(\Psi))_{[\epsilon(\Psi)/2]}] & \text{by Lemma
\ref{lemmaqlfpert}}\\ 
& = S[q^c(\kappa(\lambda_1(\Psi),\alpha_1(\Psi))] & \text{since $12$ divides
$[\epsilon(\Psi)/2]$} \\ 
& = S[q^c(\kappa_1)(\Psi)] \\
& = S[q^c(\kappa_1)(\Phi)] \\
& = q^c_{J1}(\Phi).
\end{array}$$
Hence $Id_{G_1}$ defines a $J$-isometry from $Q_J(L_{\epsilon})$ to $Q_J(P_1)$.   

Conversely, suppose that $\theta$ is a $J$-isometry from $Q^c_J(P_0)$ to $Q^c_J(P_1)$.  Let $\Phi_0 \in
\Sec(\tau_{G_0})$ be any section and choose $\Psi_0 \in \Sec(\rho_0)$ and $\Psi_1 \in \Sec(\rho_1)$ such
that $\Phi_0 = \Phi(\Psi_0)$ and $\tau_{\theta}(\Phi_0) = \Phi(\Psi_1)$.  Then for some $a \in
TG_1$,
$$q^c(\kappa_1(\Psi_1))_{12a} \circ \theta|_{TG_0} = q^c(\kappa_1)(\tau_{\theta}(\Phi_0))_{12a} \circ
\theta|_{TG_0} = q^c(\kappa_0)(\Phi_0).$$  Applying Lemma \ref{lemmaqlfpert} in reverse, choose an
$\epsilon' \in H^*_1(\Psi_1)$ which is divisible by $24$ and so that
$q^c(\kappa(\lambda_1(\Psi_1),\alpha_1(\Psi_1)+\epsilon')) =
q^c(\kappa_1(\Psi_1))_{12a}$.
It follows that 
$$q^c(\kappa(\lambda_1(\Psi_1),\alpha_1(\Psi_1)+\epsilon')) \circ \theta|_{TG_0}
= q^c(\kappa_0)(\Phi_0)$$ 
and in particular that $\theta(\beta(q^c(\Phi_0)) = [\alpha(\Psi_1) + \epsilon']$.
Now $\beta_0 = \Phi_0(\tau_{G_0}(\beta_0)) + \beta(q^c(\Phi_0))$ and so 
$\theta(\beta_0) - [\alpha(\Psi_1) + \epsilon'] = \theta(\Phi_0(\tau_{G_0})(\beta_0)) =
\tau_{\theta}(\Phi_0)(\beta_1).$  Now $\alpha_1 - \alpha_1(\Psi_1) \in F_1^*(\Psi_1)$ and so
we may choose $\epsilon^{''} \in F_1^*(\Psi_1)$ such that 
$\theta(\beta_0) = [\alpha_1 + \epsilon' + \epsilon^{''}]$ and $\epsilon^{''}$ is necessarily 
divisible by $24$ as $(\theta(\beta_0) \tensor 1) = \beta_1 \tensor 1 \in G_1 \tensor \Z/24$. 
Setting $\epsilon = \epsilon' + \epsilon^{''} \in 24 \cdot H^*_1$, we see that
$(\alpha_1+\epsilon)(\Psi_1) = \alpha_1(\Psi_1) + \epsilon'$ and we calculate that
$$\begin{array}{cl}
q^c(\kappa(\lambda_1,\alpha_1 + \epsilon))(\tau_{\theta}(\Phi_0)) \circ \theta|_{TG_0} 
& = q^c(\kappa(\lambda_1,\alpha_1 + \epsilon))(\Psi_1) \circ \theta|_{TG_0} \\
& = q^c(\kappa(\lambda_1(\Psi_1),\alpha_1(\Psi_1)+\epsilon')) \circ \theta|_{TG_0} \\
& = q^c(\kappa_0)(\Phi_0).
\end{array}$$
By the equivariance property of characteristic quadratic linking families, $\theta$ defines an isometry
from $Q^c(P_0)$ to $\delta^c(\kappa(\lambda_1,\alpha_1 + \epsilon))$.  By Theorem \ref{thmb}, there is an
almost diffeomorphism $f':P_0 \ra \del L(\kappa(\lambda_1,\alpha_1 + \epsilon))$ such that $f'_! =
\theta$.  Composing $f'$ with $f_{\epsilon}$ we obtain a homotopy equivalence $f_{\epsilon} \circ f'$
such that $(f_{\epsilon} \circ f')_! = f_{\epsilon !} \circ f'_! = Id_{G_1} \circ \theta = \theta.$ 
\end{proof}

\section{Smooth $S^3$ bundles over $S^4$} \label{secsphere}
As we remarked in the introduction, the recent interest in $2$-connected $7$-manifolds was generated by
Grove and Ziller's geometric result that the total space of every $S^3$-bundle over $S^4$ carries a
metric of non-negative sectional curvature \cite{GZ}.  These manifolds have a lengthy history in
differential topology which is covered in \cite{CE} where the manifolds are classified up to
diffeomorphism, $PL$-homeomorphism, homeomorphism and homotopy.  The proofs in \cite{CE} used
the results of Wilkens \cite{Wi1,Wi2} and calculations of the invariants $s_1$ and $\bs$.  By Hatcher's
result on the homotopy type of $\rm{Diff}(S^3)$ \cite{Ha} the classification holds equally well for
smooth bundles over $S^4$ with fiber a differentiable $3$-sphere.  We now present the classification of
these total spaces using the results and notation of this thesis.

Recall from Section \ref{secsmoothings} the oriented $\R^4$-bundles over $S^4$, $\xi(m,n)$.   Let
$L_{m,n}$ be the handle body which is the total space of the disc bundle of $\xi(m,n)$ and let $P_{m,n}$
be the boundary of $L_{m,n}$.  We orient $L_{m,n}$ and $P_{m,n}$ by fixing an orientation on $S^4$ and
using the orientation of the vector bundle on the fiber.\footnote{We warn the reader that this
orientation convention differs from the one use in \cite{CE}.  There, following James and Whitehead, we
fixed an orientation on $S^4$ and oriented the fiber so that the Euler number of $\xi(m,n)$ was
positive.  The two conventions agree for $n \geq 0$.  The later one seems more natural.}
There is a diagram of commuting fiber bundles
\[
\divide\dgARROWLENGTH by2
   \begin{diagram}
\node{D^4} \arrow{e}  \node{L_{m,n}}  \arrow{e,t}{\pi_L} \node{S^4}  \\
\node{S^3} \arrow{n} \arrow{e} \node{P_{m,n}} \arrow{n,r}{i} \arrow{e,t}{\pi_P} \node{S^4.}
\arrow{n,r}{Id}
\end{diagram}
\]
\noindent
Let $s_4$ and $s^4$ denote the preferred generators of $H_4(S^4)$ and $H^4(S^4)$ respectively.  The group
$H_4(L_{m,n})$ is infinite cyclic with a preferred generator $s$ such that $\pi_{L*}(s) = s_4$.  The
bundle $\xi(m,n)$ has Euler number $n$ and stable invariant $n+2m$, thus $s \cdot s = n$ and the
tangential invariant of $L_{m,n}$, $\alpha_{m,n}$, satisfies $\alpha_{m,n}(s)=(n+2m)$.\footnote{The
cohomology class $\alpha_{m,n} \in H^4(L_{m,n})$ is also equal to the Spin characteristic class,
$\halfp(L_{m,n})$.}   The group $H^4(P_{m,n})$ is isomorphic to $\Z/n$ with preferred generator $e
=\pi^{*}_P(s^4)$.  By applying Proposition \ref{propfree} and Remark
\ref{rempltophb} it is easy to see that $L_{m_0,0}$ (resp. $P_{m_0,0}$) is diffeomorphic to $L_{m_1,0}$
(resp. $P_{m_1,0}$) if and only if $m_0 = \pm m_1$.  We assume now that $n \neq 0$.   The
algebraic invariants of $L_{m,n}$ and $P_{m,n}$ are as follows (where $j$ and $k$ are integers, we
have used Proposition \ref{props_1P} to compute $s_1(P_{m,n})$ and $\bs(P_{m,n})$ and $\epsilon(n) =
n/|n|$):
$$\begin{array}{cl}
\alpha_{m,n}: \Z \ra \Z & j\cdot s = (n+2m) j, \\
\lambda_{m,n} = \lambda_n:  \Z \times \Z  \lra  \Z & (j \cdot s,k \cdot s)  \mapsto  njk, \\
\kappa_{m,n}: \Z  \lra  \Z &
 j \cdot s \mapsto   n j^2 + (n+2m)  j, \\
\beta_{m,n} \in \Z/n & \beta_{m,n} = 2m \cdot e\\
b_{m,n}= b_n: \Z/n \times \Z/n \lra \Q/\Z & (j \cdot e,k \cdot e) \mapsto \frac{1}{n}(jk) ~ \modu ~\Z\\
q_{m,n}:  \Z/n  \lra  \Q/\Z & 
 j e \mapsto \frac{1}{2n}(j^2 + 2m  j) ~ \modu~ \Z\\
s_1(P_{m,n}) = \frac{1}{28} \cdot \frac{1}{8n}({(n+2m)^2} - n\epsilon(n)) & 
s_1(P_{m,n}) = \frac{1}{8n}({(n+2m)^2} - n\epsilon(n)).
\end{array}$$

\noindent
When $n \neq 0$, $P_{m,n}$ is a rational homology sphere which falls under the classification 
Theorem \ref{thma}.  There ({\em 1} and {\em 5}) it is stated that the isometry class of $q_{m,n}$
determines the topological classification and ({\em 2}) that the isometry class of $q_{m,n}$ along with
$s_1(P_{m,n})$ determines the smooth classification.  Moreover by ({\em 4}), the quadratic function
$q_{m,n}$ is determined by the pair $(b(n),\beta_{m,n})$ and $\bs(P_{m,n}) \in \Q/\Z$.  
Lastly, by Theorem \ref{thma} ({\em 6}), $P_{m_0,n_0}$ is homotopy equivalent to $P_{m_1,n_1}$ if and
only if $q_{m_0,n_0} \cong (q_{m_1,n_1})_{12je}$ for some integer $j$. 
It is then a matter of definitions and some simple algebra that the following classification holds.

\begin{Theorem}[\cite {CE}]
Let $n \neq 0$ and $\epsilon = \pm 1$.  The total spaces $P_{m_0,n}$ and $P_{m_1,\epsilon n}$ are 
(1.) diffeomorphic, (2.) homeomorphic or (3.) homotopic if and only if the
following (appropriately numbered) condition holds:
\begin{enumerate}
\item{(a) $4m_0(n+m_0) +(n^2-n) \equiv \pm[4m_1(n_1+\epsilon m_1) +\epsilon(n^2-n)]~ (\modu ~ 224
n)$,\\ (b)~$2m_0 \equiv 2\alpha m_1 ~ (\modu ~ n)$ for some $\alpha \in \Z/n$ such that 
       $\alpha^2 \equiv \pm\epsilon ~ (\modu ~ n).$}
\item{(a) $4m_0(n+m_0) +(n^2-n) \equiv \pm[4m_1(n_1+\epsilon m_1)+ \epsilon(n^2-n)]~ (\modu ~ 8 ~ n)$,\\
(b) $2m_0 \equiv 2\alpha m_1 ~ (\modu ~ n)$ for some $\alpha \in \Z/n$ such that 
       $\alpha^2 \equiv \pm \epsilon ~ (\modu ~ n).$}
\item{(a) $m_0 \equiv \alpha m_1 ~ (\modu(n,12))$ where $\alpha^2 \equiv \pm \epsilon ~ (\modu(n, 12))$.}
\end{enumerate}
The diffeomorphism/homeomorphism/homotopy equivalence is orientation preserving or reversing as 
the relevant equations hold with either a $+$ sign or $-$ sign.
\end{Theorem}

\noindent
There are more detailed arithmetical calculations in \cite{CE} concerning
how many of the smooth structures within a given topological type are realized by $S^3$-bundles over
$S^4$ and of the numbers
$\alpha \in \Z/n$ such that $\alpha^2=\pm 1 ~ (\modu ~n)$.  We conclude by noting that the manifolds
$P_{1,8}$ and $P_{5,8}$ are tangentially homotopy equivalent but not homeomorphic as we now show.  One
computes that $\bs(P_{1,8}) = 9/16 \neq 1/16 =\bs(P_{5,8})$ and hence $P_{1,8}$ and $P_{5,8}$ are not
homeomorphic by Theorem \ref{thma} {\em 3}.  The quadratic functions $P_{1,8}$ and $P_{5,8}$ are
respectively $q^3(1)_{5e}$ and $q^3(1)_{e}$.  Since, $q^3(1)_{5e+12e} = q^3(1)_{e}$, it follows by
Theorem \ref{thma} {\em 6} that $P_{1,8}$ and $P_{5,8}$ are homotopy equivalent via a homotopy
equivalence $f$ which induces the identity on $H^4(P_{1,8}) = \Z/8 = H^{4}(P_{5,8})$.  Since 
$\beta(P_{1,8}) = \beta(P_{5,8}) = 2 \in \Z/8$, it follows that
$f$ is tangential.\footnote{Note that the same conclusion holds for every 
pair of manifolds $(P^7,*P^7)$ for which $(b(P),\beta(P)) \cong (b(*P),\beta(*P))$ but $\bs(P) \neq
\bs(*P)$.  There is such a pair for every finite abelian group with order divisible by $8$.} 
This pair of manifolds yields a counter example to Theorems C and Theorem 5.10 of
\cite{MTW} where it is claimed, amongst other things, that all tangentially equivalent
$2$-connected
$7$-manifolds are homeomorphic.

\end{document}